\documentclass[11pt]{article}

\usepackage{amsmath}
\usepackage{amssymb}
\usepackage{amsthm}
\usepackage{color}
\usepackage{graphicx}
\usepackage{color}

\usepackage{comment}

\usepackage{fancyhdr}
\usepackage{lastpage}
\def\be{\begin{eqnarray}}
\def\ee{\end{eqnarray}}

\def\b*{\begin{eqnarray*}}
\def\e*{\end{eqnarray*}}

\newtheorem{Theorem}{Theorem}[section]
\newtheorem{Remark}{Remark}[section]
\newtheorem{Definition}[Remark]{Definition}
\newtheorem{Proposition}[Theorem]{Proposition}

\newtheorem{Assumption}[Remark]{Assumption}

\newtheorem{Lemma}[Theorem]{Lemma}
\newtheorem{Corollary}[Theorem]{Corollary}
\newtheorem{Example}[Remark]{Example}

\makeatletter \@addtoreset{equation}{section}

\newcommand*{\TitleFont}{
      \fontsize{19}{25}%
      \selectfont}

\usepackage{color}

\def\bit{\begin{itemize}}
\def\eit{\end{itemize}}

\def \D{\mathbb{D}}
\def \E{\mathbb{E}}
\def \F{\mathbb{F}}
\def \G{\mathbb{G}}
\def \P{\mathbb{P}}
\def \Ph{\widehat \mathbb{P}}
\def \Q{\mathbb{Q}}
\def \R{\mathbb{R}}
\def \S{\mathbb{S}}
\def \T{\mathbb{T}}

\def\Ac{{\cal A}}
\def\Bc{{\cal B}}

\def\Dc{{\cal D}}

\def\Fc{{\cal F}}
\def\Gc{{\cal G}}

\def\Lc{{\cal L}}

\def\Pc{{\cal P}}
\def\Pcb{\overline{\Pc}}
\def\Pch{\widehat{\Pc}}

\def\Tc{{\cal T}}
\def\Uc{{\cal U}}

\def\Pb{\overline{\P}}
\def\Qb{\overline{\Q}}
\def\Rb{\overline{\R}}

\def\Dcb{\overline{\Dc}}

\def\xb{\mathbf{x}}
\def\Ab{\overline{A}}
\def\Mb{\overline{\M}}
\def\Fh{\widehat{\F}}
\def\Fch{\widehat{\Fc}}
\def\Ft{\widetilde{\F}}

\def \Ph{\widehat{\P}}

\def \taub{\bar{\tau}}

\def \M{\mathbb{M}}
\def \Mh{\widehat{M}}

\def \1{{\bf l}}

\def \O{\Omega}
\def \Om{\Omega}
\def \Ob{\overline{\Omega}}
\def \Omb{\overline{\Om}}
\def \Omh{\widehat{\Om}}
\def \Omt{\widetilde{\Om}}

\def \o{\omega}
\def \om{\o}

\def \omb{ \bar{\om}}
\def \omt{\tilde{\om}}

\def \eps{\varepsilon}

\def \Fcb{\overline{{\cal F}}}
\def \Fbb{\overline{\mathbb{F}}}

\def \w{\mathsf{w}}

\def \wb{\bar{\w}}
\def \mut{\tilde{\mu}}
\def \sigmat{\tilde{\sigma}}

\def \Pct{\widetilde{\Pc}}
\def \Pt{\widetilde{\P}}
\def \Qt{\widetilde{\Q}}

\def \Vt{\widetilde{V}}

\def \C{\mathbb{C}}
\def \D{\mathbb{D}}
\def \0{\mathbf{0}}
\def \Gh{\widehat{\G}}

\def \Fct{\widetilde{\Fc}}

\def \bb {\mathbf{b}}
\newcommand{\rmi}{{\rm (i)$\>\>$}}
\newcommand{\rmii}{{\rm (ii)$\>\>$}}
\newcommand{\rmiii}{{\rm (iii)$\>\,    \,$}}
\newcommand{\rmiv}{{\rm (iv)$\>\>\,$}}
\newcommand{\rmv}{{\rm (v)$\>\>\,$}}

\newcommand{\rma}{{\rm a)$\>\>$}}
\newcommand{\rmb}{{\rm b)$\>\>$}}

\def\x{\times}
\def\ox{\otimes}

\def \proof{{\noindent \sc Proof. }}
\def\edoc{\end{document}}

\title{\TitleFont Capacities, Measurable Selection \& Dynamic Programming\\
		Part II: Application in Stochastic Control Problems 
		}

\author{ El Karoui Nicole
\footnote{Laboratoire de Probabilit\'es, Statistique et Mod\'elisation, Sorbonne-Univerisit\'e Paris, France.}
\\
\and Tan Xiaolu
\footnote{Department of Mathematics, The Chinese University of Hong Kong, Hong Kong SAR.}
}

\date{\today}

\begin{document}

\maketitle

\abstract{
	We provide an overview on how to use the measurable selection techniques to derive the dynamic programming principle for a general stochastic optimal control/stopping problem.
	By considering its martingale problem formulation on the canonical space of paths, one can check the required measurability conditions. 
	This covers in particular the most classical controlled/stopped diffusion processes problems.
	Further, we study the approximation property of the optimal control problems by piecewise constant control problems.
	As a byproduct, we obtain an equivalence result of the strong, weak and relaxed formulations of the controlled/stopped diffusion processes problem.
}

\vspace{0.5em}


\vspace{0.5em}

\vspace{3mm}


\section{Introduction and examples}

\subsection{Introduction}

	The theory of stochastic control has been largely developed since 1970s, and plays an important role in engineering, physics, economics and finance, etc.
	In particular, with the development of financial mathematics since 1990s, it becomes an important subject and a powerful tool in many applications.
	A general optimal control/stopping problem can be described as follows:
	``The time evolution of some stochastic process is affected by `action' taken by the controller. 
	The action taken at every time depends on the information available to the controller. 
	The control objective is to choose actions as well as a time horizon that maximize some quantity, for example the expectation of some functional of the controlled/stopped sample path ...'' (Fleming (1986, \cite{Fleming_1986}).

	\vspace{0.5em}

	In the stochastic control theory, the controlled diffusion processes problem seems to be the most popular and most studied subject, especially motivated by its applications in finance.
	In particular, due to different motivations and applications, different (strong, weak or relaxed) formulations have been introduced, as in the theory of stochastic differential equations (SDEs).
	In the control theory, much effort has been devoted to establish rigorously the dynamic programming principle (DPP).
	The DPP consists in splitting a global time optimization problem into a series of local time optimization problems in a recursive manner, 
	and it has a very intuitive meaning, that is, a globally optimal control is also locally optimal at any time.
	This can also be seen as an extension of the tower property of Markov process in the optimization context.
	As applications, it allows one to characterize the optimal controlled/stopped process, to obtain a viscosity solution characterization of the value function, to derive the numerical algorithms, etc.

	\vspace{0.5em}
	
	The main objective of the paper is first to give a global study to the DPP of the continuous time stochastic control/stopping problems,
	and then to study its approximation by piecewise constant control problems.
	In particular, we obtain the DPP for different formulations of the controlled/stopped diffusion processes problem as well as their stability and equivalence.

	\vspace{0.5em}

	For the discrete time stochastic control problems, the DPP has been well studied by many authors, see e.g. Bertsekas and Shreve (1978, \cite{Bertsekas_1978}), or Dellacherie (1985, \cite{Dellacherie_MaisonJeux}), etc.
	However, the continuous time case becomes much more technical.
	One of the main difficulties is to show the measurability of the set of controls on the space of continuous time paths.
	To overcome this difficulty, a classical approach is to impose continuity or semi-continuity conditions on the value function of the control problem, or to consider its semi-continuous envelope, 
	and then to utilize the separability property of the time-state space
	(see e.g. Fleming and Rishel (1975, \cite{FlemingRishel}), Krylov (1980, \cite{Krylov}), Fleming and Soner (1993, \cite{FlemingSoner}), Touzi (2012, \cite{Touzi}), Bouchard and Touzi (2011, \cite{BouchardTouzi}), etc.).
 	In the 1980s,  many authors (e.g. El Karoui (1981, \cite{ElKarouiSF}), El Karoui and Jeanblanc (1988, \cite{ElKarouiJeanblanc}), etc.) studied controlled/stopped Markov processes problem where only the drift part is controlled, using measure change techniques with Girsanov theorem.
	The existence of reference probability measure simplifies the questions on the null sets, 
	and allows one to model, in a very general setting, the action of the controller through a family of martingale likelihood processes.
	At the same time, another approach is to consider the martingale problem formulation of the control problem, see e.g. Haussmann (1985, \cite{Haussmann_1985}), Lepeltier and Marchal  (1977, \cite{Lepeltier_1977}), El Karoui, Huu Nguyen and Jeanblanc (1987, \cite{ElKaroui_1987}), etc.
	In \cite{ElKaroui_1987} (see in particular Theorems 6.2, 6.3 and 6.4), the authors considered a (possibly degenerate) controlled diffusion (or diffusion-jump) processes problem, 
	where they interpreted the control processes as Young measures,
	and then derived the DPP by using measurable selection techniques without any regularity conditions. 
	Using similar ideas, but in a non-Markovian context and with a more modern presentation, Nutz and van Handel (2013, \cite{NutzHandel}), Neufeld and Nutz (2013, \cite{NeufeldNutz}) and Zitkovic (2014, \cite{Zitkovic}) provided the DPP for a class of control problems by considering their law on the canonical space of paths.
	Following these works, we formulated an abstract framework to derive the DPP for a general stochastic control/stopping problem in our accompanying paper \cite{ElKarouiTan1}.
 	Let us also notice that by the so-called stochastic Perron's method, one can obtain the viscosity solution characterization of a stochastic control problem without using DPP, and then deduce DPP posteriorly, see e.g. Bayraktar and Sirbu (2013, \cite{BayraktarSirbu}), etc.

 	\vspace{0.5em}

	In our accompanying paper \cite{ElKarouiTan1}, we have revisited the way how to deduce the measurable selection theorem by the capacity theory, where one of the basic ideas is to extend properties on the compact sets of a metric space to the Borel measurable sets by approximations.
	In the context of stochastic control/stopping problems, we are interested in its approximation by piecewise constant controls, which can be considered as a stability problem.
	A piecewise constant control process is in fact a sequence of adapted random variables along some (deterministic or stochastic) time instants,
	which is a natural extension of the discrete-time control,
	and is also closely related to the stochastic impulse control (or switching) problems (see e.g. Lepeltier and Marchal \cite{LepeltierMarchal}, Bismut \cite{Bismut}, 
	etc.).
	The idea to approximate a continuous time model by piecewise constant models has been largely used by Krylov (1980, \cite{Krylov}).
	And it is very similar to Donsker's theorem where the discrete time random walk converges weakly to a continuous time process, and also to Kushner and Dupuis's (1992, \cite{KushnerDupuis}) idea to approximate the continuous time control problem by discrete time controlled Markov chains in their numerical methods.

	\vspace{0.5em}

	Restricted to the  controlled diffusion processes problem with piecewise constant controls, it is then easy to prove the equivalence of the strong and weak formulations (see e.g. Dolinsky, Nutz and Soner (2012, \cite{DNS})),
	then a by-product of this stability result is the equivalence of different formulations of the continuous time control problems.
	 We also notice that such an equivalence is well-known for the optimal stopping problems under the so-called K-property
	(see e.g. Szpirglas and Mazziotto (1977, \cite{SzpirglasMazziotto}), and El Karoui, Lepeltier and Millet (1992, \cite{ElKaroui_Lepeltier}).

	\vspace{0.5em}

	The rest of the paper is organized as follows.
	In Section \ref{subsec:diffusion}, we provide a first discussion on the class of controlled/stopped diffusion processes problems, as examples, since it consists of a class of the most interesting and studied problems.
	Next, in Section \ref{sec:DPP}, we give an overview on how to deduce the DPP of a general stochastic control/stopping problem using measurable selection techniques under some measurability and stability conditions.
	Then in Section \ref{sec:MartPb}, we study a general controlled/stopped martingale problem 
	and show how to check the measurability and the stability conditions to obtain the DPP.
	Under this framework, we obtain easily the DPP for different formulations of the controlled/stopped diffusion processes problems.
	Finally, we study the stability of the control/stopping problem in Section \ref{sec:stability}.
	As a by-product, we obtain the equivalence of different formulations of the controlled/stopped diffusion processes problem.

\vspace{0.5em}

\noindent {\bf Notations.} \rmi Let $d \ge 1$ be an integer, we denote by $\S^d$ the collection of all $d \x d$ dimensional matrices, and define $\Rb_+ := \R_+ \cup \{ +\infty \} = [0, \infty]$
and $\Rb := \R \cup \{ - \infty, +\infty\} = [-\infty, \infty]$.
For $c, c' \in \R^d$ and $A, A' \in \S^d$, we denote the scalar product by $c \cdot c' := \sum_{i=1}^d c_i c_i'$ and $A : A' := \mbox{Tr}(A (A')^T)$, the corresponding norm are then denoted by $|c|$ and $\| A\|$.

\vspace{0.5em}

\noindent \rmii Let $E$ and $U$ be two (non-empty) Polish spaces, we denote by $\Om = \D(\R_+, E)$ the space of all c\`adl\`ag $E$-valued paths on $\R_+$, 
and by $\F = (\Fc_t)_{t \ge 0}$ the canonical filtration generated by the canonical process $X$.
We also introduce an enlarged canonical space by $\Omh := \Rb_+ \x \Om$ and $\Omb := \Rb_+ \x \Om \x \M$, where $\M$ denotes the collection of all $\sigma$-finite measures on $\R_+ \x U$ whose marginal distribution on $\R_+$ coincides with the Lebesgue measure.
Given $\sigma$-finite measure $m \in \M$, it follows by disintegration/conditioning that one has the representation $m(du, dt) = m_t(du) dt$ with $m_t \in \Pc(U)$ for all $t \in \R_+$, where $\Pc(U)$ denotes the collection of all (Borel) probability measure on $U$.

\vspace{0.5em}

\noindent \rmiii When studying controlled diffusion processes problem, we fix $E \equiv \R^d$ so that $\Om = \D(\R_+, \R^d)$. 
In this context, we denote by $B$ the canonical process, and by $\P_0$ the Wiener measure under which $B$ is a standard Brownian motion, and $\F^a$ the associated augmented filtration.
In this context, we also consider the enlarged canonical space $\Omt := \Om_0 \x \Rb_+ \x \Om \x \M$, with $\Om_0 := \Om$.

\vspace{0.5em}

\noindent \rmiv In some cases, we also consider an abstract filtered probability space, denoted by $(\Om^*, \Fc^*, \P^*, \F^*=(\Fc_t^*)_{t \ge 0})$.

\vspace{0.5em}

\noindent \rmv  For a random variable $\xi$ taking value in $\Rb$, let define its expectation by $\E[ \xi] := \E[ \xi^+] - \E[ \xi^-]$, with the convention that $\infty - \infty = -\infty$ to avoid the integrability problem.

\subsection{Examples: controlled/stopped diffusion processes problems}
\label{subsec:diffusion}

	In the optimal control/stopping theory, most of the literature has been focused on the diffusion processes case due to its complexity and its importance in applications, see e.g. Krylov \cite{Krylov}, Fleming and Soner \cite{FlemingSoner}, Borkar \cite{Borkar2}, Yong and Zhou \cite{YongZhou}, Pham \cite{Pham}, Touzi \cite{Touzi},
	El Karoui et al. \cite{ElKaroui_1987} and also the survey paper of Borkar \cite{Borkar}, etc.

	\vspace{0.5em}
	
	For the controlled/stopped diffusion processes problems, different formulations have been studied in the literature.
	Let us stay in a general path-dependent setting and recall these formulations.
	Let $\Om = \D(\R_+, \R^d)$ denote the canonical space of c\`adl\`ag paths on $\R_+$, $U$ be a (non-empty) Polish space,
	we shall consider the controlled diffusion processes with (Borel) measurable coefficient functions $(\mu, \sigma) : \R_+ \x \Om \x U \longrightarrow \R^d \x \S^d$,
	as well as reward functions $L: \R_+ \x \Om \x U \longrightarrow \Rb$ and $\Phi:  \Rb_+ \x \Om \longrightarrow \Rb$.
	To avoid possible integrability problems, we also assume that, for all $\om \in \Om$ and $T \ge 0$,
	\begin{equation} \label{eq:mu_sigma_integ}
		\int_0^T \sup_{u \in U} 
		\Big( | \mu(t, \om, u) | + \| \sigma(t, \om, u) \|^2 \Big) dt < \infty. 
	\end{equation}
	The above technical integrability condition can nevertheless be relaxed (see e.g. Section \ref{subsubsec:unbound_mu_sigma}).

\paragraph{A  strong formulation of the optimal control/stopping problem}

	Let $(\Om^*, \Fc^*, \P^*)$ be a probability space equipped with a $d$-dimensional standard Brownian motion $B$,
	let $\F^* = (\Fc^*_t)_{t \ge 0}$ be the augmented Brownian filtration generated by  $B$ (with completion), and $\Tc$ denote the collection of all $\F^*$-stopping times.
	We denote by $\Uc$ the collection of all $U$-valued $\F^*$-predictable processes.

	\vspace{0.5em}

	Given the initial condition $x_0 \in \R^d$ and the control process $\nu \in \Uc$, the controlled process $X^{\nu}$ is defined as the strong solution to the controlled stochastic differential equation (SDE):
	\be \label{eq:controlledSDE}
		X^{\nu}_t 
		~=~ 
		x_0 + \int_0^t \mu \big(s, X^{\nu}_{s \wedge \cdot}, \nu_s \big) ds 
		+ \int_0^t \sigma\big(s, X^{\nu}_{s \wedge \cdot}, \nu_s \big) dB_s,
		~~t \ge 0.
	\ee
	In practice, sufficient conditions (such as Assumption \ref{assu:Lip_coef}) will be assumed on $\mu$ and $\sigma$ to ensure that SDE \eqref{eq:controlledSDE} has a unique strong solution,
	which is an adapted continuous process in the fixed filtered probability space.
	Then a general optimal control/stopping problem is given by
	\be \label{eq:VS}
		V_S ~:=~ \sup_{\tau \in \Tc}~ \sup_{\nu \in \Uc} ~\E \Big[ \int_0^{\tau} L\big(t, X^{\nu}_{t \wedge \cdot}, \nu_t \big) dt  + \Phi \big(\tau, ~X^{\nu}_{\tau \wedge \cdot} \big) \Big].
	\ee
	
	\begin{Remark}
	{\rm	\rmi When $U$ is a singleton, i.e. $U = \{ u_0 \}$, the above control/stopping problem reduces to a pure optimal stopping problem.

	\vspace{0.5em}

	\noindent \rmii When the reward function satisfies $\Phi(t, \om) = -\infty$ for all $t \in [0,\infty)$, so that the optimal stopping time is clearly $\hat \tau \equiv \infty$,
	the above control/stopping problem reduces to a pure optimal control problem.
	
	\vspace{0.5em}

	\noindent \rmiii With $T > 0$, if the reward functions satisfy $\Phi(t, \om) = \Phi(T, \om_{T \wedge \cdot})$ and $L(t, \om, u) \equiv 0$ for all $(t, \om) \in (T, \infty] \x \Om$, 
	the initial infinite horizon control/stopping problem reduces to a finite horizon problem on $[0,T]$.
	}
	\end{Remark}

\paragraph{A piecewise constant control problem}

	Recall that $\Uc$ denotes the collection of all $U$-valued $\F^*$-predictable processes.
	A more elementary problem is to consider the piecewise constant control, i.e. the control process $\nu$ stays constant over some (deterministic or stochastic) intervals.
	From a practical point of view, it seems more natural and important in applications;
	and it is also closely related to the stochastic impulse control/switching problems (but with a null switching cost).
	More precisely, a piecewise constant mixed control-stopping problem is given by
	\be \label{eq:piecewise_ctrl}
		V_{S_0} ~:=~ \sup_{\tau \in \Tc}~ \sup_{\nu \in \Uc_0} 
		~\E \Big[ \int_0^{\tau} L\big(t, X^{\nu}_{t \wedge \cdot}, \nu_t \big) dt  + \Phi \big(\tau, ~X^{\nu}_{\tau \wedge \cdot} \big) \Big],
	\ee
	where $\Uc_0$ is is the set of all $\nu \in \Uc$ such that $\nu_t := \sum_{k \ge 0} \hat \nu_k \1_{(\tau_k, \tau_{k+1}]}(t)$ with a sequence of finite stopping times $0 = \tau_0 < \tau_1 < \cdots < \tau_k < \cdots$.	
	
	\vspace{0.5em}
	
	One can naturally expect to approximate a general control process $\nu \in \Uc$ by a sequence of elementary controls in $\Uc_0$, which can be seen as a stability result.
	Notice that such an approximation method is also a key technique to construct weak solutions to SDEs (see e.g. Stroock and Varadhan \cite{StroockVaradhan}).

	\begin{Example}[Nisio semi-group problem] \label{exam:Nisio}
	{\rm
		The above piecewise constant control problem has been studied in a much more general formulation, named Nisio semi-group problem (see e.g. El Karoui, Lepeltier and Marchal \cite{ElKarouiLepeltierMarchal}).
	Let us consider a simplified case, where $\mu$ and $\sigma$ are Markovian and time homogeneous, 
	i.e. $(\mu, \sigma)(s, \om, u) = (\mu_0, \sigma_0)(\om_s,u)$ for some  function $(\mu_0, \sigma_0): \R^d \x U \to \R^d\x\S^d$.
	For every fixed $u \in U$, we denote by $X^{x_0, u}$ the unique strong solution of SDE \eqref{eq:controlledSDE} with initial condition $x_0$ and constant control $\nu_s \equiv u$.
	Under Lipschitz conditions on the coefficients, it is easy to deduce that  $X^{x_0,u}$ is a Markov process, we denote by $(P^u_{\tau})_{\tau \in \Tc}$ the corresponding (transition) semi-group defined by:
	$$
		P^u_{\tau} f(x_0) ~:=~ \E \big[ f(X^{x_0,u}_{\tau}) \big].
	$$
	We next define a simple optimal stopping problem, together with constant control, by
	$$
		R \phi(x) ~:=~ \sup\big\{ P^u_{\tau} \phi(x) ~: u \in U, \tau \in \Tc \big\}. 
	$$
	It is then shown in \cite{ElKarouiLepeltierMarchal} that the operator $R$ maps a positive upper semi-analytic function to a positive upper semi-analytic function (see Section \ref{sec:DPP} for a precise definition of upper semi-analytic functions).
	In this context, one can further show that, the optimal control/stopping problem defined in \eqref{eq:piecewise_ctrl} is equivalent to $R^{\infty} \Phi := \lim_{n \to \infty} R^n \Phi$, 
	which is in turn a gambling house model studied by Dellacherie \cite{Dellacherie_MaisonJeux}.
	We nevertheless insist that \cite{ElKarouiLepeltierMarchal} considers a more general framework with a class of semi-groups $(P^u_{\tau})_{\tau \in \Tc}$.
	}
	\end{Example}

\paragraph{A weak formulation of the optimal control/stopping problem}

	In the strong formulation \eqref{eq:VS}, the solution of the controlled SDE \eqref{eq:controlledSDE} is given in a fixed probability space, equipped with a fixed Brownian motion.
	When the probability space (and the associated Brownian motion) is no longer fixed, one obtains a weak formulation of the optimal control/stopping problem.

	\begin{Definition} \label{def:weak_control}
		A term $\alpha = (\Om^{\alpha}, \Fc^{\alpha}, \P^{\alpha}, \F^{\alpha} = (\Fc^{\alpha}_t)_{t \ge 0}, \tau^{\alpha}, X^{\alpha}, B^{\alpha}, \nu^{\alpha})$ is called a weak control with initial condition $x_0 \in \R^d$,
		if $(\Om^{\alpha}, \Fc^{\alpha}, \P^{\alpha}, \F^{\alpha})$ is a filtered probability space, equipped with a stopping time $\tau^{\alpha}$, a $d$-dimensional Brownian motion $B^{\alpha}$, and a $U$-valued predictable process $\nu^{\alpha}$, together with an adapted continuous process $X^{\alpha}$ such that
		$$
			X^{\alpha}_t ~=~ x_0 + \int_0^t \mu \big(s, X^{\alpha}_{s \wedge \cdot}, \nu^{\alpha}_s \big) ds 
			+ \int_0^t \sigma\big(s, X^{\alpha}_{s \wedge \cdot}, \nu^{\alpha}_s \big) dB^{\alpha}_s,
			~~t \ge 0,
			~\mbox{a.s.}
		$$
	\end{Definition}
	Notice that the stochastic integral term in the above definition is implicitly assumed to be well defined.
	Let us denote by $\Ac_W$ the collection of all weak control with fixed initial condition $x_0$, then a weak formulation of the optimal control/stopping problem is given by
	$$
		V_W 
		~:=~
		\sup_{\alpha \in \Ac_W} \E \Big[\int_0^{\tau^{\alpha}} L\big(t, X^{\alpha}_{t \wedge \cdot}, \nu^{\alpha}_t \big) dt  + \Phi \big(\tau^{\alpha}, X^{\alpha}_{\tau^{\alpha} \wedge \cdot} \big)  \Big]. 
	$$

\paragraph{A relaxed formulation of the optimal control/stopping problem}

	The relaxed formulation of the controlled diffusion processes problem has been introduced by Fleming \cite{Fleming_1978}, El Karoui, Huu Nguyen and Jeanblanc \cite{ElKaroui_1987},
	where the main idea is to relax the $U$-valued control process $\nu^{\alpha}$ to be a $\Pc(U)$-valued process, with $\Pc(U)$ denoting the space of all (Borel) probability measures on $U$.
	Namely, the controller takes no longer a fixed action in the space $U$, but a randomized action of different elements in $U$ following some distribution.
	The Brownian motion will also be replaced by a continuous martingale measure in the corresponding SDE.

	\begin{Definition} \label{def:relaxed_control}
		\rmi Let $(\Om, \Fc, \F, \P)$ be a filtered probability space satisfying the usual condition, $(M_t)_{t \ge 0}$ be a $\Pc(U)$-valued predictable process, and $\Bc(U)$ denote the Borel $\sigma$-field of $U$. 
		Then $(\Mh_t(du))_{t \ge 0}$ is called a continuous martingale measure with intensity $(M_t)_{t \ge 0}$ if
		\begin{itemize}
			\item $(\Mh_t(A))_{t \ge 0}$ is continuous martingale with $\Mh_0(A) = 0$, for all $A \in \Bc(U)$;
			
			\item $(\Mh_t(A))_{t \ge 0}$ and $(\Mh_t(B))_{t \ge 0}$ are orthogonal whenever $A, B \in \Bc(U)$ satisfy $A \cap B = \emptyset$;
			
			\item the quadratic variation processes satisfy $\langle \Mh(A) \rangle_t = \int_0^t M_s(A) ds$ for all  $t \ge 0$ and $A \in \Bc(U)$.
		\end{itemize}
				
		
		\noindent \rmii A term $\alpha = (\Om^{\alpha}, \Fc^{\alpha}, \P^{\alpha}, \F^{\alpha} = (\Fc^{\alpha}_t)_{t \ge 0}, \tau^{\alpha}, X^{\alpha}, M^{\alpha}, \Mh^{\alpha})$ is called a relaxed control with initial condition $x_0 \in \R^d$,
		if $(\Om^{\alpha}, \Fc^{\alpha}, \P^{\alpha}, \F^{\alpha})$ is a filtered probability space, equipped with a stopping time $\tau^{\alpha}$, a $\Pc(U)$-valued predictable process $M^{\alpha}$, 
		and a continuous martingale measure $\Mh^{\alpha}$ with intensity $M^{\alpha}$, together with an adapted continuous process $X^{\alpha}$ such that
		$$
			X^{\alpha}_t = x_0 + \int_0^t \!\! \int_U \mu \big(s, X^{\alpha}_{s \wedge \cdot}, u \big) M^{\alpha}_s(du) ds 
			+ \int_0^t \!\! \int_U \sigma\big(s, X^{\alpha}_{s \wedge \cdot}, u \big) \Mh^{\alpha}(du, ds),
			~~t \ge 0,
			~\mbox{a.s.}
		$$
	\end{Definition}
	
	The martingale measure has been initially introduced in a very general setting (with more general intensity measure), we nevertheless only recall its definition in a setting enough for our uses.
	For the stochastic integration w.r.t. the  martingale measure, as well as their basic properties, let us refer to El Karoui and M\'el\'eard \cite{ElKaroui_Meleard} and the references therein.
	Let us denote by $\Ac_R$ the collection of all relaxed control with fixed initial condition $x_0$, we then obtain the following relaxed formulation of the optimal control/stopping problem:
	$$
		V_R
		~:=~
		\sup_{\alpha \in \Ac_R} \E \Big[\int_0^{\tau^{\alpha}} \!\!\!\! \int_U L\big(t, X^{\alpha}_{t \wedge \cdot}, u \big) M^{\alpha}_t(du) dt  + \Phi \big(\tau^{\alpha}, X^{\alpha}_{\tau^{\alpha} \wedge \cdot} \big)  \Big]. 
	$$
	Notice that a weak control $\alpha \in \Ac_R$ can be considered as a relaxed control by setting $M^{\alpha}_s(du) := \delta_{\nu^{\alpha}_s}(du)$ and $\Mh^{\alpha}(du,ds) := \delta_{\nu^{\alpha}_s}(du) d B^{\alpha}_s$.

\paragraph{Strong, weak and relaxed formulations on the canonical space}

	In the SDE theory, it is classical to study the weak solution by considering the distribution of the stochastic processes,
	which is a probability measure on the canonical space of paths (see e.g. Stroock and Varadhan \cite{StroockVaradhan}).
	Similarly, one can define equivalently the weak and relaxed formulation of the optimal control/stopping problem on an appropriate canonical space.
	The natural candidate of the canonical space for the controlled diffusion processes is $\Om = \D(\R_+, E)$ with $E = \R^d$, and that for stopping times is $\Rb_+$.
	As for the control processes, we follow El Karoui, Huu Nguyen and Jeanblanc \cite{ElKaroui_1987} to consider a space of measure valued processes.
	Let us denote by $\Mb(\R_+ \x U)$ the collection of all $\sigma$-finite (Borel) measure on $\R_+ \x U$, and then define $\M$ as subset of all measures on $\R_+ \x U$ whose marginal distribution on $\R_+$ is the Lebesgue measure $ds$, i.e.
	\be \label{eq:def_M}
		\M ~:=~ \big\{ m \in \Mb(\R_+ \x U) ~: m(ds, du) = m(s, du) ds \big\}.
	\ee
	Notice that $m(s, du)$ is a measurable kernel of the disintegration of $m(ds, du)$ in $ds$.
	
	\begin{Remark}
	{\em
		Let us define the following topology on $\M$:
		we say $m_n \longrightarrow m_0$ in $\M$ if and only if 
		$$
			\int_0^{\infty} \int_U \phi(s, u) e^{-s} m_n(s, du)  ds
			~\longrightarrow~
			\int_0^{\infty} \int_U \phi(s, u) e^{-s} m_0 (s, du)  ds
		$$
		for every $\phi \in C_b(\R_+ \x U)$, i.e. the class of all bounded continuous functions defined on $\R_+ \x U$.
		Then $\M$ is a Polish space.
	}
	\end{Remark}

	\begin{Remark}
	{\rm The space $\M$ has been largely used in the literature of deterministic control theory, to introduce the so-called relaxed control.
	It is also called the Young measure since its marginal distribution is fixed.
	More importantly, the inherited weak convergence topology on $\M$ implies better convergence properties than the classical ones.
	We would like to refer to Young \cite{Young} and Valadier \cite{Valadier} for a presentation of Young measure as well as its applications, 
	and also to Jacod and M\'emin \cite{JacodMemin} for a more probabilistic point of view with the so-called stable convergence topology.
	}
	\end{Remark}
	
	Let us consider the canonical space $\Omb := \Rb_+ \x \Om \x \M$ with canonical element $(\Theta, X = (X_t)_{t \in \R_+}, M)$ defined by
	$$
		\Theta_{\infty}(\omb) := \theta,
		~~X_t(\omb) := \om_t,
		~~M(\omb) := m,
		~~\mbox{for all}~ \omb = (\theta, \om, m) \in \Omb.
	$$
	For each weak (resp. relaxed) control $\alpha$, let us define a weak control rule (resp. relaxed control rule) $\Pb^{\alpha}$ by
	\begin{equation} \label{eq:def_Pb_alpha}
		\Pb^{\alpha} := \P^{\alpha} \circ \big(\tau^{\alpha}, X^{\alpha}, \delta_{\nu^{\alpha}_s}(du)ds \big)^{-1}
		~~\mbox{\big(resp.}~
		\Pb^{\alpha} := \P^{\alpha} \circ \big(\tau^{\alpha}, X^{\alpha}, M_s^{\alpha}(du) ds \big)^{-1} \big),
	\end{equation}
	and then
	$$
		\Pcb_W ~:=~ \big\{ \Pb^{\alpha} ~: \alpha \in \Ac_W \big\}
		~~\mbox{(resp.}~
		\Pcb_R ~:=~ \big\{\Pb^{\alpha} ~: \alpha \in \Ac_R \big\}
		\mbox{)}.
	$$
	It follows immediately that
	$$
		V_W = \sup_{\Pb \in \Pcb_W} J(\Pb)
		~~\mbox{and}~
		V_R = \sup_{\Pb \in \Pcb_R} J(\Pb),
	$$
	with
	$$
		J(\Pb) 
		~:=~ 
		\E^{\Pb} \Big[\int_0^{\Theta_{\infty}} \!\!\!\! \int_U L\big(t, X_{t \wedge \cdot}, u \big) M(du, dt)  + \Phi \big(\Theta_{\infty}, X_{\Theta_{\infty}\wedge \cdot} \big)  \Big]. 
	$$
	For the strong formulation, one similarly has that
	$$
		V_S = \sup_{\Pb \in \Pcb_S} J(\Pb),
		~~\mbox{with}~
		\Pcb_S := \big\{ \Pb := \P^* \circ \big(\tau, X^{\nu}, \delta_{\nu_s}(du) ds \big)^{-1} ~: \tau \in \Tc, ~ \nu \in \Uc \big\}.
	$$
	Further, by their definition, it is clear that 
	$$
		\Pcb_S \subseteq \Pcb_W \subseteq \Pcb_R,
		~~\mbox{so that}~
		V_S \le V_W \le V_R.
	$$

\paragraph{Weak and relaxed formulations by martingale problem}

	In the classical SDE theory, the weak solution can be defined equivalently by the corresponding martingale problem on the canonical space.
	Similarly, we can define equivalently the set $\Pcb_W$ and $\Pcb_R$ of weak and relaxed control rules by the corresponding martingale problems.
	For this purpose, let us introduce the canonical filtration on the canonical space $\Omb := \Rb_+ \x \Om \x \M$.
	Let
	$$
		\Fcb_t ~:=~ \sigma \big\{
		X_s, ~M_s(\phi), ~\{\Theta_{\infty} \le s\}  ~: \phi \in C_b(\R_+ \x U), ~s \le t
		\big\}, 
		~~t \ge 0,
	$$
	and $\Fcb_{\infty} := \bigvee_{t \ge 0} \Fcb_t$,
	where $C_b(\R_+ \x U)$ denotes the set of all bounded continuous function on $\R_+ \x U$,
	and 
	\be \label{eq:def_M_phi}
		M_s(\phi) ~:=~ \int_0^s  \int_U \phi(r,u) e^{-r} M(du, dr).
	\ee
	Let $\Fbb =(\Fcb_t)_{t \ge 0}$ be the canonical filtration on $\Omb$.
	Notice that $\Theta_{\infty}$ is a $\Fbb$-stopping time.
	For every $\varphi \in C^2_b(\R^d)$, we introduce a $\Fbb$-adapted process $S^{\varphi} = (S^{\varphi}_t)_{t \ge 0}$ by
	$$
		S^{\varphi}_t
		~:=~
		\varphi(X_t) - \int_0^t \int_U \Lc^{s,u} \varphi (X_{s \wedge \cdot}) m(ds, du),
		~~t \ge 0,
	$$
	where $\Lc^{s, u}$ is the infinitesimal generator of the controlled diffusion process defined by
	\be
		\Lc^{s, u} \varphi(\om) 
		~:=~ 
		\mu(s, \om_{s \wedge \cdot}, u) \cdot D \varphi(\om_s)
		~+~
		\frac{1}{2}
		\sigma \sigma^T (s, \om_{s \wedge \cdot}, u) : D^2 \varphi(\om_s);
	\ee
	Further, let us denote by $B(\R_+, U)$ the set of all Borel measurable functions from $\R_+$ to $U$, 
	and introduce
	\be  \label{eq:def_M0}
		\M_0 ~:=~  \{ m \in \M ~: m(ds, du) = \delta_{\phi(s)} (du) ds ~\mbox{for some}~ \phi \in B(\R_+, U) \}.
	\ee
	Notice that $\M_0$ is a Borel subset of $\M$ (see e.g. Appendix of \cite{ElKaroui_1988}).
	We can now redefine equivalently $\Pcb_W$ and $\Pcb_R$ by the corresponding martingale problems.

	\begin{Proposition} \label{prop:weak_reformulation}
		One has
		$$
			\Pcb_R 
			=
			\Big\{ \Pb \in \Pc(\Omb)~:  
				S^{\varphi}~\mbox{is a}~ (\Pb, \Fbb) \mbox{-local martingale for all $\varphi \in C^2_b(\R^d)$,}
				~\Pb[X_0 = x_0] = 1
			\Big\},
		$$
		and
		$$
			\Pcb_W 
			~=~ 
			\big\{
				\Pb \in \Pcb_R ~: \Pb \big[ M \in \M_0 \big] = 1
			\big\}.
		$$
	\end{Proposition}
	\proof
	\rmi Let us first consider the relaxed formulation.
	First, it is easy to check that,
	for each $\alpha \in \Ac_R$, the induced probability measure $\Pb^{\alpha}$ in \eqref{eq:def_Pb_alpha} solves the corresponding martingale problem on $\Omb$, so that
	$$
		\Pcb_R 
		\subseteq
		\Big\{ \Pb \in \Pc(\Omb)~:  
			S^{\varphi}~\mbox{is a}~ (\Pb, \Fbb) \mbox{-local martingale for all}~ \varphi \in C^2_b(\R^d),
			~\Pb[X_0 = x_0] = 1
		\Big\}.
	$$
	Next, let $\Pb \in \Pc(\Omb)$ such that $S^{\varphi}$ is a $(\Pb, \Fbb)$-local martingale for all $\varphi \in C^2_b(\R^d)$.
	By El Karoui and M\'el\'eard \cite[Theorem IV-2]{ElKaroui_Meleard}, one can then construct (in a possibly enlarged space) a continuous martingale measure $\widehat M^{\Pb}$ with quadratic variation $M(du,dt)$ such that
	$$
		X_t = X_0 + \int_0^t \int_U \mu(s, X_{s \wedge \cdot}, u) M(ds, du) 
		+
		\int_0^t \int_U \sigma(s, X_{s \wedge \cdot}, u) \widehat M^{\Pb} (d s, du),
		~t \ge 0,
		~~\Pb \mbox{-a.s.}
	$$
	It follows that $(\Omb, \Fcb_{\infty}, \Pb, \Fbb, \Theta_{\infty}, X, M, \Mh^{\Pb})$ is a relaxed control in $\Ac_R$, so that $\Pb \in \Pcb_R$.
	
	\vspace{0.5em}
	
	\noindent \rmii For the weak control, one can easily check that for any $\alpha \in \Ac_W$, the induced $\Pb^{\alpha}$ belong to $\Pcb_R$ and satisfies $\Pb \big[ M \in \M_0 \big] = 1$.
	Hence $\Pcb_W \subset \big\{
				\Pb \in \Pcb_R ~: \Pb \big[ M \in \M_0 \big] = 1
			\big\}$.
	
	\vspace{0.5em}
	
	On the other hand, given $\Pb \in \Pcb_R$ such that $\Pb \big[ M \in \M_0 \big] = 1$,
	let us construct a weak control as follows.
	Notice that any Polish space is isomorphic to a Borel subset of $[0,1]$, let $\psi: U \to [0,1]$ be the bijection between $U$ and $\psi(U) \subseteq [0,1]$.
	Let
	\be \label{eq:def_nuM}
		\nu^M_t := \psi^{-1}(a_t),
		~~\mbox{where}~~ a_t := \frac{d}{dt} \int_0^t \int_U \psi(u) M(ds, du),
	\ee
	so that $\nu^{M}$ is $\Fbb$-predictable.
	Since  $\Pb \big[ M \in \M_0 \big] = 1$, one has $\Pb \big[ M(ds, du) = \delta_{\nu^M_s}(du) ds \big] = 1$.
	Moreover, by Strook and Varadhan \cite[Theorem 4.5.1]{StroockVaradhan}, one can construct (in a possibly enlarged space) a Brownian motion $W^{\Pb}$ such that
	$$
		X_t = X_0 + \int_0^t \mu(s, X_{s \wedge \cdot}, \nu^{M}_s) ds 
		+
		\int_0^t \sigma(s, X_{s \wedge \cdot}, \nu^{M}_s) d W^{\Pb}_s,
		~~t \ge 0,
		~\Pb \mbox{-a.s.}
	$$
	It follows that $(\Omb, \Fcb_{\infty}, \Pb, \Fbb, \Theta_{\infty}, X, W^{\Pb}, \nu^{M})$ is a weak control in $\Ac_W$,
	and hence $\Pb \in \Pcb_W$.
	\qed

	\vspace{0.5em}

	The strong formulation \eqref{eq:VS} can also be defined by an appropriate martingale problem, but on another enlarged canonical space.
	As we shall see later, these reformulations of the optimal control/stopping problem (in different formulations) on the canonical space will play an essential role to prove the dynamic programming principles, and to deduce the approximation as well as the equivalence results.
	
	\begin{Remark}
		Let us finally mention that, in the Markovian setting, a more relaxed formulation of the controlled diffusion processes problem is the linear programming formulation, which consists in considering the occupation measures induced by the controlled diffusion processes.
		We can refer to Stockbridge \cite{Stockbridge_1990, Stockbridge_1990b}, and also to Buckdahn, Goreac and Quincampoix \cite{BuckdahnGoreacQuincampoix} for a recent development of this formulation.
	\end{Remark}

\section{An overview on the dynamic programming principle}
\label{sec:DPP}

	Let us present an overview of our accompanying paper \cite{ElKarouiTan1}, on how to deduce the dynamic programming principle by measurable selection techniques.
	The approach is the same as in El Karoui, Huu Nguyen and Jeanblanc \cite{ElKaroui_1987},
	or Nutz and van Handel \cite{NutzHandel}, but we will present it in a more general setting.
	The main idea is to interpret the control as a probability measure on the canonical space,
	and then to use the notion of conditioning and concatenation of probability measures.

	\vspace{0.5em}

	Recall that $U$ and $E$ are both (non-empty) Polish spaces, 
	and $\Om := \D(\R_+, E)$ denotes the space of all $E$-valued c\`adl\`ag paths on $\R_+$,
	which is also a Polish space under the Skorokhod topology.
	The space $\M$ and $\M_0$ are introduced in \eqref{eq:def_M} and \eqref{eq:def_M0}, equipped with the weak convergence topology.

\paragraph{Canonical space, measurable selection theorem}

	As defined above, we use the canonical space $\Omb := \Rb_+ \x \Om \x \M$ to study a general optimal control/stopping problem, where the canonical element are defined by
	$$
		\Theta_{\infty}(\omb) := \theta,
		~~
		X(\omb) := \om,
		~~
		M(\omb) := m,
		~~\mbox{for all}~ \omb =(\theta, \om, m) \in \Omb.
	$$
	For every $t \in \R_+$ and $\omb =(\theta, \om, m) \in \Omb$, let us define
	$$
		\Theta_t(\omb) := \theta_t := \infty \1_{\{t < \theta\}} + \theta \1_{\{t \ge \theta\}},
		~~~
		X_t(\omb) := \om_t,
	$$
	and $\omb_{t \wedge \cdot} = (\theta_t, \om_{t \wedge \cdot}, m_{t \wedge \cdot})$, where $\om_{t\wedge \cdot} = (\om_{t \wedge s})_{s \ge 0}$ and
	$$
		m_{t \wedge \cdot} (ds, du) :=  \1_{\{s \in [0,t]\}} m(ds, du) +  \1_{\{s \in (t, \infty)\}} \delta_{u_0}(du) ds,
		~\mbox{for some fixed}~u_o \in U.
	$$
	For $t = \infty$, let us similarly define $\omb_{\infty \wedge \cdot} := \omb$, for all $\omb = (\theta, \om, m) \in \Omb$.
	Let $\Fbb = (\Fcb_t)_{t \ge 0}$ be the canonical filtration defined by, with $M(\phi)$ being defined in \eqref{eq:def_M_phi},
	$$
		\Fcb_t ~:=~ \sigma \big\{
			\Theta_s, X_s,  M_s(\phi)~ : \phi \in C_b(\R_+ \x U),~ s \le t
		\big\},
		~~~\mbox{for all}~ t \ge 0.
	$$
	Notice that $\Fcb_{\infty} := \bigvee_{t \ge 0} \Fcb_t = \Bc(\Omb)$ is clearly countably generated,
	and $\Theta_{\infty}$ is a $\Fbb$-stopping time.

	\vspace{0.5em}
	
	Notice also that $X$, $\Theta$ and $M(\phi)$ are all c\`adl\`ag processes, for any $\phi \in C_b(\R_+ \x U)$.
	Then a process $Z:  \Omb \x \Rb_+ \longrightarrow \R$ is $\Fbb$-progressively measurable (or equivalently $\Fbb$-optional) if and only if $Z$ is $\Fcb_{\infty} \otimes \Bc(\Rb_+)$-measurable, and satisfies $Z_t (\omb) = Z_t (\omb_{t \wedge \cdot})$ for all $t \ge 0$.
	Further let $\taub$ be a $\Fbb$-stopping time, a random variable $Y$ (defined on $\Omb$) is $\Fcb_{\taub}$-measurable if and only if there is some $\Fbb$-optional process $Z$ such that $Y = Z_{\taub}$.
	This implies that the $\sigma$-field $\Fcb_{\taub}$ is that generated by the map $\omb \in \Omb \longmapsto (\omb_{\taub(\omb) \wedge \cdot}, \taub(\omb)) \in \Omb \x  \Rb_+$, where the latter is equipped with the Borel $\sigma$-field $\Bc( \Omb \x \Rb_+)$.
	In particular, $\Fcb_{\taub}$ is countably generated, since $\Bc(\Omb \x \Rb_+)$ is.

	\vspace{0.5em}
	
	In the above framework, a control will be expressed equivalently as a probability measure on the canonical space $\Omb$, we then need to introduce the notion of conditioning as well as concatenation on $\Omb$.
	For all $(t, \wb) \in \Rb_+ \x \Omb$, let us denote
	$$
		\Dcb^t_{\wb} := \big\{ \omb \in \Omb ~: ( \Theta_t, X_t) (\omb) = (\Theta_t, X_t) (\wb) \big\},
		~~~~
		\Dcb_{t, \wb} := \big\{ \omb \in \Omb ~: \omb_{t \wedge \cdot} = \wb_{t \wedge \cdot} \big\}.
	$$
	When $t = \infty$, let
	$$
		\Dcb^{\infty}_{\wb} := \big\{ \omb \in \Omb ~: \Theta_{\infty} (\omb) =\Theta_{\infty} (\wb) \big\}
		~~\mbox{and}~~
		\Dcb_{\infty, \wb} := \{ \wb \}.
	$$
	Then, given fixed $t \in \Rb_+$ and $\wb \in \Omb$, for all $\omb \in \Dcb^t_{\wb}$, we define the concatenated path $\wb \otimes_t \omb$ to be such that,
	for all $\phi \in C_b(\R_+ \x U)$,
	$$
		\big(\Theta_s, X_s, M_s(\phi) - M_t(\phi)\big)(\wb \otimes \omb) =
		\begin{cases}
			\big(\Theta_s, X_s, M_s(\phi) - M_t(\phi) \big)(\wb), ~~~ s \in [0,t); \\
			\big(\Theta_s, X_s, M_s(\phi) - M_t(\phi) \big)(\omb),~~~ s \in [t, \infty).
		\end{cases}
	$$
	
	Let $\P$ be a (Borel) probability measure on $\Omb$, and $\taub$ be a $\Fbb$-stopping time,	
	there is a family of regular conditional probability distribution (r.c.p.d.) $ \big( \Pb_{\wb} \big)_{\wb \in \Omb}$ w.r.t. $\Fcb_{\taub}$ such that 
	the $\Fcb_{\taub}$-measurable probability kernel $(\Pb_{\wb})_{\wb \in \Omb}$ satisfies
	$\Pb_{\wb} \big( \Dcb_{\taub(\wb), \wb} \big) = 1$ for every $\wb \in \Omb$.
	On the other hand, given a probability measure $\Pb$ defined on $(\Omb, \Fcb_{\taub})$ as well as a family of probability measures $\big( \Qb_{\wb} \big)_{\wb \in \Omb}$ such that $\wb \mapsto \Qb_{\wb}$ is $\Fcb_{\taub}$-measurable and $\Qb_{\wb} \big( \Dcb^{\taub(\wb)}_{\wb} \big) = 1$ for each $\wb \in \Omb$, 
	we can then define a unique concatenated probability measure $\Pb \otimes_{\taub} \Qb_{\cdot}$ by
	$$
		\Pb \otimes_{\taub} \Qb_{\cdot} (A) ~:= \int_{\Omb} \Pb( d \wb) \int_{\Omb} \1_A ( \wb \otimes_{\taub(\wb)} \omb ) \Qb_{\wb}( d\omb),
		~~\mbox{for all}~ A \in \Fcb_{\infty}.
	$$
	
	Next, let us recall some basic results about the (analytic) measurable selection theorem.
	In a Polish space $E$, a subset $A \subseteq E$ is called analytic if there is another Polish space $F$ and a Borel set $B \subseteq E \x F$ such that $A = \pi_E(B) = \{x ~: (x,y) \in B \}$.
	Notice that an analytic set is in general not Borel, but universally measurable, 
	i.e., it belongs to the $\sigma$-field obtained by completing the Borel $\sigma$-field under any arbitrary probability measure, 
	then it still makes sense to define the probability measure on the analytic sets.
	The class of all analytic sets is not a $\sigma$-field, we then also denote by $\Ac(E)$ the $\sigma$-field generated by all analytic sets.
	Next, a function $f: E \longrightarrow \Rb$ is said to be upper semi-analytic (u.s.a.) if $\{ x \in E ~: f(x) > c \}$ is analytic for every $c \in \R$.
	Let $F$ be some Polish space, a map $g: E \longrightarrow F$ is analytically measurable iff $g^{-1}(B) \in \Ac(E)$ for all Borel sets $B \in \Bc(F)$.

	\vspace{0.5em}

	With the above notions, we recall the following measurable selection theorem.

	\begin{Theorem} \label{thm:meas_selec}
		\rmi {\em Let $A \subseteq E \x F$ be analytic, $f : E \x \F \longrightarrow \Rb$ be u.s.a.
		Then the projection set $\pi_E(A)$ is still analytic and the function $x \longmapsto g(x) := \sup_{(x,y) \in A} f(x,y)$ is also u.s.a.

		\vspace{0.5em}

		\noindent \rmii For every $\eps > 0$, there is an analytically measurable map $\varphi_{\eps} : E\longrightarrow F$ 
		such that $\forall x \in \pi_{E}(A)$, $\varphi_{\eps}\in A_x$, 
		and  $f(x, \varphi_{\eps}(x)) \ge (g(x) - \eps) 1_{\{g(x) < \infty\}}+\frac{1}{\eps} 1_{\{g(x) = \infty\}}$.
		It follows that for any probability measure $\lambda$ on $E$},
		$$
			\int_E \! g(x) \lambda(dx)
			=
			\sup \Big\{ \int_E\,f(x, \varphi(x))\, \lambda(dx), ~\varphi \in \Ac_{usa} ~\mbox{s.t.}~ (x,\varphi(x)) \in A, ~\forall x \in \pi_E(A) \Big\}.
		$$
	\end{Theorem}
	Notice that $g$ is defined as the supremum of $f$, then the above equality is somehow an exchange property between the supremum and the integral, which is also the essential property appearing in the dynamic programming principle.

\paragraph{Optimization and dynamic programming principle}

	As the canonical space formulation of the optimal control/stopping problem in Section \ref{subsec:diffusion},
	we formulate the optimization problem on the canonical space $\Omb$, where a control (rule) is interpreted as  a probability measure on $\Omb$.

	\vspace{0.5em}

	Let $(\Pcb_{t, \xb})_{(t, \xb) \in \Rb_+ \x \Om}$ be a family of sets of (Borel) probability measures on $\Omb$,
	that is, $\Pcb_{t,\xb} \subset \Pc(\Omb)$ where $\Pc(\Omb)$ denotes the space of all (Borel) probability measures on $\Omb$.
	Namely, a probability measure $\Pb \in \Pcb_{t, \xb}$ is interpreted as a control/stopping rule, where $(t,\xb)$ is the initial condition, 
	and $\Pb$ describes the distribution of the controlled process, the stopping time, and also the control process itself.
	Given the {\bf reward functions} $L: \R_+ \x \Om \x U \longrightarrow \Rb$ and  $\Phi: \Rb_+ \x \Om \longrightarrow \Rb$, 
	the {\bf value function} $V$ of the optimization problem is then defined by, for all $(t, \xb) \in \Rb_+ \x \Om$,
	\be \label{eq:defV}
		V(t, \xb) 
		~:=
		\sup_{\Pb \in \Pcb_{t, \xb}} \E^{\Pb} \Big[ \int_t^{\Theta_{\infty}} \!\!\! \int_U L(s, X, u) M_s(du) ds +  \Phi \big(\Theta_{\infty}, X\big) \Big].
	\ee

	To obtain the dynamic programming principle, we will assume the following measurability condition, together with the stability conditions on the family $\big(\Pcb_{t, \xb} \big)_{(t, \xb) \in \Rb_+ \x \Om}$,
	which can be considered as an extension of the Markov property to the multi-valued probability measures $\Pcb_{t,\xb}$ case.

	\begin{Assumption}\label{assum:stability_Pch}
	\rmi For each $(t, \xb) \in \Rb_+ \x \Om$,
		the set $\Pcb_{t,\xb} = \Pcb_{t, \xb_{t \wedge \cdot}}$ is non-empty, 
		and $\Pb \big(X_{t \wedge \cdot} = \xb_{t \wedge \cdot}, \Theta_{\infty} \ge t \big) = 1$ for all $\Pb \in \Pcb_{t, \xb}$.
		Moreover, the graph set
		$$
			\big[ \big[ \Pcb \big]\big] := \big\{ (t, \xb, \Pb) \in \Rb_+ \x \Omb \x \Pc(\Omb) ~: \Pb \in \Pcb_{t, \xb} \big\}
			~\mbox{is analytic.}
		$$
	
	\noindent \rmii For all $(t_0, \xb_0) \in \Rb_+ \x \Om$, $\Pb \in \Pcb_{t_0, \xb_0}$ and $\taub$ a $\Fbb$-stopping time taking value in $[t_0, \infty]$, 
		with $A_{\taub} := \{ \omb \in \Omb ~: \Theta_{\infty}(\omb) > \taub(\omb) \}$, the following holds true.

		\noindent \rma There is a family of r.c.p.d. $\big(\Pb_{\wb} \big)_{\wb \in \Omb}$ of $\Pb$ w.r.t. $\Fcb_{\taub}$
		such that 
		$$
			\Pb_{\wb} \in \Pcb_{\taub(\wb), \w},
			~~\mbox{for $\Pb$-a.e.}~
			\wb = (\eta, \w, m) \in A_{\taub}.
		$$
		\noindent \rmb Let $(\Qb_{\wb})_{\wb \in \Omb}$ be a probability kernel from $(\Omb, \Fcb_{\taub})$ to $(\Omb, \Fcb_{\infty})$ such that $\wb \longmapsto \Qb_{\wb}$ is $\Fcb_{\taub}$-measurable,
		$\Qb_{\wb} = \Pb_{\wb}$ for $\Pb$-a.e. $\wb \in \Omb \setminus A_{\taub}$ with a family of r.c.p.d. $(\Pb_{\wb})_{\wb \in \Omb}$ of $\Pb$ w.r.t. $\Fcb_{\taub}$,
		and
		$\Qb_{\wb} \in \Pcb_{\taub(\wb), \w}$ 
		for $\Pb$-a.e. $\wb \in A_{\taub}$.
		Then $\Pb \ox_{\taub} \Qb_{\cdot} \in \Pcb_{t_0, \xb_0}$.
	\end{Assumption}




	\begin{Theorem} \label{theo:DPP}
		Let $(\Pcb_{t,\xb})_{(t,\xb) \in \Rb_+ \x \Om}$ be the family given above satisfying Assumption \ref{assum:stability_Pch}.
		Suppose in addition that the reward function $\Phi: \Rb_+ \x \Om \longrightarrow \Rb$ is upper semi-analytic, 
		and satisfies that $\Phi(t,\xb) = \Phi(t, \xb_{t \wedge \cdot})$ for all $(t, \xb) \in \Rb_+ \x \Om$.
		
		\vspace{0.5em}
		
		\noindent \rmi Then the value function $V:\Rb_+ \x \Om \longrightarrow \Rb$ defined by \eqref{eq:defV} is upper semi-analytic and in particular universally measurable,
		and $V(t, \xb) = V(t, \xb_{t \wedge \cdot})$ for all $(t, \xb) \in \Rb_+ \x \Om$.
		
		\vspace{0.5em}
		
		\noindent \rmii	For every $(t,\xb) \in \R_+ \x \Om$ and every $\Fbb$-stopping time $\taub$ taking value in $[t, \infty]$, one has the DPP
		\be \label{eq:DPP}
			V(t,\xb)
			\!\!&=&\!\!
			\sup_{\Pb \in \Pcb_{t, \xb}} \E^{\Pb}
			\Big[ 1_{\Theta_{\infty} \le \taub} \Big( \int_t^{\Theta_{\infty}} \!\!\!\! \int_U L(s, X, u) M_s(du) ds + \Phi \big( \Theta_{\infty}, X\big) \Big) \nonumber \\
			&&~~~~~~~~~~~~~+~
			1_{\Theta_{\infty} > \taub} \Big( \int_t^{\taub} \!\!\!\! \int_U L(s, X, u) M_s(du) ds + V \big( \taub, X \big) \Big) \Big].
		\ee
	\end{Theorem}
	\noindent {\sc Sketch of Proof.}
	\rmi Notice that with u.s.a. reward functions $L$ and $\Phi$, 
	the map 
	$$
		(t, \Pb) \longmapsto \E^{\Pb} \Big[ \int_t^{\Theta_{\infty}}\!\!\! \int_U L(s, X, u) M_s(du) ds +  \Phi \big(\Theta_{\infty}, X\big) \Big]
	$$
	is upper semi-analytic (see e.g. \cite[Corollary 7.48]{Bertsekas_1978}).
	Further, every $\Pcb_{t,\xb}$ is in fact a section set of the graph $[[\Pcb]]$, 
	and the supremum in \eqref{eq:defV} can be considered as a projection operator from functional space on $\Rb_+ \x \Om \x \Pc(\Omb)$ to that on $\Rb_+ \x \Om$.
	Then the measurability of $V$ follows by Theorem \ref{thm:meas_selec}.
	
	\vspace{0.5em}
	
	\noindent \rmii For the DPP in \eqref{eq:DPP}, by taking the conditioning and using Assumption \ref{assum:stability_Pch} (ii.a), it follows the inequality ``$\le$'' part of \eqref{eq:DPP}.
	To prove the reverse inequality ``$\ge$'', it is enough to take an arbitrary $\Pb \in \Pcb_{t,\xb}$,
	then to apply the measurable selection theorem to choose a ``measurable'' family of $\eps$-optimal control/stopping rules $(\Qb^{\eps}_{\wb})_{\wb \in A_{\taub}}$ for problems $V(\tau(\wb), \wb_{\tau(\wb) \wedge \cdot})$.
	Let $\Qb^{\eps}_{\wb} := \Pb_{\wb}$ for all $\wb \in \Omb\setminus A_{\taub}$ with a family of r.c.p.d. $(\Pb_{\wb})_{\wb \in \Omb}$ of $\Pb$ w.r.t. $\Fcb_{\taub}$.
	Applying the concatenation technique under Assumption \ref{assum:stability_Pch} $\mathrm{(ii.b)}$,
	one obtain $\Pb \otimes_{\taub} \Qb^{\eps}_{\cdot} \in \Pcb_{t,\xb}$ so that
	\b*
		V(t,\xb)
		\!\!\!&\ge&\!\!\!
		\E^{\Pb \otimes_{\taub} \Q^{\eps}_{\cdot} }
		\Big[ \int_t^{\Theta_{\infty}}  \!\!\!  \int_U L(s, X, u) M_s(du) ds + \Phi \big( \Theta_{\infty}, X\big) \Big] \\
		\!\!\!&\ge&\!\!\!
		\E^{\Pb}
		\Big[ 1_{\Theta_{\infty} \le \taub} \Big( \int_t^{\Theta_{\infty}} \!\!\!\! \int_U \! L(s, X, u) M_s(du) ds + \Phi \big( \Theta_{\infty}, X\big) \Big) \\
		&&~~~~~~~+
		1_{\Theta_{\infty} > \taub}  \Big( \int_t^{\taub} \!\!\!\! \int_U \! L(s, X, u) M_s(du) ds + V^{\eps} \big( \taub, X \big) \Big) \Big],
	\e*
	where $V^{\eps} : = (V-\eps) \1_{\{V<\infty\}} + \frac{1}{\eps} \1_{\{V= \infty \}}$.
	This concludes the proof of \eqref{eq:DPP} by arbitrariness of $\eps > 0$.
	\qed
 
\paragraph{Some direct consequences of the DPP}
 
 	As direct consequences of the dynamic programming principle, one obtains some characterizations of the value function as well as the optimal control/stopping rules.
	In particular, by choosing the stopping time $\taub$ in a local way, one can obtain local characterization of the value function, such as the viscosity solution property (see e.g Touzi \cite{Touzi}).
	
	\vspace{0.5em}
	
	Further, one can consider $\big(V(t, X) \big)_{t \ge 0}$ as process defined on $\Omb$, and the map $\Phi \longmapsto V$ as functional operator to explore their properties.
	For simplicity, let us assume that $L \equiv 0$ so that the DPP turns to be
	$$
		V(t,\xb)
		~=~
		\E^{\Pb}
		\Big[ 1_{\Theta_{\infty} \le \taub} \Phi \big( \Theta_{\infty}, X \big)
		+
		1_{\Theta_{\infty} > \taub} V \big( \taub, X_{\taub \wedge \cdot} \big)  \Big].
	$$
	Let us denote by  $\Ac^0_{usa}(\Rb_+ \x \Om)$ the set of all upper semi-analytic function bounded from below,
	and we say a function $\Psi \in \Ac^0_{usa}(\Rb_+ \x \Om)$ is ${\Pcb}$-super-median if $V(\Psi) \leq \Psi$ on $\Rb_+ \x \Om$.
	Let us also write $V(\Phi)$ in place of $V$ to emphasis its dependence on $\Phi$.
	\begin{Proposition} 
		\rmi The operator $\Phi \longmapsto V(\Phi)$ on $\Ac^0_{usa}(\Rb_+ \x \Om)$ is sub-linear.
 
 		\vspace{0.5em}
 
		\noindent \rmii For all $\Phi \in \Ac^0_{usa}(\Rb_+ \x \Om)$, $V(\Phi )$ is the smallest $\Pcb$-super-median function in $\Ac^0_{usa}(\Rb_+ \x \Om)$ greater than $\Phi$.

		\vspace{0.5em}

		\noindent\rmiii Assume in addition that $(V(t, X))_{t \ge 0}$ is a measurable process.
		Then it is a supermartingale under every probability measure $\P \in \Pcb_{0, \xb}$.
		Moreover, any probability measure $\Pb^* \in \Pcb_{0,\xb}$, under which $V(\Phi)$ is a martingale on $[0, \Theta_{\infty}]$, 
		is an optimal control/stopping rule for the optimization problem \eqref{eq:defV} with initial condition $(0, \xb)$.

	\end{Proposition}

	The above results can be derived easily from the DPP \eqref{eq:DPP}, by adapting the classical arguments (see e.g. \cite{ElKarouiSF}), whose proof is hence omitted.

\section{Dynamic programming principle of the optimal control and stopping problem}
\label{sec:MartPb}

	For an optimal control/stopping problem formulated on the canonical space,
	the essential point is to check the measurability and stability conditions in Assumption \ref{assum:stability_Pch} to deduce the dynamic programming principle.
	In the following, we will study an optimal control/stopping problem with a martingale problem formulation, and then check Assumption \ref{assum:stability_Pch} in this framework.
	In particular, it covers the controlled/stopped diffusion processes problem illustrated in Section \ref{subsec:diffusion}.
	
	\vspace{0.5em}
	
	Recall that $U$ and $E$ are both (non-empty) Polish spaces, $C(E)$ denotes the class of all continuous functions defined on $E$, and $C_b(E)$ is the subset of all bounded continuous functions.

\subsection{Generators and a controlled/stopped martingale problem}
\label{subsec:generator}

	We will first recall some basic facts on the Markov process, its generator as well as the associated martingale problem,
	and then introduce a general optimal control/stopping problem with a martingale problem formulation.

\paragraph{Markov process and generator}

	Let $P = (P_t)_{t \ge 0}$ be a family of homogeneous transition kernels on $(E, \Bc(E))$, which forms a semi-group on an appropriate functional space.
	Then on a filtered probability space $\big(\Om^*, ~\Fc^*, ~\F^*=(\Fc^*_t)_{t \ge 0} \big)$ rich enough, with any probability measure $\lambda$ on $(E, \Bc(E))$, 
	one can construct a continuous-time Markov process $(X^*, \P^*_{\lambda})$ w.r.t. $\F^*$ with transition kernels $(P_t)_{t \ge 0}$ and initial distribution $\lambda$, i.e., for every bounded measurable function $f: E \longrightarrow \R$,
	\b*
		\E^{\P^*_{\lambda}} \big[ f(X^*_t) \big| \Fc_s^* \big] = P_{t-s} f(X^*_s)
		&\mbox{and}&
		\P^*_{\lambda} \circ (X^*_0)^{-1} = \lambda,
	\e*
	for every $0 \le s \le t$.
	When the initial distribution is given by the Dirac measure on $x \in E$, we denote $\P^*_x := \P^*_{\delta_x}$.
	For the Markov process $X^*$,  its ``infinitesimal'' generator $\Lc$ is defined by
	\b*
		\Lc f &:=& \lim_{t \searrow 0} \frac{1}{t} \big(P_t f - f \big),
	\e*
	where $f \in C(E)$ is said to lie in the domain $\Dc_{\Lc}$ of the generator $\Lc$ whenever the above limit is well defined.
	Following the language of Ethier and Kurtz \cite{Ethier_Kurtz}, we also call its graph $\G := \{ (f, \Lc f) ~: f \in \Dc_{\Lc} \}$ as the ``full'' generator.
	It follows that for every $f \in \Dc_{\Lc}$ (equivalently for every $(f,g) \in \G$), the process
	\begin{equation} \label{eq:fg_mgt_pb}
		f(X^*_t) - \int_0^t \Lc f(X^*_s) ds
		~~\Big(\mbox{or equivalently}~~
		f(X^*_t) - \int_0^t g(X^*_s) ds \Big)
	\end{equation}
	is a $\F^*$-martingale under $\P^*_{\lambda}$ for every initial distribution $\lambda$.
	Then the martingale problem with the ``infinitesimal'' generator $\Lc$ (resp. ``full'' generator $\G$)
	consists in finding a probability space together with a process $X^*$ such that the process in \eqref{eq:fg_mgt_pb} is a (local) martingale for all $f \in \Dc_{\Lc}$ (resp. for all $(f,g) \in \G$).
	On the other hand, given the existence and uniqueness of solutions to the martingale problems,
	one can also construct the associated Markov process from solutions of the martingale problems (see Ethier and Kurtz \cite{Ethier_Kurtz} for more details).
	In the context of control problems, it seems to be more convenient to use the martingale problem formulation comparing to the semi-group formulation (see Example \ref{exam:Nisio}).

	\vspace{0.5em}

	Let us provide below some examples of the Markov processes as well as the associated martingale problems.

	\begin{Example}[Continuous-time Markov chain]
	{\rm Let $E$ be a countable space, for a $E$-valued continuous-time Markov chain with transition rate matrix $Q$, the infinitesimal generator of $X^*$ is given by
	$
		\Lc^1 \varphi(x) := \sum_{y \neq x} Q(x,y) \big( \varphi(y) - \varphi(x) \big),
	$
	where the domain $\Dc_{\Lc^1}$ is the class of all bounded functions from $E$ to $\R$,
	and hence the full generator is given by $\big\{ (\varphi, \Lc^1 \varphi) ~: \varphi \in \Dc_{\Lc^1} \big\}$.
	}
	\end{Example}

	\begin{Example}[Diffusion process]
	{\rm The diffusion process is an important example of a Markov process.
	Let $E = \R^d$, $\mu: \R^d \longrightarrow \R^d$ and $\sigma: \R^d \longrightarrow \S^d$, 
	and $X$ be the diffusion process defined by the SDE
	$$
		dX_t = \mu(X_t) dt + \sigma(X_t) dW_t,
	$$
	for some Brownian motion $W$.
	Its generator is then given by
	\begin{equation} \label{eq:def_Lc_diffusion}
		\Lc^2 \varphi (x) := \mu(x) \cdot D \varphi(x) + \frac{1}{2} \sigma \sigma^T(x) : D^2\varphi(x),
	\end{equation}
	with the domain $\Dc_{\Lc^2} := C_b^2(\R^d)$, i.e. the class of all bounded continuous functions admitting bounded continuous first and second order derivatives.
	Similarly, its full generator is provided by $\big\{ (\varphi, \Lc^2 \varphi) ~: \varphi \in \Dc_{\Lc^2} \big\}$.
	When $\mu$ and $\sigma$ are both bounded continuous, the corresponding martingale problem has existence of solutions.
	While in general the uniqueness fails, one can apply the Markovian selection approach to construct a Markov process as solution (see e.g. \cite{StroockVaradhan} for details).
	}
	\end{Example}

	\begin{Example}[Reflected diffusion process] \label{exam:reflect_diffusion}
	{\rm
	Let $O \subset \R^d$ be a bounded open set with smooth boundary $\partial O$.
	Let $0 < \beta \le 1$, denote by  $C^{1, \beta}(\partial O)$ the class of all continuous functions defined on $\partial O$ having $\beta$-H\"older first order derivatives,
	and by $C^{2,\beta}(\overline O)$ the collection of all $C^2(\overline O)$ functions $\varphi$ such that $D^2 \varphi$ is $\beta$-H\"older.
	We consider a reflected diffusion process, which is a diffusion process with generator \eqref{eq:def_Lc_diffusion} in $O$ and reflects on $\partial O$ with reflection direction given by $c\in C^{1, \beta}(\partial O)$ satisfying $\inf_{x \in \partial O} \langle c(x), n(x) \rangle > 0$, where $n(x)$ denotes the outward unit normal to $\partial O$ at $x$.
	Under sufficient regularity conditions on $\partial O$ as well as on $\mu$, $\sigma$ and $c$,  then the closure of
	\b*
		\big\{( \varphi, \Lc^2 \varphi) ~: \varphi \in C^{2,\beta}(\overline O), 
		~ c \cdot \nabla \varphi = 0 ~\mbox{on}~ \partial O \big\}
	\e*
	in $C(O) \x C(O)$ under the $L^{\infty}$-norm provides a full generator for the associated  reflected diffusion process (see e.g. Chapter 8.1 of Ethier and Kurtz \cite{Ethier_Kurtz}).
	}
	\end{Example}

	\begin{Example}[Branching Brownian motion] \label{exam:branching}
	{\rm
	Let $\beta > 0$, $(p_k)_{k \ge 0}$ be a probability sequence, i.e. $p_k \ge 0$ for every $k \ge 0$ and $\sum_{k=0}^{\infty} p_k = 1$.
	We consider a particle system, where each particle moves as a Brownian motion in $\R^d$, at exponential time of intensity $\beta$, it branches into $k$ (conditional) independent particles with probability $p_k$.
	Assume further that the increments of all particles are taken to be independent and independent to the lifetime and the numbers of offspring particles.
	By considering the measure induced by all particles in the system, 
	one obtains a measure-valued (branching) process, whose state space is given by
	$$
		E ~~:=~~ \Big\{ \sum_{i=1}^k \delta_{x_i} ~ : k = 0, 1,2, \cdots, x_i \in \R^d \Big\}.
	$$
	Notice that $E$ is clearly a closed subset of the space of finite, positive, Borel measures on $\R^d$ under the weak convergence topology.
	Then following Chapter 9.4 of \cite{Ethier_Kurtz}, a full generator of the above branching Brownian motion is given by
	\b*
		\Big\{ \Big( e^{\langle \log \varphi, \cdot \rangle}, ~ e^{\langle \log \varphi, \cdot \rangle} \Big\langle \frac{\frac{1}{2}\Delta \varphi + \beta \big( \sum_{k=0}^{\infty} p_k \varphi^k - \varphi \big)}{\varphi}, \cdot \Big \rangle \Big)
		&:
		\varphi \in C^{2,+}_b(\R^d), ~ | \varphi|_{\infty} < 1 \Big\},
	\e*
	where $\Delta$ is the Laplacian and $C^{2,+}_b(\R^d)$ denotes the collection of all strictly positive functions in $C^2_b(\R^d)$.
	}
	\end{Example}

	\begin{Remark}{\rm
		Since the transition kernels are linear operators on the functional space on $E$, it follows that the ``infinitesimal'' generator is also linear.
		Therefore, the ``full'' generator is generally composed by couples $(f,g)$ of functions, where $g$ depends linearly on $f$.
		Nevertheless, for some Markov processes, it is more convenient to use the ``full'' generator formulation, such as the reflected diffusion process in Example \ref{exam:reflect_diffusion}.
	}
	\end{Remark}

\paragraph{A controlled/stopped martingale problem}

	One of the most classical control problems is the controlled Markov processes problem (see e.g. \cite{Krylov}, etc.), which can be obtained by adding a control component in the generator of the Markov processes.
	For ease of presentation, we shall use the notion of ``full'' generator.
	More importantly, we shall present the control problem in a time and path dependent setting,
	which leads to the fact that the ``full'' generator $\G$ being a subset of $C_b(E) \x B(\R_+ \x \Om \x U \x E)$,
	where $B(\R_+ \x \Om \x U \x E)$ denotes the space of all measurable functions $g: \R_+ \x \Om \x U \x E \to \R$ such that 
	\begin{equation}
		\int_0^T \sup_{u \in U} \big| g(t, \xb_{t \wedge \cdot}, u, \xb_t) \big| ~ dt ~<~ \infty,
		~~\mbox{for all}~ T \ge 0 ~\mbox{and}~ \xb \in \Om.
	\end{equation}

	As illustrated in Section \ref{subsec:diffusion}, we will formulate the problem directly on the canonical space $\Omb = \Rb_+ \x \Om \x \M$, i.e. the control rules are interpreted as probability measures on $\Omb$.
	Given $(f, g) \in C_b(E) \x B(\R_+ \x \Om \x U \x E)$, let us define
	\be
		C_t(f,g) ~:=~ f(X_t) - \int_0^t \int_U g(s, X_{s \wedge \cdot}, u, X_s) M(s, du) ds,
	\ee
	which is clearly a right-continuous $\Fbb$-adapted process.
	For any $(f,g) \in  C_b(E) \x B(\R_+ \x \Om \x U \x E)$, let us define also a sequence of localized (bounded) process by
	$$
		C^n_t (f,g)  ~:=~ C_{\tau_n \wedge t} (f,g),
		~~\mbox{with}~
		\tau_n ~:=~ \inf \big\{ t \ge 0 ~: |C_t(f,g)| \ge n \big\}.
	$$

	\begin{Definition} \label{def:weak_relaxed_control}
		Let $\G \subset C_b(E) \x B(\R_+ \x \Om \x U \x E)$ be a ``full'' generator of the control problem,
		and $(t, \xb) \in \Rb_+ \x \Om$.
		
		\vspace{0.5em}
		
		\noindent \rmi A relaxed control/stopping rule, associated with generator $\G$ 
		and initial condition $(t, \xb) \in \R_+ \x \Om$, is a probability measure $\Pb$ on $(\Omb, \Fcb_{\infty})$ such that
		$\Pb \big[ \Theta_{\infty} \ge t, ~ X_s = \xb_s, ~0\le s \le t \big] = 1$,
		and under which the process
		$\big( C^n_s(f,g) \big)_{s \ge t}$ is a $\Fbb$-martingale (and hence a martingale w.r.t. the augmented fitlration $\Fbb^{\Pb}_+$) for every $(f,g) \in \G$ and all $n \ge 1$.
		Further, when $t = \infty$, a probability measure $\Pb$ is called a relaxed control/stopping rule with inital condition $(t, \xb)$ if $\Pb \big[ \Theta_{\infty} = \infty, ~X_s = \xb_s,~ s\in \R_+ \big] = 1$.
		Denote
		$$
			\Pcb_{t, \xb}
			~:=~
			\big\{
				\mbox{All relaxed  rules with generator $\G$ and initial condition $(t, \xb)$}
			\big\}.
		$$

		\noindent \rmii A weak control/stopping rule associated with generator $\G$ and initial condition $(t, \xb)$ is a probability measure $\Pb \in \Pcb_{t, \xb}$ such that $\Pb \big[ M \in \M^0 \big]  = 1$ (see \eqref{eq:def_M0} for the definition of $\M_0$). Denote
		$$
			\Pcb^0_{t, \xb} ~:=~ \big\{ \Pb \in \Pcb_{t,\xb} ~: \Pb \big[ M \in \M^0 \big]  = 1 \big\}.
		$$

		\vspace{0.5em}

		\noindent \rmiii We say $\G$ is countably generated, if there exists a countable subset $\G_0\subseteq \G$ such that
		every $\G_0$-relaxed control/stopping rule is a $\G$-relaxed control/stopping rule.
	\end{Definition}
	Let $\Phi: \Rb_+ \x \Om \longrightarrow \Rb$ and $L: \R_+ \x \Om \x U \longrightarrow \Rb$
	be upper semi-analytic satisfying $\Phi(t, \xb) = \Phi (t, \xb_{t \wedge \cdot})$ and $L(t, \xb, u) = L(t, \xb_{t \wedge \cdot}, u)$ for all $(t, \xb, u) \in \R_+ \x \Om \x U$, we then define
	\be \label{eq:value_function_V}
		V(t,\xb) ~:=~ \sup_{\Pb \in \Pcb_{t,\xb}} 
		\E^{\Pb} \Big[ \int_t^{\Theta_{\infty}} \!\! \int_U L(s, X, u) M_s(du) ds +  \Phi(\Theta_{\infty}, X) \Big],
	\ee
	and
	\be
		V_0(t,\xb) ~:=~ \sup_{\Pb \in \Pcb^0_{t,\xb}} 
		 \E^{\Pb} \Big[ \int_t^{\Theta_{\infty}} \!\! \int_U L(s, X, u) M_s(du) ds +  \Phi(\Theta_{\infty}, X) \Big].
	\ee

	In the above abstract formulation,
	we do not discuss the conditions on the generator $\G$ to make the problem well-posed.
	It is possible, in general, that the martingale problem in Definition \ref{def:weak_relaxed_control} has no solution or has multiple solutions with an arbitrary generator.
	For concrete problems, one can formulate more explicit conditions to ensure the existence of solutions to the martingale problem, such as the controlled diffusion processes problem in Section \ref{subsec:dpp_diffusion} below.
	In any case, with the convention that $\sup_{\emptyset} = - \infty$, the above function $V$ and $V_0$ are well defined.

\paragraph{More discussions on the weak/relaxed formulation}

	The above weak or relaxed control problem is usually formulated in a different but equivalent way.
	Given a generator $\G$ and initial condition $(t, \xb) \in \R_+ \x \Om$, a weak (resp. relaxed) control/stopping term $\alpha$ is a term
	\b*
		\alpha &=& \big( \O^{\alpha}, \Fc^{\alpha}, \F^{\alpha}, \P^{\alpha} , X^{\alpha}, \tau^{\alpha}, \nu^{\alpha} ~(\mbox{resp.}~ m^{\alpha})\big),
	\e*
	where $(\O^{\alpha}, \Fc^{\alpha}, \F^{\alpha}, \P^{\alpha})$ is a filtered probability space equipped with
	an adapted $E$-valued c\`adl\`ag process $X^{\alpha}$ such that $X^{\alpha}_{t \wedge \cdot} = \xb_{t \wedge \cdot}$,
	and a stopping time $\tau^{\alpha}$ taking value in $[t, \infty]$,
	together with a $U$-valued (resp. $\Pc(U)$-valued) progressively measurable control process $\nu^{\alpha} = (\nu^{\alpha}_t)_{t \ge 0}$
	(resp. $m^{\alpha} = (m^{\alpha}_t)_{t \ge 0}$),
	such that the process $(C^{\alpha}_{s \land \tau^{\alpha}}(f,g))_{s \ge t}$ given below is a local martingale for every couple $(f,g) \in \G$,
	\b*
		C^{\alpha}_s(f,g)
		~:=~
		f( X^{\alpha}_s) - \int_0^s g \big(r,X^{\alpha}_{r \wedge \cdot}, \nu^{\alpha}_r ~(\mbox{resp.}~ m^{\alpha}_r), X_r^{\alpha} \big) dr,
	\e*
	where $g(r, \xb, m^{\alpha}_r, x) := \int_U g(r, \xb, u, x) m^{\alpha}_r(du)$.
	To see their equivalence, it is enough to notice that any weak (resp. relaxed) term $\alpha$ induces a weak (resp. relaxed) rule probability on $\Omb$; and in contrast, any weak (resp. relaxed) rule $\Pb$ together with the canonical space $\Omb$ and the augmented filtration $\Fbb^{\Pb}_+$ is a weak (resp. relaxed) term
	(see e.g. Proposition \ref{prop:weak_reformulation}).

	\begin{Remark}[On the relaxed control]
	{\rm The relaxed control/stopping rule consists in replacing the $U$-valued control process by a $\Pc(U)$ measure-valued processes. 
	This technique has been largely used in deterministic control problem to obtain the closeness and convexity of set of controls.
	In the stochastic control of diffusion processes setting, the relaxed formulation has initially been introduced by Fleming \cite{Fleming_1978}, and by El Karoui, Huu Nguyen and Jeanblanc \cite{ElKaroui_1987} in order to obtain the existence of optimal control rules.
	}
	\end{Remark}

	\begin{Remark}[Comparison with Nisio semi-group formulation]
	{\rm
	The ``full'' generator $G$ is fixed in the above martingale problem formulation;
	restricted to the controlled Markov processes case, this implies that the domain of generator should be the same for all controls.
	From this point of view, the above formulation is more restrictive comparing to the Nisio semi-group formulation illustrated in Example \ref{exam:Nisio}, where one can consider a larger class of different generators (or equivalently semi-groups) for the controlled Markov processes.
	}
	\end{Remark}

\subsection{The dynamic programming principle}
\label{subsec:DPP}

	We now show that the  family $\Pcb_{t, \xb}$ (resp. $\Pcb^0_{t,\xb}$) in Definition \ref{def:weak_relaxed_control}
	satisfies Assumption \ref{assum:stability_Pch}, which implies the corresponding dynamic programming principle.
	Moreover, let $\lambda$ be a (Borel) probability measure on $E$, similar to Definition \ref{def:weak_relaxed_control},
	we say that a probability $\Pb$ on $\Omb$ is a relaxed control/stopping rule with initial distribution $\lambda$,
	if $X_0 \sim^{\Pb} \lambda$ and $\big(C^n_s(f,g) \big)_{s \ge 0}$ is a martingale for every $(f,g) \in \G$ and $n \ge 1$,
	and $\Pb$ is a weak control/stopping rule if it satisfies in addition that $\Pb \big[ M \in \M^0 \big] = 1$.
	Let us denote by $\Pcb(\lambda)$ (resp. $\Pcb^0(\lambda)$) the collection of all relaxed (resp. weak) control/stopping rules with initial distribution $\lambda$,
	and then define
	$$
		V(\lambda) ~:= \sup_{\Pb \in \Pcb(\lambda)} 
		\E^{\Pb} \Big[ \int_0^{\Theta_{\infty}} \!\! \int_U L(s, X, u) M_s(du) ds +  \Phi(\Theta_{\infty}, X) \Big],
	$$
	and
	$$
		V_0(\lambda) ~:= \sup_{\Pb \in \Pcb^0(\lambda)} 
		\E^{\Pb} \Big[ \int_0^{\Theta_{\infty}} \!\! \int_U L(s, X, u) M_s(du) ds +  \Phi(\Theta_{\infty}, X) \Big].
	$$	

	\begin{Theorem} \label{theo:ctrl_stopping_DPP}
		Assume that $\G$ is countably generated, and $\Phi$ and $L$ are upper semi-analytic and such that
		$\Phi(t, \xb) = \Phi (t, \xb_{t \wedge \cdot})$ and $L(t, \xb, u) = L(t, \xb_{t \wedge \cdot}, u)$ for all $(t, \xb, u) \in \R_+ \x \Om \x U$.
		
		\noindent \rmi  Then the value function $V_0: \Rb_+ \x \Om \longrightarrow \Rb$ is upper semi-analytic,
		and for every $\Fbb$-stopping time $\taub$ taking value in $[t,\infty]$, one has
		\be \label{eq:DPP0}
			V_0(t,\xb)
			&\!\!\!=\!\!\!\!&
			\sup_{\Pb \in \Pcb^0_{t,\xb}} ~ \E^{\Pb}
			~\Big[
				\Big( \int_t^{\Theta_{\infty}} \!\!\!\!\int_U L(s, X, u) M_s(du) ds+
				\Phi( \Theta_{\infty}, X)  \Big) 1_{\Theta_{\infty} \le \taub} \nonumber \\
			&&~~~~~~~~~~~~~~~~+~
				\Big( \int_t^{\taub} \!\!\!\int_U L(s, X, u) M_s(du) ds+ V_0(\taub, X) \Big) 1_{\Theta_{\infty} > \taub}
			\Big].~~~
		\ee
		Moreover, assume in addition that $\Pcb^0_{0, \xb}$ is nonempty for all $\xb \in \Om$, then the set $\Pcb^0(\lambda)$ is nonempty for a Borel probability measure $\lambda$ on $E$, and
		$$
			V_0(\lambda) = \int_E V_0(0, x) \lambda(dx).
		$$

		\noindent \rmii The results hold true if one replaces $(V_0, \Pcb^0)$ by $(V, \Pcb)$ in the above statement.
	\end{Theorem}

	For the proof, we will only consider the statement for $V_0$ since the arguments are the same for $V$.
	Notice that it is clear that the family $(\Pcb^0_{t,\xb})_{(t,\xb) \in \Rb_+ \x \Om}$ satisfies that $\Pcb^0_{t,\xb} = \Pcb^0_{t,\xb_{t\wedge \cdot}}$ for all $(t,\xb) \in \Rb_+ \x \Om$.
	Then in view of Theorem \ref{theo:DPP}, it is enough to prove the following two lemmas 
	(Lemmas \ref{lemm:ctrl_stopping_measurability} and \ref{lemm:ctrl_stopping_stability})
	in order to conclude the proof of Theorem \ref{theo:ctrl_stopping_DPP}.

	\begin{Lemma} \label{lemm:ctrl_stopping_measurability}
		Suppose that $\G$ is countably generated. 
		Then $\big[\big[ \Pcb^0 \big]\big]$ defined below is Borel measurable in the Polish space 
		$\Rb_+ \x \Om \x \Pc(\Omb)$,
		\b*
			\big[\big[ \Pcb^0 \big]\big]
			&:=&
			\Big \{ (t, \xb, \Pb) \in \Rb_+ \x \Om \x \Pc(\Omb)
			~: (t,\xb) \in \Rb_+ \x \Om,~ \Pb \in \Pcb^0_{t, \xb} \Big\}.
		\e*
	\end{Lemma}
	\proof Let $0 \le r \le s$, $\xi \in C_b(\Ob, \Fcb_r)$ and $(f,g) \in \G$,
	we introduce some subsets in $\Rb_+ \x \Om \x \Pc(\Omb)$ as follows.
	Let $\Ab^0 := \big\{ (t,\xb,\Pb) \in \Rb_+ \x \Om \x \Pc(\Ob) ~: \Pb (M \in \M^0,~ \Theta_{\infty} \ge t) = 1 \big \}$,
	$\Ab^1_s := \big \{ (t, \xb, \Pb) ~: \Pb \big(X_{s \land t} = \xb ( s \land t) \big) = 1 \big \}$
	and
	\b*
		\Ab^{2,n}_{r, s, \xi, f,g}
		&:=&
		\Big\{ (t, \xb,\Pb) ~: 
		\E^{\Pb} \Big[ \big( C^n_{s \land t}(f,g) - C^n_{r \land t}(f,g) \big) \xi \Big] = 0 
		\Big \},
	\e*
	which are all Borel measurable since $\M^0$ is a Borel measurable set in $\M$ and $C^n(f,g)$ is c\`adl\`ag $\Fbb$-progressively measurable.
	It follows that $\big[\big[ \Pcb^0 \big]\big]$ is also Borel measurable
	since it is the intersection of $\Ab^0$, $\Ab^1_s$ and $\Ab^{2,n}_{r, s, \xi, f,g}$,
	where $n \ge 1$, $r \le s$ vary among rational numbers in $\R_+$,
	$\xi$ varies among a countable dense subset of $C_b(\Ob, \Fcb_r)$ and $(f,g)$ varies among the countable set $\G_0$ which generates $\G$. \qed

	\begin{Lemma} \label{lemm:ctrl_stopping_stability}
		Suppose that $\G$ is countably generated, and $\Pcb^0_{t,\xb}$ is nonempty for every $(t,\xb) \in \Rb_+ \x \Om$.
		Let $(t_0,\xb_0) \in \R_+ \x \Om$, $\Pb \in \Pcb^0_{t_0,\xb_0}$
		and $\taub$  be a $\Fbb$-stopping time taking value in $[t_0, \infty]$, denoting $A_{\taub} := \{ \omb \in \Omb ~:  \Theta_{\infty} (\omb) > \taub(\omb) \}$. \\
		\rmi Then there exists a family of r.c.p.d. $(\Pb_{\omb})_{\omb \in \Omb}$ of $\Pb$ w.r.t. $\Fcb_{\taub}$
		such that $\Pb_{\omb} \in \Pcb^0_{\taub(\omb), \om}$ for $\Pb$-almost every $\omb = (\theta, \om, m) \in A_{\taub}$.\\
		\rmii Let $(\Qb_{\omb})_{\omb \in \Omb}$ be such that $\omb \mapsto \Qb_{\omb}$ is $\Fcb_{\taub}$-measurable,
		$\Qb_{\omb} = \Pb_{\omb}$ for $\Pb$-a.e. $\omb \in \Omb \setminus A_{\taub}$ with a family of r.c.p.d. $(\Pb_{\omb})_{\omb \in \Omb}$ of $\Pb$ w.r.t. $\Fcb_{\taub}$,
		and $\Qb_{\omb} \in \Pcb^0_{\taub(\omb), \om}$ for $\Pb$-a.e. $\omb = (\theta, \om, m) \in A_{\taub}$,
		then $\Pb \ox_{\taub} \Qb_{\cdot} \in \Pcb^0_{t_0, \xb_0}$.
	\end{Lemma}
	\proof
		Let $(t_0, \xb_0) \in \R_+ \x \Om$, $\taub$ be a $\Fbb$-stopping time taking value in $[t_0, \infty]$ and $\Pb \in \Pcb^0_{t_0,\xb_0}$. \\
		\rmi Since $\Fcb_{\taub}$ is countably generated,
		there is a family of r.c.p.d. $ (\Pb_{\omb} )_{\omb \in \Omb}$ of $\Pb$ w.r.t. $\Fcb_{\taub}$.
		In particular, $\Pb_{\omb} [ \Theta_{\infty} \ge \taub(\omb), ~M \in \M_0 ] = 1$ for $\Pb$-a.e. $\omb \in A_{\taub}$, and
		\b*
			\Pb_{\omb} \big[ X_s=\om_s, ~ s \in [0, \taub(\omb)] \cap \R_+ \big] ~=~ 1, 
			~~\mbox{for every}~~ \omb = (\theta, \om, m) \in \Omb.
		\e*
		Moreover, since $( C^n_t(f,g) )_{t \ge t_0}$ is a $\Pb$-martingale on $[t_0, \infty)$ for every $(f,g) \in \G_0$,
		it follows by Theorem 1.2.10 of Stroock and Varadhan \cite{StroockVaradhan} that
		there is $\Pb$-null set $N^n_{f,g} \in \Fcb_{\taub}$
		such that $C^n (f,g)$ is $\Pb_{\omb}$-martingale after time $\taub(\omb)$ for every $\omb \notin N^n_{f,g}$ such that $\taub(\omb) < \infty$.
		Using the fact that $\G_0$ is countable, $N^n := \cup_{(f,g) \in \G_0} N^n_{f,g}$ is $\Pb$-null set such that
		$C^n(f,g)$ is a $\Pb_{\omb}$-martingale after time $\taub(\omb)$ for every $\omb \in \Ob \setminus  N^n$ and every $(f,g)\in \G_0$.
		And hence $\Pb_{\omb} \in \Pcb^0_{\taub(\omb), \om}$ for every $\omb = (\theta, \om, m) \in A_{\taub} \setminus N$ with $N = \cup_{n \ge 1} N^n$.\\
		\rmii By the definition of $(\Pcb^0_{t,\xb})_{(t,\xb)\in \Rb_+ \x \Om}$, 
		we notice that $\Qb_{\omb} \in \Pcb^0_{\taub(\omb), \om}$ implies that $\delta_{\omb} \ox_{\taub(\omb)} \Qb_{\omb} \in \Pcb^0_{\taub(\omb), \om}$
		for all $\omb = (\theta, \om, m) \in A_{\taub}$.
		In particular, $(\delta_{\omb} \ox_{\taub(\omb)} \Qb_{\omb})_{\omb \in \Omb}$
		is a family of r.c.p.d. of $\Pb \ox_{\taub} \Qb_{\cdot}$ w.r.t. $\Fcb_{\taub}$,
		and under each $ \delta_{\omb} \ox_{\taub(\omb)} \Qb_{\omb}$, $\big(C^n_s(f,g) \big)_{s \ge \taub(\omb)}$ is a bounded c\`adl\`ag martingale, for every $(f,g) \in \G$.
		Then still by Theorem 1.2.10 of \cite{StroockVaradhan}, it follows that $\Pb \ox_{\taub} \Qb_{\cdot}$ solves the martingale problem, and hence $\Pb \ox_{\taub} \Qb_{\cdot} \in \Pcb_{t_0, \xb_0}$. \qed

\subsection{The controlled/stopped diffusion processes problem}
\label{subsec:dpp_diffusion}

	Let us now apply the results in Theorem \ref{theo:ctrl_stopping_DPP} to the controlled/stopped diffusion processes problem with 
	coefficient functions $(\mu, \sigma) : \R_+ \x \Om \x U \longrightarrow \R^d \x \S^d$
	(see Section \ref{subsec:diffusion}),
	where $\Om := \D(\R_+, E)$ with $E := \R^d$.
	Recall also that $\Omb := \R_+ \x \Om \x \M$.
	We will first study the problem under the following technical integrability condition \eqref{eq:mu_sigma_integ}, that is, for all $\xb \in \Om$ and $T \ge 0$,
	\begin{equation} \label{eq:mu_sigma_integ_p}
		\int_0^T \sup_{u \in U} \Big( | \mu(t, \xb, u) + \| \sigma(t, \xb, u) \|^2 \Big) dt < \infty.
	\end{equation}
	Then in Section \ref{subsubsec:unbound_mu_sigma}, we also discuss how to relax this technical condition.

\subsubsection{The weak and relaxed formulation}
\label{subsubsec:weak_relax_form}

	Let $(t, \xb) \in \Rb_+ \x \Om$ be the initial condition,
	we follow Definition \ref{def:weak_control} in Sections \ref{subsec:diffusion} to introduce the weak control in the controlled diffusion processes setting.
	Concretely, for $(t, \xb) \in \R_+ \x \Om$, a weak control (of diffusion process) with initial condition $(t, \xb)$ is a term 
	$\alpha = (\Om^{\alpha}, \Fc^{\alpha}, \P^{\alpha}, \F^{\alpha}, \pi^{\alpha}, X^{\alpha}, B^{\alpha}, \nu^{\alpha})$,
	where $(\Om^{\alpha}, \Fc^{\alpha}, \P^{\alpha}, \F^{\alpha})$ is a filtered probability space, 
	equipped with a stopping time $\pi^{\alpha}$, a $d$-dimensional Brownian motion $B^{\alpha}$, and a $U$-valued predictable process $\nu^{\alpha}$, together with a continuous adapted process $X^{\alpha}$ such that
	$X^{\alpha}_{t \wedge \cdot} = \xb_{t \wedge \cdot}$, a.s. and
	$$
		X^{\alpha}_s = \xb_t + \int_t^s \mu(r, X^{\alpha}, \nu^{\alpha}_r) dr + \int_t^s \sigma(r, X^{\alpha}, \nu^{\alpha}_r) d B^{\alpha}_r,
		~~s \ge t, ~\mathrm{a.s.}
	$$
	When $t = \infty$, we say a term $\alpha$ is a weak control (of diffusion process) with initial condition $(t, \xb)$ if $X^{\alpha} = \xb$, a.s.
	Let us denote by $\Ac_W(t, \xb)$ the collection of all weak controls (of diffusion process) with initial condition $(t,\xb) \in \Rb_+ \x \Om$.
	For $\alpha \in \Ac_W(t, \xb)$, we denote $M^{\alpha}(du, ds) := \delta_{\nu^{\alpha}_s}(du)ds$.
	
	\vspace{0.5em}
	
	Comparing to Definition \ref{def:weak_control}, 
	one just replaces the initial condition $(0, x_0)$ in Definition \ref{def:weak_control} by $(t, \xb)$ in above.
	Similarly, by changing the initial condition in Definition \ref{def:relaxed_control}, one can define the relaxed control (of diffusion process) with initial condition $(t, \xb)$,
	and denote the corresponding set by $\Ac_R(t, \xb)$ for all $(t, \xb) \in \Rb_+ \x \Om$.
	Then with the reward functions $\Phi: \Rb_+ \x \Om \longrightarrow \Rb$ and $L: \R_+ \x \Om \x U \longrightarrow \Rb$,
	let us introduce the value functions of the weak and relaxed formulation of the controlled diffusion processes problem:
	$$
		V_W(t,\xb) ~:= \sup_{\alpha \in \Ac_W(t,\xb) } \E^{\P^{\alpha}} \Big[
			\int_t^{\pi^{\alpha}} L(s, X^{\alpha}, \nu^{\alpha}_s) ds
			+
			\Phi( \pi^{\alpha}, X^{\alpha}) 
		\Big],
	$$
	and
	$$
		V_R(t,\xb) ~:= \sup_{\alpha \in \Ac_R(t,\xb) } \E^{\P^{\alpha}} \Big[
			\int_t^{\pi^{\alpha}} \!\!\! \int_U L(s, X^{\alpha}, u) M^{\alpha}_s(du) ds
			+
			\Phi( \pi^{\alpha}, X^{\alpha}) 
		\Big].
	$$

	\begin{Theorem} \label{theo:DPP_VWVR}
		Assume that the coefficient functions $\mu$ and $\sigma$ are Borel measurable and satisfy \eqref{eq:mu_sigma_integ_p},
		and the reward functions $\Phi$ and $L$ are upper semi-analytic 
		and satisfy $\Phi(t, \xb) = \Phi(t, \xb_{t \wedge \cdot})$, $L(t, \xb, u) = L(t, \xb_{t \wedge \cdot}, u)$, for all $(t,\xb, u) \in \R_+ \x \Om \x U$.
		Then both value functions $V_W$ and $V_R$ are also upper semi-analytic.
		Moreover, for any $(t, \xb) \in \Rb_+ \x \Om$ and $\Fbb$-stopping time $\taub: \Omb \longrightarrow [t, \infty]$,
		by denoting $\tau^{\alpha} := \taub(\pi^{\alpha}, X^{\alpha}, M^{\alpha})$,
		one has the dynamic programming principle:
		\b*
			V_W(t,\xb) &\!\!=\!\!& \sup_{\alpha \in \Ac_W(t,\xb) } \E^{\P^{\alpha}} \Big[
			\Big(\int_t^{\pi^{\alpha}} \!\!\! L(s, X^{\alpha}, \nu^{\alpha}_s) ds + \Phi(\pi^{\alpha}, X^{\alpha})  \Big) \1_{\pi^{\alpha} \le \tau^{\alpha}} \\
			&&~~~~~~~~~~~~~~~~~~~~~+~
			\Big(\int_t^{\tau^{\alpha}} \!\!\! L(s, X^{\alpha}, \nu^{\alpha}_s) ds + V_W(\tau^{\alpha}, X^{\alpha}) \Big) \1_{\pi^{\alpha} > \tau^{\alpha}}
			\Big],
		\e*
		and
		\b*
			V_R(t,\xb) &\!\!=\!\!& \sup_{\alpha \in \Ac_R(t,\xb) } \E^{\P^{\alpha}} \Big[
			\Big(\int_t^{\pi^{\alpha}} \!\!\! \int_U L(s, X^{\alpha}, u) M^{\alpha}_s(du) ds + \Phi(\pi^{\alpha}, X^{\alpha})  \Big) \1_{\pi^{\alpha} \le \tau^{\alpha}} \\
			&&~~~~~~~~~~~~~~~~~+
			\Big(\int_t^{\tau^{\alpha}} \!\!\! \int_U L(s, X^{\alpha}, u) M^{\alpha}_s(du) ds + V_R(\tau^{\alpha}, X^{\alpha}) \Big) \1_{\pi^{\alpha} > \tau^{\alpha}}
			\Big].
		\e*
	\end{Theorem}
	\proof We only prove the results for the weak formulation.
	Let us consider the probability measures on  the canonical space $\Omb := \Rb_+ \x \Om \x \M$ induced by the weak controls:
	for all $(t, \xb) \in \Rb_+ \x \Om$,
	$$
		\Pcb_W(t,\xb) := \big\{ \P^{\alpha} \circ (\pi^{\alpha}, X^{\alpha}, \delta_{\nu^{\alpha}_s}(du)ds )^{-1}~:\alpha \in \Ac_W(t,\xb) \big\}.
	$$
	By Proposition \ref{prop:weak_reformulation}, 
	one notices that $\Pcb_W(t,\xb)$ is equal to the collection of all weak control/stopping rules (in the sense of Definition \ref{def:weak_relaxed_control})
	associated with generator $\Gh$ and initial condition $(t,\xb)$,
	where $\Gh := \big\{ (\varphi, \Lc^{t, \xb, u}\varphi)~: \varphi \in C^2_b(E) \big\}$ with $E := \R^d$ and
	\be \label{eq:diff_generator}
		\Lc^{t,\xb,u} \varphi(x) := \mu(t, \xb_{t \wedge \cdot}, u) \cdot D \varphi(x)
		+
		\frac{1}{2} \sigma\sigma^T(t, \xb_{t \wedge \cdot}, u) : D^2 \varphi(x),
		~~\mbox{for all}~x \in \R^d.
	\ee
	Further, by considering a countable dense subset of $C^2_b(\R^d)$
	(under the point-wise convergence of $\varphi$, $D \varphi$ and $D^2 \varphi$),
	it is clear that $\Gh$ is countably generated.
	One can then directly apply Theorem \ref{theo:ctrl_stopping_DPP} to conclude the proof.
	\qed

	\begin{Remark}
	{\rm	When $(\mu, \sigma)(t, \xb, u)$ is continuous in $\xb$, then using classical localization technique and compactness arguments, it can be deduced that $\Pcb^R_{t,\xb}$ is non-empty for every $(t, \xb) \in \Rb_+ \x \Om$
		(see e.g. Stroock and Varadhan \cite{StroockVaradhan}).
	}
	\end{Remark}

\subsubsection{The strong formulation}
\label{subsubsec:strong_form}

	We now consider the strong formulation of the controlled/stopped diffusion processes problem (see Section \ref{subsec:diffusion}),
	which needs a little more work to be reformulated as the general framework in Section \ref{subsec:DPP}.
	
	\vspace{0.5em}
	
	Recall that when $E = \R^d$, we also denote by $B$ the canonical process on $\Om_0 = \D(\R_+, \R^d)$ with 
	canonical filtration $\F$, et denote by $\P_0$ the Wiener measure under which $B$ is a standard Brownian motion.
	Let $\F^0 = (\Fc^0_t)_{t \ge 0}$ be the canonical filtration with $\Fc^0_t := \sigma(B_s ~:s \le t)$ and $\Fc^0_{\infty} := \sigma(B_s ~: s \ge 0)$,
	and $\F^a$ denote the augmented Brownian filtration on $\Om_0$ under $\P_0$.
	Further, let $\Uc$ denote the class of all control processes (i.e. all $U$-valued $\F^0$-predictable processes).
	For $t \in \Rb_+$ and $\om^0 \in \D(\R_+, \R^d)$, 
	we denote by $\Uc_t$ the subclass of all control processes independent of $\sigma(B_s~:s \le t)$ (under $\P_0$),
	and by $\P_0^{t, \om^0}$ the measure on $\Om_0$ under which 
	$\P_0^{t, \om^0} [ B_{t \wedge \cdot} = \om^0_{t \wedge \cdot}] = 1$ 
	and $(B_s - B_t)_{s \ge t}$ is a standard Brownian motion.

	\vspace{0.5em}

	In additional to the integrability condition \eqref{eq:mu_sigma_integ_p}, let us assume the following Lipschitz condition.
	\begin{Assumption} \label{assu:Lip_coef}
		For any $T > 0$, there is some constant $L_0 > 0$ such that,
		for all $(t, \om, \om', u) \in [0,T] \x \Om \x \Om \x U$,
		$$
			|\mu(t, \om, u) - \mu(t, \om', u)| ~+~ \|\sigma(t, \om, u) - \sigma(t, \om', u) \| 
			~\le~
			 L_0 \|\om - \om'\|_T,
		$$
		where $\| \om\|_T := \sup_{0 \le t\le T} |\om_t|$.
	\end{Assumption}	
	Then given a control $\nu \in \Uc$ and an initial condition $(t, \xb) \in \R_+ \x \Om$, the controlled SDE
	\begin{equation} \label{eq:controlled_SDE}
		X^{t,\xb,\nu}_s ~=~
		\xb_t + \int_t^s \mu(r, X^{t,\xb,\nu}_{r \wedge \cdot}, \nu_r) dr 
		+  \int_t^s \sigma(r, X^{t,\xb,\nu}_{r \wedge \cdot}, \nu_r) d B_r,
		~~\P_0 \mbox{-a.s.},
	\end{equation}
	with initial condition $X^{t,\xb,\nu}_s := \xb_s$ for all $s \in [0,t]$,
	has a unique strong solution (under  \eqref{eq:mu_sigma_integ_p} and Assumption \ref{assu:Lip_coef}).
	The value function $V_S$ of the strong formulation of the optimal controlled/stopped diffusion processes problem is given by
	\be \label{eq:VS1}
		V_S(t, \xb)
		~:=~
		\sup_{\nu \in \Uc} ~\sup_{\pi \in \Tc_t}~
		\E \Big[ \int_t^{\pi} L(s, X^{t, \xb, \nu}, \nu_s) ds + \Phi \big(\pi, X^{t,\xb,\nu}_{\cdot} \big) \Big],
	\ee
	where $\Tc_t$ denotes the collection of all $\F^a$-stopping times taking value in $[t, \infty]$.
	
	\vspace{0.5em}
	
	To study the above strong formulation in the framework of Section \ref{subsec:DPP},
	we need to consider an enlarged canonical space $\Omt := \Om_0 \x \Omb$  with $\Om_0 := \Om$.
	Let $(B, \Theta, X, M)$ be the canonical process on $\Omt$, defined by
	$B_t(\omt) := \om^0_t$, $X_t(\omt) := \om_t$, $\Theta_{\infty}(\omt) := \theta$ and $M(\omt) := m$,
	for all $t \in \R_+$ and $\omt = (\om^0, \theta, \om, m)$.
	Let $\widetilde \F = (\Fct_t)_{t \ge 0}$ denote the canonical filtration, defined by
	$\Fct_t := \sigma\big( B_s, X_s, M_s(\phi), \{\Theta_{\infty} \le s\}, ~ s \le t, \phi \in \C_b(\R_+ \x U) \big)$
	and $\Fct_{\infty} := \bigvee_{t \ge 0} \Fct_t$.
	Given a $\Ft$-stopping time  $\tilde \tau$, then for every $(t,\xb) \in \R_+ \x \Om$, $\pi \in \Tc$ and $\nu \in \Uc$,
	we define a $\F^a$-stopping time by
	\begin{equation} \label{eq:def_tau_tildetau}
		\tau^{\nu,\pi} (\om) ~:=~  \tilde \tau \big(\om, \pi(\om), X_{\cdot}^{t,\xb,\nu}(\om), \delta_{\nu_s(\om)}(du) ds \big).
	\end{equation}
	Our main DPP result of the strong formulation of the optimal controlled/stopped diffusion process is given as follows.

	\begin{Theorem} \label{theo:DPP_strong_formu}
		Assume that the coefficient functions $\mu$ and $\sigma$ satisfy Assumption \ref{assu:Lip_coef} and  \eqref{eq:mu_sigma_integ_p},
		and the reward functions $L$ and $\Phi$ are upper semi-analytic and satisfy
		$L(t, \xb, u) = L(t, \xb_{t \wedge \cdot}, u)$ and $\Phi(t, \xb) = \Phi(t, \xb_{t\wedge \cdot})$ for all $(t,\xb, u) \in \R_+ \x \Om \x U$.
		Then the value function $V_S: \Rb_+ \x \Om \longrightarrow \Rb$ defined in \eqref{eq:VS1} is also upper semi-analytic,
		and for every $(t, \xb) \in \Rb_+ \x \Om$ and $\Ft$-stopping time $\tilde \tau$ larger than $t$, together with the induced stopping times $(\tau^{\nu, \pi})$ in \eqref{eq:def_tau_tildetau},
		one has
		\b*
			V_S(t, \xb)
			&\!\!\!=\!\!\!&
			\sup_{\nu \in \Uc} ~\sup_{\pi \in \Tc_t}
			\E ~\Big[ \Big( \int_t^{\pi} L(s, X^{t,\xb,\nu}, \nu_s) ds + \Phi \big(\pi, X^{t,\xb,\nu} \big) \Big) \1_{\pi \le \tau^{\nu,\pi}} \\
			&&~~~~~~~~~~~~~~~~~~+
				\Big( \int_t^{\tau^{\nu,\pi}} \!\!\!\! L(s, X^{t,\xb,\nu}, \nu_s) ds + V_S \big(\tau^{\nu,\pi}, X^{t,\xb,\nu} \big) \Big) \1_{\pi > \tau^{\nu,\pi}} \Big].
		\e*
	\end{Theorem}

	To prepare the proof of Theorem \ref{theo:DPP_strong_formu}, we will reformulation the strong formulation \eqref{eq:VS1} of the optimal controlled/stopped diffusion processes problem
	on the enlarged canonical space $\Omt$ as a controlled/stopped martingale problem.
	With the given coefficient functions $\mu$ and $\sigma$, 
	let us define two coefficient functions $\mut: \R_+ \x \Om \x U \to \R^{2d}$ and $\sigmat: \R_+ \x \Om \x U \to \S^d \x \S^d$ by
	\b*
		\mut(t,\xb,u) :=
		\begin{pmatrix}
			0\\
			\mu(t,\xb,u)
		\end{pmatrix}
		&\mbox{and}&
		\sigmat(t,\xb,u) :=
		\begin{pmatrix}
			\text{Id}_d\\
			\sigma(t,\xb,u)
		\end{pmatrix}.
	\e*
	The full generator of the control problem is then given by
	$$
		\widetilde{\G} := \big\{ (\varphi, \tilde \Lc \varphi) ~: \varphi \in C^2_b(\R^{d} \x \R^d) \big\},
	$$
	where $\widetilde \Lc$ is the infinitesimal generator defined by, for all $ \varphi \in C^2_b(\R^d \x \R^d)$,
	\b*
		~\tilde \Lc \varphi(t,\xb,u, x)
		:=
		\mut(t,\xb,u) \cdot D \varphi(x)
		+
		\frac{1}{2} \sigmat \sigmat^T(t,\xb,u) : D^2 \varphi(x).
	\e*
	Similar to the case of generator $\widehat \G$, it is easy to see that $\widetilde \G$ is also countably generated.
    	We next equip $\Uc$ with the following $H_2$-norm $\| \cdot \|_{H_2}$ by
	\b*
		\| \nu^1 - \nu^2 \|_{H_2}^2
		&:=&
		\E \Big[ \int_0^{\infty} e^{-\beta t} \big( d(\nu^1_t, \nu^2_t) \big)^2 dt \Big],
		~~\mbox{for some constant}~ \beta > 0,
	\e*
	so that $\Uc$ is a Polish space.
	Further, every $\nu \in \Uc$ induces a probability measure on $\Om_0 \x \M$ by
	\b*
		\Pi(\nu)
		&:=&
		\P_0 \circ ( B, m^{\nu})^{-1},
		~~\mbox{where}~
		m^{\nu} := \delta_{\nu_t(B)}(du) dt \in \M, ~\mbox{for all}~ \nu \in \Uc.
	\e*
	Notice that the operator $\Pi : \Uc \longrightarrow \Pc(\Om_0 \x \M)$ is continuous and injective,
	it follows that $\Pi(\Uc):= \{ \Pi(\nu) ~: \nu \in \Uc\}$ is a Borel set in the Polish space $\Pc(\Om_0 \x \M)$.
	We shall also consider the set $\Pi(\Uc_t) := \{ \Pi(\nu) ~: \nu \in \Uc_t\}$.
	
	\begin{Remark} \label{rem:SDE_t_omega}
		The controlled SDE \eqref{eq:controlled_SDE} is driven by the increment of the Brownian motion $B$ after time $t$.
		For a control process $\nu \in \Uc_t$ independent of $\sigma(B_s, s \le t)$, the solution of \eqref{eq:controlled_SDE} under $\P_0$
		and that under $\P_0^{t, \om^0}$ have the same distribution.
		In this setting, by abus of notation, we always denote it by $X^{t, \xb, \nu}$.
	\end{Remark}

	With the above preparation, we can then reformulate the control/stopping problem \eqref{eq:VS1} on $\Omt$ as a controlled martingale problem.
	For every $(t,\xb, \om^0) \in \R_+ \x \Om \x \Om_0$, let
	$$
		\Pct_S(t,\xb) 
		~:=~
		\big\{ \P_0 \circ (B, \pi, X^{t, \xb, \nu}, m^{\nu})^{-1} ~: \pi \in \Tc_t, ~\nu \in \Uc \big\},
	$$
	$$
		\Pct_S (t,\xb, \om^0)
		~:=~
		\big\{ \P^{t, \om^0}_0 \circ (B, \pi, X^{t, \xb, \nu}, m^{\nu})^{-1} ~: \pi \in \Tc_t, ~\nu \in \Uc_t \big\},
	$$
	and
	\b*
		\Pct (t,\xb)
		&\!\!\!\!:=\!\!\!\!&
		\Big\{
		\Pt \in \Pc(\Omt) ~: \Pt \big[ \Theta_{\infty} \ge t, X_s = \xb_s, ~0 \le s \le t \big] = 1, ~\Pt|_{\Om_0 \x \M} \in \Pi(\Uc_t),\\
		&& ~~~~~~~~~\big( C_s^n(f,g)\big)_{s \ge t} ~\mbox{is a}~ (\Pt, \widetilde{\F}) \mbox{-martingale}, ~\mbox{for all}~ n \ge 1~\mbox{and}~ (f,g) \in \widetilde{\G}
		\Big\},
	\e*
	\b*
		\Pct (t,\xb, \om^0)
		&\!\!\!\!:=\!\!\!\!&
		\Big\{
		\Pt \in \Pc(\Omt) ~: \Pt \big[ \Theta_{\infty} \ge t, X_s = \xb_s, B_s = \om^0_s, ~0 \le s \le t \big] = 1,\\
		&& ~~~~~~~~~~~~~~~~~~ \Pt \big[ M(du, ds) = \delta_{\nu_s(B)}(du) ds \big] = 1, ~\mbox{for some}~\nu \in \Uc_t,\\
		&& ~~~~~~ \big( C_s^n(f,g)\big)_{s \ge t} ~\mbox{is a}~ (\Pt, \widetilde{\F}) \mbox{-martingale}, ~\mbox{for all}~ n \ge 1~\mbox{and}~ (f,g) \in \widetilde{\G}
		\Big\}.
	\e*
	When $t = \infty$, we let
	$$
		\Pct (\infty, \xb) := \big\{ \Pt \in \Pc(\Omt) ~: \Pt \big[ \Theta_{\infty} = \infty, X = \xb \big] = 1, ~\Pt|_{\Om_0 \x \M} \in \Pi_{\infty}(\Uc) \big\},
	$$
	
	Namely, $\Pct_S(t, \xb)$ is the set of control/stopping rules induced by the control with all possible control processes $\nu \in \Uc$,
	and $\Pct_S(t, \xb, \om^0)$ is induced by those with control processes $\nu \in \Uc_t$ (i.e. independent of the Brownian motion before time $t$).
	We observe that, the canonical variable $\Theta_{\infty}$ is a stopping time w.r.t. the canonical filtration $\widetilde{\F}$.
	However, while it is still a stopping time w.r.t. the augmented Brownian filtration under a control/stopping rule in $\Pct_S(t, \xb)$ and $\Pct_S (t,\xb, \om^0)$,
	it may not be so under a control/stopping rule in $\Pct(t, \xb)$ and $\Pct (t,\xb, \om^0)$.
	Thus a control/stopping rule in  $\Pct(t, \xb)$ (resp. $\Pct(t, \xb, \om^0)$) may not be a rule in $\Pct_S(t,\xb)$ (resp. $\Pct_S(t, \xb, \om^0)$).
	
	\vspace{0.5em}
	
	Let us define the value functions $\Vt$ and $\Vt_S$ by
	$$
		\Vt (t, \xb)
		:=\!\!
		\sup_{\Pt \in \Pct(t,\xb)} \!\!\!\! J(\Pt),
		~~~
		\Vt (t, \xb, \om^0)
		:=\!\!
		\sup_{\Pt \in \Pct (t,\xb, \om^0)} \!\!\!\! J(\Pt),
		~~~
		\Vt_S (t, \xb, \om^0)
		:=\!\!
		\sup_{\Pt \in \Pct_S(t,\xb, \om^0)} \!\!\!\! J(\Pt),
	$$
	where
	$$
		J(\Pt)
		~:=~
		\E^{\Pt} \Big[  \int_t^{\Theta_{\infty}} \!\!\! \int_U L(s, X,u) M_s(du) ds + \Phi \big(\Theta_{\infty}, X \big) \Big].
	$$

	\begin{Lemma} \label{lemm:equiv_dpp}
		Let us stay in the setting of Theorem \ref{theo:DPP_strong_formu}.
		
		\vspace{0.5em}
		
		\noindent \rmi For all $(t,\xb) \in \R_+ \x \Om$ and $\om^0 \in \Om_0$, one has
		\be \label{eq:VS_equiv}
			\Vt (t, \xb)
			~=~
			\Vt (t, \xb, \om^0)
			~=~
			\Vt_S (t, \xb, \om^0)
			~=~
			V_S(t, \xb)
			= \sup_{\Pt \in \Pct_S(t, \xb)} J(\Pt).
		\ee
		\rmii Moreover, the graph set $[[\Pct]]$ of the family $(\Pct(t, \xb))_{(t, \xb) \in \Rb_+ \x \Om}$ is Borel,
		so that $V_S: \Rb_+ \x \Om \longrightarrow \Rb$ is upper semi-analytic.
		Further, for all $ \Ft$-stopping time $\tilde \tau$ taking value in $[t, \infty]$, one has the DPP
		\be \label{eq:DPP_Vt}
			V_S(t, \xb)
			&\!\!=\!\!\!\!&
			\sup_{\Pt \in \Pct_S(t,\xb)}\!\!\!
			\E^{\Pt} \Big[ \Big( \int_t^{\Theta_{\infty}} \!\!\! \int_U L(s, X,u) M_s(du) ds + \Phi \big(\Theta_{\infty}, X \big) \Big)
				\1_{\Theta_{\infty} \le \tilde\tau} \nonumber \\
			&&~~~~~~~~~~~~~~~~+
				\Big( \int_t^{\tilde \tau} \!\! \int_U L(s, X,u) M_s(du) ds + V_S(\tilde \tau, X) \Big) \1_{ \Theta_{\infty} > \tilde \tau } 
			\Big].~~~
		\ee
	\end{Lemma}
	\proof \rmi First, we notice that $\Vt_S(t, \xb, \om^0)$ clearly does not depend on $\om^0 \in \Om_0$,
	and in view of Remark \ref{rem:SDE_t_omega}, 
	a control/stopping rule in $\Pct_S(t, \xb, \om^0)$ can be considered as special control/stopping rule in $\Pct_S(t, \xb)$ which depends only on the increment of the Brownian motion after time $t$.
	Therefore, one has 
	$$
		\Vt_S(t, \xb, \om^0) ~\le~ V_S(t, \xb) ~=~ \sup_{\Pt \in \Pct_S(t, \xb)} J(\Pt).
	$$
	On the other hand, given $\Pt \in \Pct_S(t, \xb)$, 
	its (regular) conditional probability $(\Pt_{\om^0})_{\om^0 \in \Om_0}$ knowing $\sigma(B_s~:s \le t)$ satisfies $\Pt_{\om^0} \in \Pct_S(t, \xb, \om^0)$ for $\Pt$-a.e. $\om^0 $ (see also \cite[section 2]{JDX} for a more detailed argument).
	It follows that $J(\Pt) \le \E \big[ \Vt_S(t, \xb, B) \big] = \Vt_S(t, \xb, \om^0)$ for an arbitrary $\om^0$.
	This proves that $V_S(t, \xb) = \Vt_S(t, \xb, \om^0)$.
	By exactly the same argument,
	one can prove that $\Vt(t, \xb) = \Vt(t, \xb, \om^0)$.
	
	\vspace{0.5em}

	\noindent Next, we observe that $\Pct_S(t, \xb, \om^0) \subset \Pct(t, \xb, \om^0)$ so that $\Vt_S(t, \xb, \om^0) \le \Vt(t, \xb, \om^0)$.
	However, let $\Pt \in \Pct(t,\xb, \om^0)$, the control/stopping rule $\Pt$ is not necessarily in $\Pct_S(t, \xb, \om^0)$ since $\Theta_{\infty}$ is not necessarily a stopping time w.r.t. the augmented filtration $\widetilde{\F}^{B,\Pt}_+$ generated by the Brownian motion $B$.
	Nevertheless, $B$ is a $(\Pt,\widetilde{\F})$-Brownian motion, 
	and there is some control process $\nu$ which is $\widetilde {\F}^B$-predictable such that $M(ds, du) = \delta_{\nu_s(B)}(du) ds, ~\Pt$-a.s., with $\Ft^B := (\Fct^B_t)_{t \ge 0}$.
	Then as strong solution to the controlled SDE, $X$ is continuous and $\widetilde{\F}^{B,\Pt}_+$-adapted.
	Moreover, denoting by $\Ft_+$ the right-continuous version of filtration $\Ft$, 
	$B$ is a $\widetilde{\F}_+$-Brownian motion and
	$\Theta_{\infty}$ is a $\widetilde{\F}_+$-stopping time.
	Notice that the filtered space $(\Omt, \Fct_{\infty}, \Pt, \widetilde{\F}_+)$ together with the Brownian motion $B$ satisfies property (K) in the optimal stopping theory,
	then it follows by Proposition \ref{prop:eqiv_opt_stopping} (see also Remark \ref{rem:PropertyK}) that 
	$J(\Pt) \le \Vt_S (t, \xb, \om^0)$
	(see more details about property (K) and the equivalence of optimal stopping problem in Section \ref{subsec:equiv_opt_stop}).
	This proves that $\Vt(t, \xb, \om^0) = \Vt_S(t, \xb, \om^0)$.

	\vspace{0.5em}

	\noindent \rmii For the second part of the statement, let us first consider the graph set
	$$
		[[\Pct]] 
		~:=~
		\big \{(t,\xb, \Pt) ~: (t,\xb) \in \Rb_+ \x \Om, ~\Pt \in \Pct(t,\xb) \big\}.
	$$
	In view of Lemma \ref{lemm:ctrl_stopping_measurability}, in order to prove that $[[\Pct]]$ is Borel measurable,
	it is enough to prove the Borel measurability of 
	$$
		A ~:=~ \big\{ (t, \xb, \Pt) ~: (t, \xb) \in \R_+ \x \Om, ~ \Pt|_{\Om_0 \x \M} \in \Pi(\Uc_t) \big\}.
	$$
	Notice that $\Pt|_{\Om_0 \x \M} \in \Pi(\Uc_t)$ is equivalent to $\Pt|_{\Om_0 \x \M} \in \Pi(\Uc)$ and $M_{\cdot}(\phi)$ is independent of $B_{t \wedge \cdot}$ under $\Pt$ for all $\phi \in C_b(\R_+ \x U)$.
	Therefore, there exists a countable family of (bounded continuous) test functions $(\phi_n, \varphi_n, \psi_n)_{ n \ge 1}$ rich enough, such that
	\b*
		A &\!\!\!=\!\!\!&  \Big\{ 
			(t, \xb, \Pt) ~:  (t,\xb) \in \Rb_+ \x \Om, ~\Pt|_{\Om_0 \x \M} \in \Pi(\Uc), \\
		&&~~~~~~\mbox{and}~
			\E^{\Pt} \big[ \varphi_n( M_{\cdot}(\phi_n)) \psi_n(B_{t \wedge \cdot}) \big] 
			=
			\E^{\Pt} \big[ \varphi_n( M_{\cdot}(\phi_n))\big] \E^{\Pt} \big[ \psi_n(B_{t \wedge \cdot}) \big] , ~n \ge 1
		\Big\}.
	\e*
	Notice that $\Pi(\Uc)$ is a Borel set, this is enough to prove that $A$ is Borel, and hence $[[\Pct]]$ is Borel measurable.
	
	\vspace{0.5em}
	
	Next, to prove the dynamic programming result in \eqref{eq:DPP_Vt}, we follow Theorem \ref{theo:DPP} to 
	to apply the conditioning and concatenation arguments.
	First, for an arbitrary $\Pt \in \Pct_S(t, \xb)$, let $\tilde \tau$ be a  $\Ft$-stopping time taking value in $[t, \infty]$,
	we consider a family of r.c.p.d. $(\Pt_{\omt})_{\omt \in \Omt}$ of $\Pt$ knowing $\Fct_{\tilde \tau \wedge \Theta_{\infty}}$.
	One can check that (see in particular  Claisse, Talay and Tan \cite{JDX} for more detailed arguments) 
	$\Pt_{\omt} \in  \Pct_S(\tau(\omt), \omt^X, \omt^B)$, for $\Pt$-a.e. $\omt = (\omt^B, \omt^{\Theta}, \omt^X, \omt^M) \in \Omt$ such that $\tilde \tau(\omt) < \Theta_{\infty}(\omt)$.
	Together with arbitrariness of $\Pt \in \Pct(t, \xb)$ and the fact that $\Vt_S(\tau(\omt), \omt^X, \omt^B) = V_S(\tau(\omt), \omt^X)$, 
	this proves the claim, which implies that
	\b*
		V_S(t, \xb)
		&\!\!\le\!\!\!\!&
		\sup_{\Pt \in \Pct_S(t,\xb)}\!\!\!
		\E^{\Pt} \Big[ \Big( \int_t^{\Theta_{\infty}} \!\!\! \int_U L(s, X,u) M_s(du) ds + \Phi \big(\Theta_{\infty}, X \big) \Big)
			\1_{\Theta_{\infty} \le \tilde\tau} \nonumber \\
		&&~~~~~~~~~~~~~~~~~~~~~+~
			\Big( \int_t^{\tilde \tau} \! \int_U L(s, X,u) M_s(du) ds + V_S(\tilde \tau, X) \Big) \1_{ \Theta_{\infty} > \tilde \tau } 
		\Big].
	\e*
	
	To prove the inverse inequality, we follow Theorem \ref{theo:DPP} to use the concatenation arguments.
	First, let $\Pt \in \Pct_S(t, \xb)$. In view of the equivalence result in Item $\mathrm{(i)}$, one can assume w.l.o.g. that
	$\Pt[ M(du, ds) = \delta_{\nu_s(B)}(du) ds ] = 1$ for some $\nu \in \Uc_t$.
	Let us denote by $\Ft^B = (\Fct^B_t)_{t \ge 0}$ the filtration generated by $B$ on $\Omt$, 
	and by $\Ft^{B,a}$ the augmented (Brownian) filtration under $\Pt$.
	In particular, $\Theta_{\infty}$ is a $\Ft^{B,a}$-stopping time, and $(X, M)$ are $\Ft^{B,a}$-adapted.
	Thus the $\Ft$-stopping times $\tilde \tau$ taking value in $[t, \infty]$ is also a $\Ft^{B,a}$-stopping time.
	Then there exits a $\Ft^B$-stopping time $\tilde \tau'$ on $\Omt$, such that $\Pt[\tilde \tau' = \Theta_{\infty} \wedge \tilde \tau] = 1$.
	Moreover, a family of r.c.p.d. of $\Pt$ knowing $\Fct^B_{\tilde \tau'}$ is also a family of r.c.p.d. of $\Pt$ knowing $\Fct_{\tilde \tau'}$,
	since for any bounded r.v. $\xi$, one has
	$\E^{\Pt} [ \xi | \Fct_{\tilde \tau'} ] = \E^{\Pt} [ \xi | B_{\tilde \tau' \wedge} ] = \E^{\Pt} [ \xi | \Fct^B_{\tilde \tau'} ]$, $\Pt$-a.s.
	
	\vspace{0.5em}
	
	Next, as in Theorem \ref{theo:DPP}, we apply the measurable selection theorem to choose a (universally) measurable family of $(\Pt^{\eps}_{s, \xb})_{(s, \xb) \in \Rb_+ \x \Om}$ such that each $\Pt^{\eps}_{s, \xb}$ consists in a $\eps$-optimizer in $\Pct(s, \xb)$ for the optimization problem in the definition of $\Vt(s, \xb)$.
	Let us further define
	$$
		\Qt^{\eps}_{s, \xb, \om^0} ~:=~ \Pt^{\eps}_{s, \xb} \circ \big( \delta_{\om^0} \otimes_s B, \Theta_{\infty}, X, M \big)^{-1},
		~~\mbox{for all}~
		(s, \xb, \om^0) \in \Rb_+ \x \Om \x \Om_0,
	$$
	where $(\delta_{\om^0} \otimes_s B )_r := \om^0_r \1_{\{r \in [0,s]\}} + (\om^0_s + B_r - B_s ) \1_{\{r \in (r, \infty)\}}$.
	One observes that $(s, \xb, \om^0) \longmapsto \Qt^{\eps}_{s, \xb, \om^0} $ is still universally measurable
	and $\Qt^{\eps}_{s, \xb, \om^0} \in \Pct(s, \xb, \om^0)$.
	Moreover, $\Qt^{\eps}_{s, \xb, \om^0}$ is also $\eps$-optimizer in $\Pct(s, \xb, \om^0)$ for the optimization problem in the definition of $\Vt(s, \xb, \om^0)$.
	
	\vspace{0.5em}
	
	Now, let $\Qt_{\omt} := \Qt^{\eps}_{\tilde \tau'(\omt), \omt^X, \omt^B}$ for all $\omt = (\omt^B, \omt^{\Theta}, \omt^X, \omt^M) \in \Omt$.
	We consider  the concatenated probability measure $\Pt \otimes_{\tilde \tau'} \Qt_{\cdot}$
	and claim that  $\Pt \otimes_{\tilde \tau'} \Qt_{\cdot} \in \Pct(t, \xb)$.
	By similar arguments as in Lemma \ref{lemm:ctrl_stopping_stability}, 
	it is easy to see that $\Pt \otimes_{\tilde \tau'} \Qt_{\cdot}$ solves the corresponding martingale problem in the definition of $\Pct(t, \xb)$,
	and $B$ is still a Brownian motion under $\Pt \otimes_{\tilde \tau'} \Qt_{\cdot}$.
	Then, it is enough to prove that $\Pt \otimes_{\tilde \tau'} \Qt_{\cdot} |_{\Om_0 \x \M} \in \Pi(\Uc_t)$, or equivalently that
	$$
		\Pt \otimes_{\tilde \tau'} \Qt_{\cdot} 
		\big[M(du, ds) = \delta_{\nu_s(B)}(du) ds \big] 
		=1,
		~~\mbox{for some}~\nu \in \Uc_t.
	$$
	Since $\Pt \otimes_{\tilde \tau'} \Qt_{\cdot}  \big[M \in \M_0 \big] = 1$, 
	this reduces to prove that, with $\Fct^{B}_{\infty} := \sigma(B_s ~s \ge 0)$,
	\begin{equation} \label{eq:M_B_measurable}
		\E^{\Pt \otimes_{\tilde \tau'} \Qt_{\cdot}} \big[ M_s(\phi) ~\big|~ \Fct^{B}_{\infty} \big] =M_s(\phi),
		~\Pt \otimes_{\tilde \tau'} \Qt_{\cdot} \mbox{-a.s.}
		~~\mbox{for all}~
		s \in \R_+,~
		\phi \in C_b(\R_+ \x U).
	\end{equation}
	Notice that $\Pt \in \Pct_S(t, \xb)$ and $\Qt_{\omt} \in \Pct({\tilde \tau'(\omt), \omt^X, \omt^B})$, then
	$$
		\E^{\Pt} \big[ M_s(\phi) \1_{\{s \le \tilde \tau'\}} ~\big|~ \Fct^{B}_{\infty} \big] =M_s(\phi) \1_{\{s \le \tilde \tau'\}},~\Pt\mbox{-a.s.}
	$$
	and
	$$
		\E^{\Qt_{\omt}}\big[ M_s(\phi) \1_{\{s > \tilde \tau'(\omt)\}}~\big|~ \Fct^{B}_{\infty} \big] =M_s(\phi) \1_{\{s > \tilde \tau'(\omt)\}}, ~\Qt_{\omt} \mbox{-a.s. for each}~\omt \in \Omt.
	$$
	Moreover, since $(\Qt_{\omt})_{\omt \in \Omt}$ is also a r.c.p.d. of $\Pt \otimes_{\tilde \tau'} \Qt_{\cdot}$ knowing $\Fct^B_{\tilde \tau'}$,
	it follows by Lemma \ref{lemm:conditioning_conditional_proba} below that
	$$
		\E^{\Pt \otimes_{\tilde \tau'} \Qt_{\cdot}} \big[ M_s(\phi) ~\big|~ \Fct^{B}_{\infty} \big]
		=
		\E^{\Qt_{\omt}}\big[ M_s(\phi) ~\big|~ \Fct^{B}_{\infty} \big],
		~\Qt_{\omt}\mbox{-a.s.}
		~\mbox{for}~\Pt\mbox{-a.e.}~\omt.
	$$
	This is enough to prove \eqref{eq:M_B_measurable}, and hence  the claim that $\Pt \otimes_{\tilde \tau'} \Qt_{\cdot} \in \Pct(t, \xb)$ holds true.
	When $V_S$ is finite, one can then argue as in Theorem \ref{theo:DPP} to conclude that 
	\b*
		V_S(t, \xb)
		&\!\!\ge\!\!\!\!&
		\E^{\Pt} \Big[ \int_t^{\tilde \tau'} \!\!\! \int_U L(s, X,u) M_s(du) ds + V_S \big(\tilde \tau', X \big) - \eps \Big]\\
		&\!\!\ge\!\!\!\!&
		\E^{\Pt} \Big[ \Big( \int_t^{\Theta_{\infty}} \!\!\! \int_U L(s, X,u) M_s(du) ds + \Phi \big(\Theta_{\infty}, X \big) \Big)
			\1_{\Theta_{\infty} = \tilde\tau' \le \tilde\tau} \nonumber \\
		&&~~~~~~~~~~~~~+
			\Big( \int_t^{\tilde \tau} \!\! \int_U L(s, X,u) M_s(du) ds + V_S(\tilde \tau, X) \Big) \1_{ \Theta_{\infty} > \tilde \tau' = \tilde \tau } 
		\Big] - \eps.
	\e*
	so that \eqref{eq:DPP_Vt} holds by arbitrariness of $\Pt$ and $\eps > 0$.
	When $V_S$ takes possibly the value $\infty$ or $-\infty$,
	one can still proceed as in  Theorem \ref{theo:DPP} to conclude.
	\qed

	\begin{Lemma} \label{lemm:conditioning_conditional_proba}
		Let $(\Omt, \Fct, \Pt)$ be a probability space, equipped with two sub-$\sigma$-filed $\Fct_1 \subset \Fct_2$.
		Assume that $\Fct$, $\Fct_1$ and $\Fct_2$ are all countably generated,
		and $(\Pt^1_{\omt})_{\omt \in \Omt}$ be a family of r.c.p.d. of $\Pt$ knowing $\Fct_1$,
		and $(\Pt^2_{\tilde \w})_{\tilde \w \in \Omt}$ be a family of r.c.p.d. of $\Pt$ knowing $\Fct_2$.
		The for $\Pt$-a.e. $\omt \in \Omt$, the family $(\Pt^2_{\tilde \w})_{\tilde \w \in \Omt}$ is a family of r.c.p.d. of $\Pt_{\omt}$ knowing $\Fct_2$.
	\end{Lemma}
	\proof First, for all bounded random variables $\xi \in \Fct$ and $\zeta \in \Fct_2$, one has by tower property that
		\begin{eqnarray*}
			\E^{\Pt^1_{\omt}} \big[ \xi \zeta \big]
			&=&
			\E^{\Pt} \big[\xi \zeta \big| \Fct_1 \big] (\omt)
			~=~
			\E^{\Pt} \big[  \E^{\Pt}[ \xi \zeta | \Fct_2 ]  \big| \Fct_1 \big](\omt) \\
			&=&
			\E^{\Pt^1_{\omt}} \big[  \E^{\Pt}[ \xi \zeta | \Fct_2 ]  \big] 
			~=~
			 \E^{\Pt^1_{\omt}} \big[ \E^{\Pt}[ \xi | \Fct_2 ] \zeta \big],
			~~
			\mbox{for}~\Pt\mbox{-a.e.}~\om.
		\end{eqnarray*}
		Therefore, for a sequence $(\xi_n, \zeta_n)_{n \ge 1}$ rich enough, there exits $\Omt_1 \subset \Omt$,
		such that $\Pt[\Omt_1] = 1$,
		and for all $\omt \in \Omt_1$, one has
		$$
			\E^{\Pt^1_{\omt}} \big[ \xi_n \zeta_n \big]
			~=~
			\E^{\Pt^1_{\omt}} \big[ \E^{\Pt}[ \xi_n | \Fct_2 ] \zeta_n \big].
		$$
		When $(\zeta_n)_{n \ge 1}$ is rich enough, this implies that for all $\omt \in \Omt_1$ and $\xi_n$,
		$$
			\E^{\Pt^1_{\omt}} \big[ \xi_n  \big| \Fct_2 \big] (\tilde \w)
			~=~
			\E^{\Pt} [\xi_n | \Fct_2 ] (\tilde \w),
			~~\Pt^1_{\omt} \mbox{-a.e.}~ \tilde \w \in \Omt.
		$$
		Finally, when $(\xi_n)_{n \ge 1}$ is rich enough,
		it follows that, for every $\omt \in \Omt_1$,
		$(\Pt^2_{\tilde\w})_{\tilde\w \in \Omt}$ is a family of r.c.p.d. of $\Pt^1_{\om}$ knowing $\Fct_2$.
	\qed

	\vspace{0.5em}

	\noindent {\sc Proof of Theorem \ref{theo:DPP_strong_formu}.}
	It is a direct consequence of Item \rmii of Lemma \ref{lemm:equiv_dpp}.
	\qed

\subsubsection{More examples of the stochastic control problems}

	With the above results for the optimal control/stopping problem, 
	by manipulating the reward function $\Phi : \Rb_+ \x \Om \to \Rb$,
	we can easily deduce the DPP for various different formulations of pure control problems.
	Throughout this section, let us stay in the context of Theorem \ref{theo:DPP_strong_formu},
	i.e. the coefficient functions $\mu$ and $\sigma$ satisfy Assumption \ref{assu:Lip_coef} and  \eqref{eq:mu_sigma_integ_p}.
	
	\vspace{0.5em}

	Let us fix $(t,\xb) \in \Rb_+ \x \Om$,
	$\tilde \tau$ be a $\Ft$-stopping time taking value in $[t,\infty]$, $\nu \in \Uc$, 
	we denote $\tau^{\nu}(\om) := \tilde \tau \big( \om, \infty, X^{t, \xb, \nu}_{\cdot}(\om), \delta_{\nu_s(\om)}(du)ds \big)$ which is a stopping time on $(\Om_0, \Fc^0_{\infty}, \P_0)$ w.r.t. the augmented Brownian filtration.

	\begin{Corollary}[A pure control problem]
		Let  $\Phi_1: \Om \longrightarrow \Rb$ and $L_1 : \R_+ \x \Om \x U \longrightarrow \Rb$ be upper semi-analytic 
		and satisfy $L_1(t, \xb, u) = L_1(t, \xb_{t \wedge \cdot}, u)$ for all $(t,\xb, u) \in \R_+ \x \Om \x U$.
		We consider the following control problem 
		$$
			V_1^S(t, \xb)
			~:=~
			\sup_{\nu \in \Uc}
			~\E \Big[
				\int_t^{\infty} L_1(s, X, \nu_s) ds
				+
				\Phi_1 \big( X^{t,\xb,\nu} \big)
			\Big].
		$$
		Then $V_1^S: \Rb_+ \x \Om \longrightarrow \Rb$ is upper semi-analytic, and one has the dynamic programming principle:
		$$
			V^S_1(t, \xb)
			~=~
			\sup_{\nu \in \Uc}
			~\E ~\Big[
				\int_t^{\tau^{\nu}}  L_1(s, X, \nu_s) ds
				+
				V^S_1 \big(\tau^{\nu}, X^{t,\xb,\nu} \big) 
			\Big].
		$$
	\end{Corollary}
	\proof It is enough to set $\Phi(\theta, \xb) := \Phi_1(\xb) \1_{\theta = \infty} - \infty \1_{\theta < \infty}$, and then apply Theorem \ref{theo:DPP_strong_formu} to conclude the proof.
	\qed

	\begin{Corollary}[A control problem with random horizon]
		Let  $\Phi_2: \Rb_+ \x \Om \longrightarrow \Rb$ and $L_2 : \R_+ \x \Om \x U \longrightarrow \Rb$ be upper semi-analytic 
		and satisfy $\Phi_2(t,\xb) = \Phi_2(t, \xb_{t \wedge \cdot})$ and $L_2(t, \xb, u) = L_2(t, \xb_{t \wedge \cdot}, u)$ for all $(t,\xb, u) \in \R_+ \x \Om \x U$.
		Let $E_0 \subset E$ be a closed subset of $E = \R^d$, and $\pi^{\nu} := \inf \{s ~: X^{t,\xb,\nu}_{s-} \in E_0 ~\mbox{or}~ X^{t,\xb,\nu}_s \in E_0 \}$.
		We consider the following control problem 
		$$
			V_2^S(t, \xb)
			~:=~
			\sup_{\nu \in \Uc}
			~\E \Big[
				\int_t^{\pi^{\nu}} L_2(s, X^{t, \xb, \nu}, \nu_s) ds 
				+
				\Phi_2 \big( \pi^{\nu},  X^{t,\xb,\nu} \big)
			\Big].
		$$
		Then $V_2^S: \Rb_+ \x \Om \longrightarrow \Rb$ is also upper semi-analytic, and one has the dynamic programming principle:
		$$
			V^S_2(t, \xb)
			~=~
			\sup_{\nu \in \Uc}
			~\E ~\Big[
				\int_t^{\tau^{\nu} \wedge \pi^{\nu}} L_2(s, X^{t, \xb, \nu}, \nu_s) ds 
				+
				 V^S_2 \big( \tau^{\nu} \wedge \pi^{\nu}, X^{t,\xb,\nu} \big) \Big].
		$$
	\end{Corollary}
	\proof Notice that $\pi(\om) := \inf \{s ~: \om_{s-} \in E_0 ~\mbox{or}~ \om_s \in E_0 \}$ defines a $\F$-stopping time,
	set $\Phi(\theta, \xb) := \Phi_2(\pi(\xb), \xb) \1_{\theta = \infty} - \infty \1_{\theta < \infty}$, and then use Theorem \ref{theo:DPP_strong_formu}, we hence conclude the proof.
	\qed

	\begin{Corollary}[A control problem under state constraint]
		Let  $\Om_0 \subset \Om$ be a Borel subset of $\Om$,
		$\Phi_3: \Om_0 \longrightarrow \Rb$ and $L_3 : \R_+ \x \Om_0 \x U \longrightarrow \Rb$ be upper semi-analytic 
		and satisfy $L_3(t, \xb, u) = L_3(t, \xb_{t \wedge \cdot}, u)$ for all $(t,\xb, u) \in \R_+ \x \Om_0 \x U$.
		Let $\Uc^{t,\xb}_0 \subset \Uc$ be a subset of control processes $\nu$ in $\Uc$, such that $X^{t,\xb, \nu} \in \Om_0$, $\P_0$-a.s.
		We consider the following control problem 
		$$
			V_3^S(t, \xb)
			~:=~
			\sup_{\nu \in \Uc^{t,\xb}_0}
			~\E \Big[
				\int_t^{\infty} L_3 (s, X, \nu_s) ds
				+
				\Phi_3 \big( X^{t,\xb,\nu} \big)
			\Big].
		$$
		Then $V_3^S: \Rb_+ \x \Om \longrightarrow \Rb$ is also upper semi-analytic, and one has the dynamic programming principle:
		$$
			V^S_3(t, \xb)
			~=~
			\sup_{\nu \in \Uc_0^{t,\xb}}
			~\E ~\Big[ 
				\int_t^{\tau^{\nu}} L_3(s, X, \nu_s) ds
				+
				V^S_3 \big(\tau^{\nu}, X^{t,\xb,\nu} \big) 
			\Big].
		$$
	\end{Corollary}
	\proof It is enough to set $\Phi(\theta, \xb) := \Phi_3(\xb) \1_{\theta = \infty, ~\xb \in \Om_0} - \infty \1_{\theta < \infty ~\mbox{or}~ \xb \notin \Om_0}$, and then apply Theorem \ref{theo:DPP_strong_formu} to conclude the proof.
	\qed

\subsubsection{Relaxation of the integrability condition \eqref{eq:mu_sigma_integ_p}}
\label{subsubsec:unbound_mu_sigma}

	In many situation, the integrability condition \eqref{eq:mu_sigma_integ_p} becomes a little restrictive for a controlled diffusion processes problem.
	In place of \eqref{eq:mu_sigma_integ_p}, let us consider the following technical conditions:
	for some constant $C> 0$ and $u_0 \in U$,
	\begin{equation} \label{eq:linear_growth_cond}
		|\mu(t, \xb, u| + \| \sigma(t, \xb, u) \|
		\le
		C \big( 1 +  \|\xb_{t \wedge \cdot} \| + d(u, u_0) \big),
		~\mbox{for all}~(t, \xb, u) \in \R_+ \x \Om \x U.
	\end{equation}
	At the same time, for the weak/relaxed formulation of the optimal control/stopping problem, 
	we recall the definition of $\Ac_W(t, \xb)$ and $\Ac_R(t, \xb)$ in Section \ref{subsubsec:weak_relax_form},
	and define
	$$
		\Ac^2_W(t,\xb) ~:=~ 
		\Big\{
			\alpha \in \Ac_W(t, \xb)~:
			\int_t^{\infty} |\nu^{\alpha}_s - u_0|^2 ds < \infty
		\Big\},
	$$
	and
	$$
		\Ac^2_R(t,\xb) ~:=~ 
		\Big\{
			\alpha \in \Ac_R(t, \xb)~:
			\int_t^{\infty} \int_U |u - u_0|^2 M^{\alpha}_s(du) ds < \infty
		\Big\}.
	$$
	One can then define the value function of the new weak/relaxed formulation of the problem:
	$$
		V^2_W(t,\xb) ~:= \sup_{\alpha \in \Ac^2_W(t,\xb) } \E^{\P^{\alpha}} \Big[
			\int_t^{\pi^{\alpha}} L(s, X^{\alpha}, \nu^{\alpha}_s) ds
			+
			\Phi( \pi^{\alpha}, X^{\alpha}) 
		\Big],
	$$
	and
	$$
		V^2_R(t,\xb) ~:= \sup_{\alpha \in \Ac^2_R(t,\xb) } \E^{\P^{\alpha}} \Big[
			\int_t^{\pi^{\alpha}} \!\!\! \int_U L(s, X^{\alpha}, u) M^{\alpha}_s(du) ds
			+
			\Phi( \pi^{\alpha}, X^{\alpha}) 
		\Big].
	$$
	Further, for the strong formulation, we recall that $\Uc_t$ denotes the collection of all $U$-value $\F^0$-predictable processes defined on $(\Om_0, \Fc^0, \P_0)$ which is independent of $\sigma(B_s ~: s \le t)$ under $\P_0$ (c.f. Section \ref{subsubsec:strong_form}).
 	Let us introduce
	$$
		\Uc^2_t ~:=~ \Big\{ \nu \in \Uc_t ~: \int_t^{\infty} | \nu_s - u_0 |^2 ds < \infty \Big\},
	$$
	and
	$$
		V_S^2(t, \xb) ~:=~ 
		\sup_{\nu \in \Uc^2_t} ~\sup_{\pi \in \Tc_t}~
		\E \Big[ \int_t^{\pi} L(s, X^{t, \xb, \nu}, \nu_s) ds + \Phi \big(\pi, X^{t,\xb,\nu}_{\cdot} \big) \Big].
	$$
	Notice that under the Lipschtiz condition in Assumption \ref{assu:Lip_coef}, together with the linear growth condition in \eqref{eq:linear_growth_cond},
	the controlled SDE \eqref{eq:controlled_SDE} has a unique solution for every $\nu \in \Uc^2_t$.
	
	\begin{Theorem}
		Assume that the coefficient functions $\mu$ and $\sigma$ are Borel measurable and satisfy \eqref{eq:linear_growth_cond},
		the reward functions $\Phi$ and $L$ are upper semi-analytic 
		and satisfy $\Phi(t, \xb) = \Phi(t, \xb_{t \wedge \cdot})$, $L(t, \xb, u) = L(t, \xb_{t \wedge \cdot}, u)$, for all $(t,\xb, u) \in \R_+ \x \Om \x U$.

		\noindent \rmi Then the both value functions $V^2_W: \Rb_+ \x \Om \longrightarrow \Rb$ and $V^2_R: \Rb_+ \x \Om \longrightarrow \Rb$ are upper semi-analytic.
		Moreover, for any $(t, \xb) \in \Rb_+ \x \Om$ and $\Fbb$-stopping time $\taub$ defined on $\Omb$ and taking value in $[t, \infty]$,
		by denoting $\tau^{\alpha} := \taub(\pi^{\alpha}, X^{\alpha}, M^{\alpha})$,
		one has the dynamic programming principle:
		\b*
			V^2_W(t,\xb) &\!\!=\!\!& \sup_{\alpha \in \Ac^2_W(t,\xb) } \E^{\P^{\alpha}} \Big[
			\Big(\int_t^{\pi^{\alpha}} \!\!\! L(s, X^{\alpha}, \nu^{\alpha}_s) ds + \Phi(\pi^{\alpha}, X^{\alpha})  \Big) \1_{\pi^{\alpha} \le \tau^{\alpha}} \\
			&&~~~~~~~~~~~~~~~~~~~~~+~
			\Big(\int_t^{\tau^{\alpha}} \!\!\! L(s, X^{\alpha}, \nu^{\alpha}_s) ds + V^2_W(\tau^{\alpha}, X^{\alpha}) \Big) \1_{\pi^{\alpha} > \tau^{\alpha}}
			\Big],
		\e*
		and
		\b*
			V^2_R(t,\xb) &\!\!=\!\!& \sup_{\alpha \in \Ac^2_R(t,\xb) } \E^{\P^{\alpha}} \Big[
			\Big(\int_t^{\pi^{\alpha}} \!\!\! \int_U L(s, X^{\alpha}, u) M^{\alpha}_s(du) ds + \Phi(\pi^{\alpha}, X^{\alpha})  \Big) \1_{\pi^{\alpha} \le \tau^{\alpha}} \\
			&&~~~~~~~~~~~~~~~~~+
			\Big(\int_t^{\tau^{\alpha}} \!\!\! \int_U L(s, X^{\alpha}, u) M^{\alpha}_s(du) ds + V^2_R(\tau^{\alpha}, X^{\alpha}) \Big) \1_{\pi^{\alpha} > \tau^{\alpha}}
			\Big].
		\e*

		\noindent \rmii Suppose in addition that Assumption \ref{assu:Lip_coef} holds true.
		Then $V^2_S : \Rb_+ \x \Om \longrightarrow \Rb$ is also upper semi-analytic.
		Moreover, for every $(t, \xb) \in \Rb_+ \x \Om$ and $\Ft$-stopping time $\tilde \tau$ defined on $\Omt$ and taking value in $[t, \infty]$, with $(\tau^{\nu, \pi})$ be defined in \eqref{eq:def_tau_tildetau},
		one has
		\b*
			V^2_S(t, \xb)
			&\!\!\!=\!\!\!&
			\sup_{\nu \in \Uc^2_t} ~\sup_{\pi \in \Tc_t}
			\E ~\Big[ \Big( \int_t^{\pi} L(s, X^{t,\xb,\nu}, \nu_s) ds + \Phi \big(\pi, X^{t,\xb,\nu} \big) \Big) \1_{\pi \le \tau^{\nu,\pi}} \\
			&&~~~~~~~~~~~~~~~~~~+
				\Big( \int_t^{\tau^{\nu,\pi}} \!\!\!\! L(s, X^{t,\xb,\nu}, \nu_s) ds + V^2_S \big(\tau^{\nu,\pi}, X^{t,\xb,\nu} \big) \Big) \1_{\pi > \tau^{\nu,\pi}} \Big].
		\e*
	
	\end{Theorem}
	\proof We will only provide the proof for the weak formulation to illustrate the additional technique needed in this new setting.
	The proofs for the relaxed formulation and strong formulation can be easily adapted with the techniques in Sections \ref{subsubsec:weak_relax_form} and \ref{subsubsec:strong_form}.
	
	\vspace{0.5em}
	
	First, we check easily that
	$$
		[[\Pcb^2_W]] 
		:=
		\Big\{
			(t, \xb, \Pb)  \in [[\Pcb_W]]
			~:
			\int_t^{\infty} \!\!\!\! \int_U \big| u - u_0 \big|^2 M_s(du) ds \le \infty
		\Big\}
	$$
	is still Borel, so that $V^2_W$ is also upper semi-analytic.
	
	\vspace{0.5em}
		
	Next, for every $\Pb \in \Pcb^2_W(t, \xb)$, and its r.c.p.d. $(\Pb_{\omb})_{\omb \in \Omb}$ knowing $\Fcb_{\tau}$,
	it is easy to check as in Theorems \ref{theo:ctrl_stopping_DPP} and \ref{theo:DPP_VWVR} that $\Pb_{\omb} \in \Pcb^2_W(\tau(\omb), \omb^X)$ for $\Pb$-a.e. $\omb \in \Omb$.
	This is enough to deduce that
	\b*
		V^2_W(t,\xb) &\!\!\le\!\!& \sup_{\alpha \in \Ac^2_W(t,\xb) } \E^{\P^{\alpha}} \Big[
		\Big(\int_t^{\pi^{\alpha}} \!\!\! L(s, X^{\alpha}, \nu^{\alpha}_s) ds + \Phi(\pi^{\alpha}, X^{\alpha})  \Big) \1_{\pi^{\alpha} \le \tau^{\alpha}} \\
		&&~~~~~~~~~~~~~~~~~~~~~+~
		\Big(\int_t^{\tau^{\alpha}} \!\!\! L(s, X^{\alpha}, \nu^{\alpha}_s) ds + V^2_W(\tau^{\alpha}, X^{\alpha}) \Big) \1_{\pi^{\alpha} > \tau^{\alpha}}
		\Big].
	\e*
	
	For the reverse inequality, we need to use the concatenation argument.
	To this end, let us introduce, for every $K > 0$ and $(t, \xb) \in \R_+ \x \Om$,
	$$
		\Ac^2_W(t, \xb, K) ~:=~ 
		\Big\{ \alpha \in \Ac^2_W(t, \xb) ~:
			\int_t^{\infty} \big| \nu^{\alpha}_s - u_0 \big|^2 ds \le K 
		\Big\},
	$$
	and
	$$
		V^2_W(t,\xb, K) ~:= \sup_{\alpha \in \Ac^2_W(t,\xb, K) } \E^{\P^{\alpha}} \Big[
			\int_t^{\pi^{\alpha}} L(s, X^{\alpha}, \nu^{\alpha}_s) ds
			+
			\Phi( \pi^{\alpha}, X^{\alpha}) 
		\Big].
	$$
	We notice that 
	\begin{equation} \label{eq:VWK2VW}
		V_W^2(t, \xb, K) \nearrow V^2_W(t, \xb),
		~~\mbox{as}~
		K \nearrow \infty,
	\end{equation}
	and with
	$$
		\Pcb^2_W(t, \xb, K) := \Big\{ 
			\Pb \in \Pcb_W(t, \xb) ~:
			\int_t^{\infty} \int_U \big| u - u_0 \big|^2 M_s(du) ds \le K
		\Big\},
	$$
	one has
	\begin{equation} \label{eq:V2WK_supP}
		V^2_W(t, \xb, K) = 
		\sup_{\Pb \in \Pcb^2_W(t, \xb, K)} \E^{\Pb} \Big[
			\int_t^{\Theta_{\infty}} \!\! \int_U L(s, X, u) M_s(du) ds +  \Phi(\Theta_{\infty}, X)
		\Big].
	\end{equation}
	Moreover, for every $K > 0$, the following graph set is still Borel measurable:
	\b*
		[[\Pcb^2_W(K)]] 
		&:=& 
		\big\{(t, \xb, \Pb) ~: \Pb \in \Pcb^2_W(t, \xb, K) \big\}\\
		&=&
		\Big\{
			(t, \xb, \Pb)  \in [[\Pcb_W]]
			~:
			\int_t^{\infty} \!\!\!\! \int_U \big| u - u_0 \big|^2 M_s(du) ds \le K
		\Big\}.
	\e*
	Now, for $K > 0$, by measurable selection theorem, let us choose a measurable family $(\Pb^{\eps, K}_{t, \xb})_{(t, \xb) \in \R_+ \x \Om}$,
	where for each $(t, \xb)$, $\Pb^{\eps, K}_{t, \xb}$ is an $\eps$-optimal weak control rule for the problem at the r.h.s. of \eqref{eq:V2WK_supP}.
	Then, for every $\Pb \in \Pcb^2_W(t, \xb)$, we let $\Qb^{\eps, K}_{\omb} := \Pb^{\eps, K}_{\tau(\omb), \omb^X}$,
	and then consider the concatenated probability measure $\Pb \otimes_{\tau} \Qb^{\eps, K}_{\cdot}$.
	Following the arguments in Theorems \ref{theo:ctrl_stopping_DPP} and \ref{theo:DPP_VWVR}, 
	one can check directly that $\Pb \otimes_{\tau} \Qb^{\eps, K}_{\cdot} \in \Pcb^2_W(t, \xb)$, which implies that
	\b*
		V^2_W(t,\xb) &\!\!\ge\!\!& \sup_{\alpha \in \Ac^2_W(t,\xb) } \E^{\P^{\alpha}} \Big[
		\Big(\int_t^{\pi^{\alpha}} \!\!\! L(s, X^{\alpha}, \nu^{\alpha}_s) ds + \Phi(\pi^{\alpha}, X^{\alpha})  \Big) \1_{\pi^{\alpha} \le \tau^{\alpha}} \\
		&&~~~~~~~~~~~~~~~~~~~~~+~
		\Big(\int_t^{\tau^{\alpha}} \!\!\! L(s, X^{\alpha}, \nu^{\alpha}_s) ds + V^2_W(\tau^{\alpha}, X^{\alpha}, K) \Big) \1_{\pi^{\alpha} > \tau^{\alpha}}
		\Big].
	\e*
	Now, let $K \nearrow \infty$, and by \eqref{eq:VWK2VW} together with the monotone convergence theorem, 
	one can conclude the proof of the dynamic programming principle.
	\qed

	\begin{Remark}
		One can of course consider other growth conditions on $\mu$ and $\sigma$ than \eqref{eq:linear_growth_cond},
		and add other adapted integrability conditions on the admissible control process in the definition of $\Ac_W$, $\Ac_R$ and $\Uc$
		to formulate the problem,
		and then adapt the above techniques to prove the DPP.
	\end{Remark}

\section{Approximation and equivalence of different formulations of the optimal control/stopping problems}
\label{sec:stability}

	We will study an approximation problem of the relaxed control/stopping rule by weak control/stopping rules,
	which can be considered as a stability property.
	In particular, this consists of an important technical step to prove the equivalence between different formulations (strong, weak and relaxed formulations) of the optimal controlled/stopped diffusion processes problem.

\subsection{Approximation of relaxed control by weak control rules}

\subsubsection{Relaxed control rule in an abstract probability space}

	The martingale problem in Section \ref{subsec:generator} is formulated on the canonical space without fixing the equipped probability measures.
	To obtain a similar formulation of relaxed control in a fixed and abstract filtered probability space, 
	one can make use of a product space together with the notion of stable convergence topology of Jacod and M\'emin \cite{JacodMemin}.	
	
	\vspace{0.5em}

	Let  $(\Om^*, \Fc^*)$ be a fixed measurable space equipped with the filtration $\F^* = (\Fc^*_t)_{t \ge 0}$,
	we denote by $\Tc^*$ the collection of all $\F^*$-stopping times.
	Recall that $\Om := \D(\R_+, E)$ denotes the canonical space of all c\`adl\`ag $E$-valued paths on $\R_+$, equipped with the Skorokhod topology, and the canonical filtration $\F = (\Fc_t)_{t \ge 0}$. 
	Let us introduce an enlarged space  $\Omb^* := \Om^* \x \Om$,
	equipped with  the $\sigma$-field $\Fcb^*_{\infty} := \Fc^* \otimes \Fc_{\infty}$, 
	and the enlarged filtration $\Fbb^* = (\Fcb^*_t)_{t \ge 0}$ defined by $\Fcb^*_t := \Fc^*_t \otimes \Fc_t$.
	On $\Omb^*$, let $X$ be the canonical process defined by $X_t (\omb^*) := \om_t$ for all $\omb^* = (\om^*, \om) \in \Omb^*$.
	Let $B_{mc}(\Omb^*)$ denote the collection of all bounded $\Fcb^*_{\infty}$-measurable functions 
	$\xi: \Omb^* \longmapsto \R$ such that for every $\om^* \in \Om^*$, the mapping $\xb \longmapsto \xi(\om^*, \xb)$ is continuous.
	Denote also by $\Pcb(\Omb^*)$ (resp. $\Pc(\Om)$, $\Pc(\Om^*)$) the collection of all probability measures on $\big( \Omb^*, \Fcb^*_{\infty} \big)$ (resp. $(\Om, \Fc_{\infty})$, $(\Om^*, \Fc^*)$).
	Let $\P^* \in \Pc(\Om^*)$ be a fixed probability measures, 
	we define 
	$$
		\Pcb(\P^*)
		~:=~
		\big\{
			\Pb \in \Pcb(\Omb^*) ~:
			\Pb|_{\Om^*} = \P^*
		\big\}.
	$$
	
	\begin{Definition}
		The stable convergence topology on $\Pcb(\Omb^*)$ is defined as the coarsest topology for which the mapping $\Pb \longmapsto \E^{\Pb}[\xi]$ is continuous for all $\xi \in B_{mc}(\Omb^*)$.
	\end{Definition}
	
	In the following, we equip $\Pcb(\Omb^*)$ with the table convergence topology,
	and $\Pc(\Om)$ with the weak convergence topology (i.e. the coarsest topology such that $\P \longmapsto \E^{\P}[\xi]$ is continuous for all bounded continuous functions $\xi$ on $\Om$),
	and $\Pc(\Om^*)$ with the coarsest topology such that $\P \longmapsto \E^{\P}[\xi]$ is continuous for all bounded measurable functions $\xi$ on $(\Om^*, \Fc^*)$.
	One has the following results on stable convergence topology from \cite{JacodMemin}.

	\begin{Proposition} \label{prop:StableCvg}
	{\rm	\rmi A subset $\Pcb$ of $\Pcb(\Omb^*)$ is relatively compact w.r.t. the stable topology if and only if 
		$\Pcb_{\Om^*} := \{ \Pb|_{\Om^*}: \Pb \in \Pcb \}$ 
		and $\Pcb_{\Om} := \{ \Pb |_{\Om}: \Pb \in \Pcb \}$ 
		are both relatively compact in $\Pc(\Om^*)$ and $\Pc(\Om)$, respectively.
		
		\vspace{0.5em}
		
		\noindent \rmii Let $(\Pb_n)_{n \ge 1} \subset \Pcb(\Omb^*)$ be a sequence such that $\Pb_n \longrightarrow \Pb_{\infty}$ under the stable convergence topology,
		and $\xi: \Omb^* \longrightarrow \R$ be a bounded and $\Fcb^*_{\infty}$-measurable function,
		such that for every $\om^* \in \Om^*$, the mapping $\xb \longmapsto \xi(\om^*, \xb)$ is continuous.
		Then one has $\lim_{n \longrightarrow \infty} \E^{\Pb_n}[\xi] = \E^{\Pb_{\infty}}[\xi]$.

		\vspace{0.5em}
		
		\noindent \rmiii Let $(\Pb_n)_{n \ge 1} \subset \Pcb(\Omb^*)$ be a sequence such that $\Pb_n \longrightarrow \Pb_{\infty}$ under the stable convergence topology,
		and $\xi: \Omb^* \longrightarrow \R$ be a bounded $\Fcb^*_{\infty}$-measurable function,
		such that the set $\{(\om, \xb) \in \Omb^* ~: \xb' \mapsto \xi(\om, \xb') ~\mbox{is discontinuous at}~\xb \}$ is $\Pb_{\infty}$-negligible.
		Then one has $\lim_{n \longrightarrow \infty} \E^{\Pb_n}[\xi] = \E^{\Pb_{\infty}}[\xi]$.
		
		\vspace{0.5em}
		
		\noindent \rmiv  Let $\P^* \in \Pc(\Om^*)$ be a fixed probability measure, and  $(\Pb_n)_{n \ge 1} \subset \Pcb(\P^*)$ be a relatively compact sequence (under the stable convergence topology).
		Then there exists a subsequence $(\Pb_{n_k})_{k \ge 1}$ and  $\Pb_{\infty} \in \Pcb(\P^*)$ such that $\Pb_{n_k} \longrightarrow \Pb_{\infty}$.
		
		\vspace{0.5em}
		
		\noindent \rmv Assume that $\Om^*$ is a Polish space, $\Fc^*$ is its Borel $\sigma$-field and $\P^* \in \Pc(\Om^*)$.
		Then restricted on $\Pcb(\P^*)$, the stable convergence topology coincides with the weak convergence topology.
	}
	\end{Proposition}

	Now we are ready to introduce a notion of relaxed control rule by using a martingale problem on $\Omb^*$, 
	which is also in the same spirit of Jacod and M\'emin \cite{JacodMemin_SDE}.
	On the filtered probability space $(\Om^*, \Fc^*, \F^*, \P^*)$, we denote by $\Uc^*$ the set of all $\Pc(U)$-valued $\F^*$-predictable processes $m^* = (m^*_t)_{t \ge 0}$, where $\Pc(U)$ is the set of all Borel probability measures on $U$.
	By naturally extension, one can also consider $m^*$ as a $\Fbb^*$-predictable process defined on $\Omb^*$.

	\vspace{0.5em}
	
	As in Section  \ref{subsec:generator}, we consider a generator $\G$ of a control problem, which is a subset of $C_b(E) \x B(\R_+ \x \Om \x E \x U)$.
	Let $x_0 \in E$ be fixed, and $m^* \in \Uc^*$,
	a relaxed control rule with initial condition $x_0$ and control process $m^*$ is a probability measure $\Pb^* \in \Pcb(\Omb^*)$ 
	such that $\Pb^*|_{\Om^*} = \P^*$,
	$\Pb^*[X_0 = x_0] = 1$,
	and the process $\big(C^{m^*}_t(f,g) \big)_{t \ge 0}$ is a $(\Pb^*, \Fbb^*)$-martingale for every $(f,g) \in \G$, with
	\be \label{eq:Cfg_star}
		C^{m^*}_t(f,g)
		~:=~
		f(X_t) - \int_0^t \int_U g(s, X_{s \wedge \cdot}, u, X_s) m^*_s(du) ds.
	\ee
	When $m^*$ is induced by a $U$-valued $\F^*$-predictable process $\nu^*$ in the sense that $m^*(du,ds) = \delta_{\nu^*_s}(du)ds$, $\P^*$-a.s.,
	we also call $\Pb^*$ a weak control rule.
	Let us denote by $\Pcb(m^*)$ the set of all relaxed control rules with control process  $m^* \in \Uc^*$ (the initial condition $x_0$ is fixed).

	\begin{Theorem} \label{theo:equiv_0}
		Assume that, for all functions $(f,g) \in \G$ in the generator of the control problem, 
		the function $g$ is uniformly bounded and the map $(\xb, x) \longmapsto g(t, \xb,u,  x)$ is continuous for each $t \in \R_+$ and $u \in U$.
		Let $(m^n)_{n \ge 1} \subset\Uc^*$ be a sequence such that $m^n \longrightarrow m^{\infty} \in \Uc^*$, $\P^*$-a.s.,
		and $(\Pb^*_n)_{n \ge 1}$ be a sequence such that $\Pb^*_n \in \Pcb(m^n)$, for all $n \ge 1$.
		Assume in addition that
		$\Pb^*_n \longrightarrow \Pb^*_{\infty}$ (under the stable convergence topology),
		and that, for all $s \le t$ and $\Fcb_s$-measurable bounded r.v. $\xi$,
		\begin{equation} \label{eq:gmn_cvg}
			\E^{\Pb^*_n} \Big[ \Big(\! \int_s^t\!\! \int_U g(r, X, u, X_r)  \big( m^n_r(du) \!-\! m^{\infty}_r(du) \big) dr \Big) \xi \Big] 
			\longrightarrow 
			0,
			~\mbox{as}~n \longrightarrow \infty.
		\end{equation}
		Then $\Pb^*_{\infty} \in \Pcb(m^{\infty})$.
	\end{Theorem}
	\proof Notice that $\Pb^*_n|_{\Om^*} = \P^*$ and $\Pb^*_n[X_0 = x_0] = 1$ for all $n \ge 1$.
	Moreover, the map $\xb \longmapsto \xb_0$ from $\Om$ to $\R$ is continuous under the Skorokhod topology.
	Then it is clear that $\Pb^*_{\infty}|_{\Om^*} = \P^*$ and $\Pb^*_{\infty}[X_0 = x_0] = 1$.
	
	\vspace{0.5em}
	
	Further, for all $s \le t$, and any bounded $\Fcb^*_s$-measurable random variable $\xi$ such that $\xb \mapsto \xi(\om, \xb)$ is continuous,
	it follows by the martingale property that
	$$
		\E^{\Pb^*_n} \Big[
			\Big( f(X_t) - f(X_s) - \int_s^t \!\! \int_U g(r, X, u, X_r) m^n_r(du) dr \Big)
			\xi
		\Big]
		=0,
		~\mbox{for all}~ n\ge 1.
	$$
	Next, by \eqref{eq:gmn_cvg}, one obtains that
	$$
		\lim_{n \longrightarrow \infty} \E^{\Pb^*_n} \big[\zeta_{s,t} \xi \big] = 0,
		~\mbox{with}~
		\zeta_{s,t} (\om, \xb) 
		:= 
		 f(\xb_t) - f(\xb_s) - \int_s^t \!\!\! \int_U g(r, \xb, u, \xb_r) m^{\infty}_r \!(\om, du) dr.
	$$
	Further, with the given probability $\Pb^*_{\infty}$, there exists a countable set $\T \subset (0, \infty)$ such that
	$\xb \longmapsto \xb_t$ is continuous (under the Skorokhod topology) for $\Pb^*_{\infty}$-a.e. $\xb \in \Om$ (see e.g. Jacod and Shiryaev \cite[Lemma IV.3.12]{JacodShiryaev}).
	Together with the continuity of $(\xb, x) \longmapsto g(t, \xb, u ,x)$, 
	this is enough to deduce that 
	$\{ (\om, \xb) \in \Omb^* ~: \xb \longmapsto \zeta_{s,t}(\om, \xb) ~\mbox{is discontinuous} \}$ is $\Pb^*_{\infty}$-negligible whenever $s, t \in \R_+ \setminus \T$.
	Therefore, one has $\E^{\Pb^*_{\infty}} \big[ \zeta_{s,t}\xi \big] = 0$ for all $s\le t$ such that $s, t \in \R_+ \setminus \T$.
	This is enough to conclude that $(C^{m^*_{\infty}}_t(f,g))_{t \ge 0}$ is a $(\Pb^*_{\infty}, \Fbb^*)$-martingale,
	and hence $\Pb^*_{\infty} \in \Pcb(m^{\infty})$.
	\qed

	\vspace{0.5em}

	In practice, one usually fix a given $\Pc(U)$-valued (relaxed) control process $m^{\infty}$ and construct a sequence $m^n$ to approximate $m^{\infty}$.
	In particular, $m^n$ can be chosen be to a relaxed control induced by a $U$-valued (weak) control process.
	This is the so called Fleming's chattering lemma, which is recalled below.

	\begin{Lemma}[Flemming's chattering lemma] \label{lemm:chattering}
		For every relaxed control $m^{\infty} \in \Uc^*$, there is a sequence of $U$-valued control processes $(\nu^n)_{n \ge 1}$ 
		such that each $\nu^n$ is $\F^*$-adapted and piecewise constant in the sense that $\nu^n_t = \nu^n_{t_k}$ for all $t \in [t_k, t_{k+1})$ with a some discrete time grid $0 = t_0 < t_1 < \cdots$.
		Moreover, the induced measure valued process 
		$m^{\nu^n}(dt, du) := \delta_{\nu_t^n}(du) dt$ converges in $\M$ to $m^{\infty}(du, dt)$, $\P^*$-a.s.
	\end{Lemma}

	\begin{Remark} \label{rem:piece_cst_ctrl}
	{\rm
		\rmi The proof of Theorem \ref{theo:equiv_0} is in the same spirit of the classical limit arguments in the proof of existence of solutions to the (uncontrolled) martingale problem (see e.g. Stroock and Varadhan \cite{StroockVaradhan}, or Protter \cite{Protter}).
		In that setting, one can approximate functions $(f,g) \in \G$ by more regular functions $(f_n, g_n)$ (or even piecewise constant functions),
		whose martingale problems have easily solutions and the limit provides a solution to the original martingale problem.

		\vspace{0.5em}
		
		\noindent \rmii Together with Lemma \ref{lemm:chattering}, one can use Theorem \ref{theo:equiv_0} to approximate a relaxed control rule by weak control rules.
		Indeed, given a relaxed control process $m^{\infty}$, one can first approximate it by a sequence of weak control processes $(\nu^n)_{n \ge 1}$ in the sense of Lemma \ref{lemm:chattering}.
		Next, under standard conditions, it is easy to check that the sequence $(\Pb^*_n)_{n \ge 1}$ of the associated weak control rules is relatively compact, so that one can take a subsequence of weak control rules converges to some probability measure $\Pb^*_{\infty}$.
		The rest is to check Condition \eqref{eq:gmn_cvg} and that the set $\Pcb(m^{\infty})$ of relaxed control rules is unique so that the weak control rules converges to the given relaxed control rule.
		In Section \ref{subsubsec:approximation}, we will show how to check \eqref{eq:gmn_cvg} and how to obtain uniqueness of $\Pcb(m^{\infty})$ in the context of the controlled diffusion processes problem.

		\vspace{0.5em}
		
		\noindent \rmiii In Theorem \ref{theo:equiv_0}, the boundedness condition on $g$ for any $(f,g) \in \G$ is a technical condition for simplicity, which can be relaxed in concrete examples. See also discussions in Remark \ref{eq:cond_bound_compact}.
	}
	\end{Remark}

\subsubsection{Approximation of relaxed control/stopping rules in the diffusion processes setting}
\label{subsubsec:approximation}

	In this section, we stay in the controlled diffusion process setting, and provide an approximation result of relaxed control rule by weak control rules.
	More precisely, let $(\mu, \sigma): \R_+ \x \Om \x U \longrightarrow \R^d \x \S^d$ be the coefficient functions of the controlled diffusion process,
	denote
	$$
		\Lc \varphi(s,\xb, u, x) 
		:=
		\mu(s, \xb, u) \cdot D \varphi(x) + \frac12 \mbox{Tr} \big( \sigma \sigma^{\top}(s, \xb, u) D^2 \varphi(x) \big).
	$$
	Then the generator $\G$ of the controlled diffusion process problem is given by
	\begin{equation} \label{eq:diffusion_generator}
		\G ~:=~ \big\{
			(\varphi, \Lc \varphi)~: \varphi \in C^{\infty}_c(\R^d)
		\big\}.
	\end{equation}
	We make the following conditions throughout this subsection.
	
	\begin{Assumption} \label{assu:bbd_Lip_coef}
	\rmi The coefficient functions $\mu$ and $\sigma$ are uniformly bounded and $\F$-progressive in the sense that $(\mu, \sigma)(t, \xb, u) = (\mu, \sigma)(t, \xb_{t \wedge \cdot}, u)$ for all $(t, \xb, u) \in \R_+ \x \Om \x U$.
	Further, for all $T > 0$, there is some constant $L_0 > 0$ such that,
		for all $(t, \xb, \xb', u) \in [0,T] \x \Om \x \Om \x U$, with $\| \xb\|_T := \sup_{0 \le t\le T} |\xb_t|$, one has
		$$
			|\mu(t, \xb, u) - \mu(t, \xb', u)| ~+~ \|\sigma(t, \xb, u) - \sigma(t, \xb', u) \| 
			~\le~
			 L_0 \|\xb - \xb'\|_T.
		$$
		Assume in addition that $(\mu, \sigma)(t, \xb, u)$ are uniformly continuous in $t$ in the sense that, for all $\eps > 0$, there exists $\delta > 0$, such that
		for all $\xb \in \Om$, $u \in U$ and $s \le t$ satisfying $t-s \le \delta$, one has
		$$
			\big| \mu(t, \xb_{s \wedge \cdot}, u) - \mu(s, \xb_{s \wedge \cdot}, u) \big|
			+
			\big \| \sigma(t, \xb_{s \wedge \cdot}, u) - \sigma(s, \xb_{s \wedge \cdot}, u) \big \|
			~\le~
			\eps.
		$$
	\noindent \rmii The set $U$ is compact, and the map $u \mapsto (\mu, \sigma)(t, \xb, u)$ is uniformly continuous, uniformly in $(t, \xb) \in \R_+ \x \Om$.

	\end{Assumption}

	\begin{Remark} \label{eq:cond_bound_compact}
		\rmi The coefficient functions $\mu$ and $\sigma$ are assumed to be bounded for simplicity.
		One can easily consider the setting with Condition \eqref{eq:mu_sigma_integ_p} (and Assumption \ref{assu:Lip_coef}), 
		or the linear growth setting in Section \ref{subsubsec:unbound_mu_sigma}.
		In fact, by using a simple truncation technique, one can easily approximate a diffusion process by those with bounded drift and diffusion coefficient functions.
		
		\vspace{0.5em}
		
		\noindent \rmii Similarly, when $U$ is not compact, 
		one can also use truncation technique to reduce the approximation problem to the setting with compact set $U$.
		This would be quite standard if $U$ is a non-compact subset of $\R^n$, and the coefficient functions $\mu$ and $\sigma$ satisfy some growth condition in $u$.
	\end{Remark}

	In the following, let us fix a relaxed control process $m^* \in U^*$, 
	and let $\Pb^*$ be a fixed relaxed control rule associated with the generator $\G$ given in \eqref{eq:diffusion_generator}.
	By  \cite{ElKaroui_Meleard} (see also Proposition \ref{prop:weak_reformulation}), 
	there exists (in a possibly enlarged space) a continuous martingale measure $\widehat M^*$ with quadratic variation $m^*(du,dt)$ such that
	$$
		X_t = x_0 + \int_0^t \!\! \int_U \mu(s, X_{s \wedge \cdot}, u) m^*_s(du) ds 
		+
		\int_0^t  \!\! \int_U \sigma(s, X_{s \wedge \cdot}, u) \widehat M^* (du, ds),
		~t \ge 0,
		~\Pb^* \mbox{-a.s.}
	$$

\paragraph{Step 1: approximation by relaxed control rules supporting in a finite control space}

	In a first step, we will approximate the relaxed control $m^*$ by a specially constructed approximating sequence.
	Since $U$ is a compact metric space, for all $\eps > 0$, there exists a partition $(U^{\eps}_i)_{i = 1, \cdots, N_{\eps}}$ of $U$
	(i.e.  $\cup_{i=1}^{N_{\eps}} U^{\eps}_i = U$, and $U^{\eps}_i \cap U^{\eps}_j = \emptyset$ whenever $i \neq j$) 
	together with a set $(u^{\eps}_i)_{i=1, \cdots, N_{\eps}}$
	such that $u^{\eps}_i \in U^{\eps}_i$, and $d(u, u^{\eps}_i) \le \eps$ for all $u \in U^{\eps}_i$, $i=1, \cdots, N_{\eps}$.
	
	\vspace{0.5em}
	
	For every $\eps > 0$, let us define
	$$
		m^{*,\eps}_s(du) ~:=~ \sum_{i=1}^{N_{\eps}} q^{\eps,i}_s \delta_{u^{\eps}_i}(du),
		~\mbox{with}~
		q^{\eps,i}_s := m^*_s(U^{\eps}_i),
	$$
	and define $\widehat M^{*,\eps}$ by, for all compactly supported measurable functions $\phi: \R_+ \x U \longrightarrow \R$,
	$$
		\int_0^{\infty} \!\!\! \int_U \phi(s, u) \widehat M^{*,\eps}(du, ds) 
		~:=~
		\sum_{i=1}^{N_{\eps}} \int_0^{\infty} \!\!\! \int_{U^{\eps}_i} \phi(s, u^{\eps}_i) \widehat M^*(du, ds).
	$$
	We notice that $\widehat M^{*,\eps}$ is a continuous martingale measure with quadratic variation $m^{*,\eps}$, w.r.t. the same filtration that generated by $(\widehat M^*, m^*)$.
	Let us define $X^{\eps}$ by SDE
	\begin{equation} \label{eq:SDE_Neps}
		X^{\eps}_t = x_0 + \int_0^t \!\!\! \int_U \mu(s, X^{\eps}_{s \wedge \cdot}, u) m^{*,\eps}_s(du) ds 
		+
		\int_0^t  \!\!\! \int_U \sigma(s, X^{\eps}_{s \wedge \cdot}, u) \widehat M^{*,\eps} (du, ds),
		~t \ge 0,
		~\Pb^* \mbox{-a.s.}
	\end{equation}

	\begin{Remark} \label{rem:relaxed_rule_finite}
		\rmi $X^{\eps}$ can be considered as a controlled process with relaxed control process $m^{*,\eps}$, 
		which supports in a finite control space $\{ u^{\eps}_i,~i=1, \cdots, N_{\eps} \} \subset U$.
		More precisely, the probability $\Pb^{*,\eps}$ defined below can be considered as a relaxed control rule with control process $m^{*,\eps}$: for all bounded measurable $\phi: \Omb^* \longrightarrow \R$,
		$$
			\int_{\Omb^*} \phi(\om, \xb) \Pb^{*,\eps}(d\om, d\xb) 
			~:=~
			\int_{\Omb^*} \phi \big(\om, X^{\eps}(\om,\xb) \big) \Pb^*(d\om, d\xb).
		$$

		\noindent \rmii Since $m^{*,\eps}$ supported in a finite space, there exists (in a possibly enlarged space) $N_{\eps}$ independent Brownian motion $(B^{\eps,i})_{i=1, \cdots, N_{\eps}}$, and rewrite the SDE \eqref{eq:SDE_Neps} equivalently
		\begin{equation} \label{eq:SDE_NepsBM}
			X^{\eps}_t = x_0 + \sum_{i=1}^{N_{\eps}} 
			\Big(\!
			\int_0^t\!\! \mu(s, X^{\eps}_{s \wedge \cdot}, u^{\eps}_i) q^{\eps,i}_s ds 
			+
			\int_0^t \!\! \sigma(s, X^{\eps}_{s \wedge \cdot}, u^{\eps}_i) \sqrt{q^{\eps,i}_s} dB^{\eps,i}_s
			\Big),
			~t \ge 0,
			~\Pb^* \mbox{-a.s.}
		\end{equation}
	\end{Remark}

	\begin{Proposition} \label{prop:approximate_q}
		Let Assumption \ref{assu:bbd_Lip_coef} hold true.
		Let $\Pb^{*,\eps}$ be given as in Remark \ref{rem:relaxed_rule_finite} with the above construction.
		Then, as $\eps \longrightarrow 0$, one has $m^{*,\eps} \longrightarrow m^*$, $\Pb^*$-a.s. and
		$$
			\Pb^{*,\eps} \longrightarrow \Pb^*
			~\mbox{under the stable convergence topology on}~
			\Pcb(\Omb^*).
		$$
		Moreover, in the approximation sequence, one can choose $\Pb^{*,\eps}$ be such that the processes $(q^{\eps,i}))_{i=1, \cdots, N_{\eps}}$ are piecewise constant in the sense that, for a discrete time grid $0 = t_0 < t_1 < \cdots$, one has $q^{\eps,i}_t = q^{\eps,i}_{t_k}$ for $t \in [t_k, t_{k+1})$ and $k \ge 1$.
	\end{Proposition}
	\proof
	\rmi First, by its construction, one has $m^{*,\eps} \longrightarrow m^*$ as $\eps \longrightarrow 0$.
	Next, we will prove that, for all $T > 0$,
	\begin{equation} \label{eq:cvg_L2_eps}
		\E^{\Pb^*} \Big[ \sup_{0 \le t \le T} \big| X^{\eps}_t - X_t \big|^2 \Big] \longrightarrow 0, ~\mbox{as}~\eps \longrightarrow 0,
	\end{equation}
	which is enough to conclude that 
	$\Pb^{*,\eps} \longrightarrow \Pb^*$.
	To prove \eqref{eq:cvg_L2_eps}, we notice that $(\mu, \sigma)$ are uniformly continuous in $u$.
	Then for all $\delta > 0$, there exists $\eps > 0$ such that
	$$
		\big| \mu(t, \xb, u) - \mu(t, \xb, u^{\eps}_i) \big| + \big\| \sigma(t, \xb, u) - \sigma(t, \xb, u^{\eps}_i) \big\| 
		\le \delta,
		~\mbox{for all}~ u \in U^{\eps}_i, ~i=1, \cdots, N_{\eps}.
	$$
	One then obtains that
	\b*
		X_t - X^{\eps}_t &=& \int_0^t \int_U \big(\mu(s, X_s, u) - \mu(s, X^{\eps}_s, u) \big) m^*_s(du) ds \\
		&&+ \int_0^t \int_U  \big( \sigma(s, X_s, u) - \sigma(s, X^{\eps}_s, u)  \big) \widehat M^*(du, s) 
		+ R^{\eps}_t,~\Pb^*\mbox{-a.s.},
	\e*
	where $R^{\eps}$ satisfies $\E^{\Pb^*}\big[ \sup_{0 \le t \le T} |R^{\eps}_t|^2 \big] \le C \delta ^2 $ for some constant $C$ (depending on $T$).
	Using the Lipschitz property of $(\mu, \sigma)$ in $\xb$, by standard arguments in the SDE theory (with It\^o's isometry, Doob's martingale inequality, and Gromwall lemma),
	one can easily prove \eqref{eq:cvg_L2_eps}.

	\vspace{0.5em}

	\noindent \rmii
	To prove that the processes $(q^{\eps,i}))_{i=1, \cdots, N_{\eps}}$ can be chosen to be piecewise constant,
	we fix the process $X^{\eps}$ in \eqref{eq:SDE_NepsBM} and approximate it  by controlled processes with piecewise constant (relaxed) control.
	Indeed, for the given (progressively measurable) processes $(q^{\eps,i}))_{i=1, \cdots, N_{\eps}}$, 
	one can approximate it by a sequence $(q^{\eps,i, n})_{n \ge 1}$ of (adapted) piecewise constant processes 
	(see e.g. Karatzas and Shreve \cite[Lemma 3.2.4]{KaratzasShreve}) in the sense that
	\begin{equation} \label{eq:approximate_q}
		\lim_{n \longrightarrow \infty} 
		\int_0^T \Big( 
			\Big| \sqrt{q^{\eps, i, n}_t} - \sqrt{q^{\eps, i}_t} \Big |^2
			+ 
			\Big| {q^{\eps, i, n}_t} - {q^{\eps, i}_t} \Big |^2
		\Big) dt = 0.
	\end{equation}
	Moreover, by adding a renormalization step in the proof of \cite[Lemma 3.2.4]{KaratzasShreve}, one can ensure that $\sum_{i=1}^{N_{\eps}} q^{\eps,i,n}_t = 1$ for all $t \ge 0$.
	Let us now define $X^{\eps,n}$ by the SDE 
	$$
		X^{\eps,n}_t = x_0 + \sum_{i=1}^{N_{\eps}} 
		\Big(\!
		\int_0^t\!\! \mu(s, X^{\eps,n}_{s \wedge \cdot}, u^{\eps}_i) q^{\eps,i,n}_s ds 
		+
		\int_0^t \!\! \sigma(s, X^{\eps,n}_{s \wedge \cdot}, u^{\eps}_i) \sqrt{q^{\eps,i,n}_s} dB^{\eps,i}_s
		\Big),
		~t \ge 0,
		~\Pb^* \mbox{-a.s.}
	$$
	Then $X^{\eps,n}$ can be considered as a controlled diffusion process with (relaxed) control process $m^{*, \eps, n}(du,ds) := \sum_{i=1}^{N_{\eps}} q^{\eps, i,n}_s \delta_{u^{\eps}_i}(du) ds$.
	In particular, one has $m^{*,\eps, n} \longrightarrow m^{*, \eps}$, a.s. as $n \longrightarrow \infty$.
	Further, notice that $(\mu, \sigma)$ is uniformly bounded,
	by using \eqref{eq:approximate_q} together with standard arguments in the SDE theory, one can prove that
	$$
		\E^{\Pb^*} \Big[ \sup_{0 \le t \le T} \big| X^{\eps, n} - X^{\eps} \big|^2 \Big]
		~\longrightarrow~
		0.
	$$
	This is enough to conclude that the relaxed control rule induced by the piecewise constant (relaxed) control process $m^{*,\eps, n}$ (together with the associated  controlled process $X^{\eps, n}$)
	converges to $\Pb^{*,\eps}$ under the stable convergence topology.
	\qed

\paragraph{Step 2: approximation by weak control rules}

	We now approximate the relaxed control rules by whose with control processes taking value in a finite space,
	and the controlled process are given in the form of \eqref{eq:SDE_NepsBM} with piecewise constant processes $(q^{\eps, i})_{i =1 ,\cdots, N_{\eps}}$.
	Moreover, for ease of presentation, we assume $N_{\eps} = 2$ and then omit $\eps$ in the notation.
	Namely, we fix a relaxed control rule $\Pb^*$ with control process $(m^*_s(du))_{s \ge 0}$ satisfying
	$$
		m^*_s(du) = q^1_s \delta_{u_1}(du) + q^2_s \delta_{u_2} (du),
		~~
		(q^1_s, q^2_s) = (q^1_{t_k}, q^2_{t_k}),~\mbox{for}~s \in [t_k, t_{k+1}),~k\ge 0,
	$$
	where $0 = t_0 < t_1 < \cdots$ is a discrete time grid  on $[0, \infty)$.
	In particular, one has $q^1_s + q^2_s = 1$, and there exists 2 independent Brownian motion $B^1$ and $B^2$ such that
	\begin{equation} \label{eq:SDE_2BM_relax}
		X_t = x_0 + \sum_{i=1}^2 
		\Big(\!
		\int_0^t\!\! \mu(s, X_{s \wedge \cdot}, u_i) q^{i}_s ds 
		+
		\int_0^t \!\! \sigma(s, X_{s \wedge \cdot}, u_i) \sqrt{q^{i}_s} dB^{i}_s
		\Big),
		~t \ge 0,
		~\Pb^* \mbox{-a.s.}
	\end{equation}
	
	We next construct a sequence $(\nu^{*,n})_{n \ge 1}$ of $U$-valued control processes in $(\Om^*, \Fc^*, \P^*)$ to approximate $m^*$.
	For each $k \ge 0$, let us construct $\nu^{*,n}$ on time interval $[t_k, t_{k+1})$.
	First, let us consider a subdivision $t_k = t_{k,0} < t_{k,1}< \cdots < t_{k,n} = t_{k+1}$, where $t_{k,i+1} - t_{k,i} = (t_{k+1} - t_k)/n$ for each $i = 0, \cdots, n-1$.
	Next, let $\theta_{k,i} \in [t_{k,i}, t_{k,i+1}]$ be such that $(\theta_{k,i} - t_{k,i})/(t_{k,i+1} - t_{k,i}) = q^1_{t_k}$.
	Finally, let
	\begin{equation} \label{eq:def_mn_star}
		\nu^{*,n}_s
		:=
		\sum_{k =0}^{\infty} \sum_{i=0}^{n-1}
		\big( 
			u_1 \1_{\{s \in [t_{k,i}, \theta_{k,i})\}} + u_2 \1_{\{s \in [\theta_{k,i}, t_{k,i+1})\}}
		\big)
		~\mbox{and}~
		m^{*,n} (du, ds) := \delta_{\nu^{*,n}}(du) ds.
	\end{equation}
	Then $\nu^{*,n}$ is a $U$-valued $\F^*$-adapted piecewise constant control process.
	Moreover, one can check that $m^{*,n} \longrightarrow m^*$, $\P^*$-a.s.
	(see also \cite[Section 4]{ElKaroui_1987} for a detailed proof).
	We also notice that, for all measurable function $\phi: U \longrightarrow \R$, one has
	\begin{equation} \label{eq:phi_mn_m}
		\int_{t_{k,i}}^{t_{k,i+1}} \!\!\! \phi( \nu^{*,n}_s) ds 
		=
		\int_{t_{k,i}}^{t_{k,i+1}} \!\!\!\int_U \phi( u) m^{*,n}(du,ds)
		=
		\int_{t_{k,i}}^{t_{k,i+1}} \!\!\!\int_U \phi( u) m^{*}(du,ds),
		~\P^*\mbox{-a.s.}
	\end{equation}

	\begin{Proposition} \label{prop:qmn_cvg}
		Let Assumption \ref{assu:bbd_Lip_coef} hold true, and $(\nu^{*,n})_{n \ge 1}$ be given as above.
		For each $n \ge 1$,
		let $X^n$ be a controlled process corresponding to the control process $\nu^{*,n}$, so that there exists a Brownian motion $B^{*,n}$ such that
		\begin{equation} \label{eq:def_Xn}
			X^n_t = x_0 +\! \int_0^t \mu(s, X^n_{s \wedge \cdot}, \nu^{*,n}_s) ds 
			+\!
			\int_0^t \sigma(s, X^n_{s \wedge \cdot}, \nu^{*,n}_s) d B^{*,n}_s,
			~t \ge 0,
			~\Pb^* \mbox{-a.s.}
		\end{equation}
		Then for all $\varphi \in C_c^{\infty}(\R^d)$ and $s \le t$, one has
		$$
			\E^{\Pb^*} \bigg[
				\bigg| \int_s^t \!\!\! \int_U \Lc \varphi(s, X^n_{s \wedge \cdot}, u, X^n_s) \big( m^{*,n}(du,ds) - m^*(du,ds) \big) 
				\bigg| \bigg]
			\longrightarrow 
			0,
			~\mbox{as}~
			n \longrightarrow \infty.
		$$
		In particular, Condition \eqref{eq:gmn_cvg} holds true in this setting.
	\end{Proposition}
	\proof Notice that $(\mu, \sigma)$ are uniformly bounded, so that, by standard arguments,
	$$
		\lim_{\Delta \longrightarrow 0}~ \sup_{n \ge 1} \sup_{s \le u \le t} 
		\E^{\Pb^*} \Big[ \sup_{u \le r \le u+\Delta}  \big|X^n_r - X^n_u \big| \Big] 
		= 0.
	$$
	Further, by Lipschitz property of $\Lc \varphi(t, \xb, u, x)$ in $(\xb, x)$, and its uniform continuity in $t$,
	it follows that, for any $\eps > 0$, there exits a partition $s = t_0 < t_1 < \cdots < t_N = t$ for some $N \ge 1$,
	such that, for all $n \ge 1$, one has
	\be \label{eq:Lc_phi_freeztime}
		&&
		\E^{\Pb^*}
		\bigg[ \bigg| 
			\int_s^t  \!\! \int_U \Lc \varphi(s, X^n_{s \wedge \cdot}, u, X^n_s) (m^{*,n}_s(du) - m^*_s(du)) ds
		\nonumber \\
		&&~~~~~~~~~- 
			\sum_{i=0}^{N-1} \int_{s_i}^{s_{i+1}} \!\!\!\! \int_U \Lc \varphi(s_i, X^n_{s_i \wedge \cdot}, u, X^n_{s_i}) (m^{*,n}_s(du) - m^*_s(du)) ds
		\bigg| \bigg]
		\le \eps.~~~~~~~
	\ee
	At the same time, by \eqref{eq:phi_mn_m},
	one can easily deduce that, as $n \longrightarrow \infty$,
	$$
		\int_{s_i}^{s_{i+1}} \!\!\! \int_U \Lc \varphi(s_i, X^n_{s_i \wedge \cdot}, u, X^n_{s_i}) 
		\big(m^{*,n}_s(du) - m^*_s(du) \big) ds
		\longrightarrow
		0,~\Pb^*\mbox{-a.s.}
	$$
	Notice that $\Lc \varphi$ is uniformly bounded, this is enough to prove that
	$$
		\E^{\Pb^*} \Big[
			\Big|
				\int_s^t \!\! \int_U \Lc \varphi(s, X^n_{s \wedge \cdot}, u, X^n_s) (m^{*,n}_s(du) - m^*_s(du)) ds
			\Big|
		\Big]
		\longrightarrow
		0,
		~\mbox{as}~
		n \longrightarrow \infty,
	$$
	and we hence conclude the proof.
	\qed

	\begin{Remark} \label{rem:uniform_continu_sigma}
		The uniform continuity of $(\mu, \sigma)(s, \xb, u)$ in $s$ is only used in \eqref{eq:Lc_phi_freeztime}.
		In particular, when $\mu$ (resp. $\sigma$) is uncontrolled in the sense that $\mu(t, \xb, u_1) = \mu(t, \xb, u_2)$ (reps. $\sigma(t, \xb, u_1) = \sigma(t, \xb, u_2)$) for all $u_1, u_2 \in U$,
		the uniform continuity of $\mu$ (resp. $\sigma$) in time is then not needed to check  \eqref{eq:gmn_cvg} in the proof of Proposition \ref{prop:qmn_cvg}.	
	\end{Remark}

	We now consider a first case, where the diffusion process $\sigma$ is uncontrolled, to obtain the convergence result.
	
	\begin{Proposition} \label{prop:cvg_un_control_vol}
		Let Assumption \ref{assu:bbd_Lip_coef} hold true.
		Assume in addition that the diffusion coefficient $\sigma$ is uncontrolled in the sense that
		$$
			\sigma(t, \xb, u_1) = \sigma(t, \xb, u_2),
			~\mbox{for all}~
			(t,\xb) \in \R_+ \x \Om
			~\mbox{and}~
			u_1, u_2 \in U.
		$$
		Then there exists a sequence $(\nu^n)_{n \ge 1}$ of $U$-valued $\Fbb$-adapted piecewise constant control processes
		together with a sequence $(\Pb^{*,n})_{n \ge 1}$ of (weak) control rule associated with the control processes $m^{*,n}  (du, ds) := \delta_{\nu^{*,n}_s}(du) ds$,
		such that
		$$
			m^{*,n} \longrightarrow m^*,~\Pb^*\mbox{-a.s. and}~
			\Pb^{*,n} \longrightarrow \Pb^* ~\mbox{under the stable convergence topology.}
		$$
	\end{Proposition}
	\proof Let us apply Theorem \ref{eq:gmn_cvg} to deduce the convergence result.
	In fact, when the volatility coefficient is uncontrolled, one can combine the two Brownian motion $B^1$ and $B^2$ in \eqref{eq:SDE_2BM_relax} into one Brownian motion $B^*$, and rewrite the dynamic of $X$ as
	$$
		X_t = x_0 + \sum_{i=1}^2 
		\int_0^t\!\! \mu(s, X_{s \wedge \cdot}, u_i) q^{i}_s ds 
		+
		\int_0^t  \sigma(s, X_{s \wedge \cdot}) dB^*_s,
		~t \ge 0,
		~\Pb^* \mbox{-a.s.}
	$$
	To apply Theorem \ref{eq:gmn_cvg}, we need to consider an enlarged space $\Omh^* := \Om^* \x \Om \x \Om$, with canonical process $(X, B)$, and $\Ph^* (d\om, d\xb, d \bb) := \delta_{B^*(\om, \xb)}(d\bb) \Pb^*(d \om, d\xb)$.
	Namely, one has $B = B^*$, $\Ph^*$-a.s.
	We then consider the generator $\widehat \G$ of the couple $(X, B)$, defined by
	$$
		\widehat{\G} ~:=~ \big\{ (\varphi, \widehat \Lc \varphi ~: \varphi \in C^{\infty}_c(\R^d \x \R^d) \big\},
	$$
	where, with $\hat \sigma^{\top} (\cdot) := (\sigma^{\top}(\cdot), I_d)$, 
	$$
		\widehat \Lc \varphi (s, \xb, \bb, u, x, b) 
		:=
		\mu(s, \xb, u) \cdot D_x \varphi(x,b)
		+
		\frac12 \mathrm{Tr} \big( \hat \sigma \hat \sigma^{\top}(s, \xb, u) D^2 \varphi(x,b) \big).
	$$
	
	For each $n \ge 1$, let $\Ph^{*,n}$ be the weak control rule associated with the generator $\widehat{\G}$ and the control process $\nu^{*,n}$.
	Let us impose the following additional condition on $\Ph^{*,n}$:
	$$
		\Ph^{*,n} \big[ B = B^* \big] = 1.
	$$
	Then it is easy to check that the functionals in the generator $\widehat{\G}$ satisfies the required continuity condition,
	and one has $m^{*,n}(du, ds) \longrightarrow m^*$, $\Ph^*$-a.s.
	Moreover, as $\mu$ and $\sigma$ are uniformly bounded, one can easily check that $(\Ph^{*,n})_{n \ge 1}$ is relatively compact,
	so that for any subsequence, one can subtract a further subsequence such that
	$$
		\Ph^{*,n_k} 
		\longrightarrow
		\Ph^{*, \infty}.
	$$
	Further, it follows by Proposition \ref{prop:qmn_cvg} that \eqref{eq:gmn_cvg} holds true in the $\widehat{\G}$ generator setting.
	Therefore, one can apply Theorem \ref{eq:gmn_cvg} to deduce that $\Ph^{*, \infty}$ is a relaxed control rule with generator $\widehat{\G}$ and the weak control process $m^*$.
	In particular, one has
	$$
		X_t = x_0 + \sum_{i=1}^2 
		\int_0^t\!\! \mu(s, X_{s \wedge \cdot}, u_i) q^{i}_s ds 
		+
		\int_0^t  \sigma(s, X_{s \wedge \cdot}) dB_s,
		~t \ge 0,
		~\Ph^{*,\infty} \mbox{-a.s.}
	$$
	Notice that $\Ph^{*,\infty}[B = B^*] = 1$ as $\Ph^{*,n}[B = B^*] = 1$ for all $n \ge 1$, and $\Ph^{*,\infty}|_{\Om^*} = \Ph^*|_{\Om^*} = \P^*$.
	One can then deduce that $\Ph^{*,\infty} = \Ph^*$, which concludes the proof.
	\qed

	\begin{Remark}
		In view of Remark \ref{rem:uniform_continu_sigma},
		in the uncontrolled volatility coefficient context in Proposition \ref{prop:cvg_un_control_vol},
		the same result hold still true if the uniform continuity in time variable property in Assumption \ref{assu:bbd_Lip_coef} is only assumed on $\mu$ (but not on $\sigma$).
	\end{Remark}

	We now consider the general case, where both drift and diffusion coefficient functions could be controlled.
	For this purpose, with the fixed Brownian motions $B^1$ and $B^2$ in \eqref{eq:SDE_2BM_relax}, 
	let us construct a new Brownian motion $B^{*,n}$.
	For each $n \ge 1$, let $B^{*,n}_0 := 0$;
	and for each $i = 0, \cdot, n-1$, given the value $B^{*,n}_{t_{k,i}}$, we define $B_t^{*,n}$ for $t \in [t_{k,i}, t_{k,i+1}]$ as follows:
	\b*
		B^{*,n}_t - B^{*,n}_{t_{k,i}} 
		&\!\!\!:=\!\!\!\!&
		\sqrt{q^1_{t_k}} \big( B^1_{t_{k,i} + (t-t_{k,i})(t_{k,i+1}- t_{k,i})/( \theta_{k,i} - t_{k,i}) } \!\!-\! B^1_{t_{k,i}} \big) 
			\1_{\{t \in [t_{k,i}, \theta_{k,i})\}} \\
		&&\!\! + 
		\sqrt{q^2_{t_k}}  \big(B^2_{t_{k,i+1} - (t_{k,i+1} - t) (t_{k,i+1} - t_{k,i})/(t_{k,i+1} -\theta_{k,i})} \!\!-\! B^2_{t_{k,i}} \big) 
			\1_{\{t \in [\theta_{k,i}, t_{k,i+1}] \}}.
	\e*
	Namely, we compress the increment of $B^1$ on interval $[t_{k,i}, t_{k,i+1})$ into a martingale on $ [t_{k,i}, \theta_{k,i})$,
	and compress the the increment of $B^2$ on $[t_{k,i}, t_{k,i+1})$ into a martingale on $ [\theta_{k,i}, t_{k,i+1}]$,
	and then paste and renormalize them to obtain the increment of $B^{*,n}$ on $[t_{k,i}, t_{k,i+1}]$.
	Although $B^{*,n}$ is not adapted to the filtration of $(B^1, B^2)$, 
	but it is a standard Brownian motion w.r.t. the filtration generated by $(B^{*,n}, \nu^{*,n})$.
	One can then define $X^n$ by
	$$
			X^n_t = x_0 +\! \int_0^t \mu(s, X^n_{s \wedge \cdot}, \nu^{*,n}_s) ds 
			+\!
			\int_0^t \sigma(s, X^n_{s \wedge \cdot}, \nu^{*,n}_s) d B^{*,n}_s,
			~t \ge 0,
			~\Pb^* \mbox{-a.s.}
	$$
	For later uses, we also fix a $\F^*$-stopping time $\tau^*$ taking value in $\R_+$.
	
	\begin{Proposition} \label{prop:approximation_general}
		Let Assumption \ref{assu:bbd_Lip_coef} hold true.
		Then there exists a sequence $(\tau^{*,n})_{n \ge 1}$ of stopping times
		such that, with $m^{*,n}$ and $X^n$ be defined in \eqref{eq:def_mn_star} and \eqref{eq:def_Xn},
		\begin{equation} \label{eq:cvg_law_tmX}
			\big( \tau^{*,n}, m^{*,n}, X^n \big)
			\longrightarrow
			\big( \tau^*, m^*, X \big),
			~\Pb^*\mbox{-a.s., as}~
			n \longrightarrow \infty.
		\end{equation}
	\end{Proposition}
	\proof First, let us take a time discretization parameter $\Delta > 0$, and define the corresponding time freezing function $\eta_{\Delta}: \R_+ \longrightarrow \R_+$ by $\eta_{\Delta}(t) := k \Delta$ for all $t \in [k \Delta, (k+1) \Delta)$, $k = 0, 1, \cdots$.
	We then introduce a $\F^*$-stopping time $\tau^{*,\Delta}$ by
	$$
		\tau^{*, \Delta} := \sum_{k=1}^{\infty} (k+1) \Delta \1_{\{\tau^* \in (k-1)\Delta, k \Delta]\}},
	$$
	and observe that $\big | \tau^* - \tau^{*,\Delta} \big| \le 2 \Delta$, $\Pb^*$-a.s.
	Let us also define $X^{\Delta}$ and $X^{n,\Delta}$ by
	$$
		X^{\Delta}_t = x_0 + \sum_{i=1}^2 
		\Big(\!
		\int_0^t\!\! \mu \big(\eta_{\Delta}(s) , \widehat{X^{\Delta}}, u_i \big) q^{i}_s ds 
		+
		\int_0^t \!\! \sigma \big(\eta_{\Delta}(s) , \widehat{X^{\Delta}} , u_i \big) \sqrt{q^{i}_s} dB^{i}_s
		\Big),
		~t \ge 0, ~\Pb^*\mbox{-a.s.},
	$$
	and
	$$
		X^{n,\Delta}_t = x_0 +\! \int_0^t \!\! \mu(\eta_{\Delta}(s) , \widehat{X^{n,\Delta}}, \nu^{*,n}_s) ds 
		+\!
		\int_0^t \!\! \sigma(\eta_{\Delta}(s) , \widehat{X^{n,\Delta}}, \nu^{*,n}_s) d B^{*,n}_s,
		~t \ge 0,
		~\Pb^*\mbox{-a.s.},
	$$
	where $\widehat{X^{\Delta}}$ (resp. $\widehat{X^{n,\Delta}}$) denotes the continuous time process obtained from the linear interpolation of $(X^{\Delta}_{k\Delta})_{k \ge 0}$ (resp. $(X^{n, \Delta}_{k\Delta})_{k \ge 0}$.
	Namely, without taking into account the control process, $X^{\Delta}$ and $X^{n,\Delta}$ can be considered as Euler scheme of $X$ and $X^n$ respectively.
	As in the numerical analysis of the simulation of SDEs (see e.g. Graham and Talay \cite{GrahamTalay}), 
	for every $T > 0$, one has
	\begin{equation} \label{eq:EulerScheme_rate}
		\lim_{\Delta \longrightarrow 0}
		~\sup_{n \ge 1} 
		~\E \Big[ \sup_{0 \le t \le T} \Big( \big| X_t - \widehat{X^{\Delta}_t} \big|^2 + \big| X^n_t - \widehat{X^{n,\Delta}_t} \big|^2 \Big) \Big]
		=
		0.
	\end{equation}

	We next consider the processes $X^{\Delta}$ and $X^{n, \Delta}$ on a time interval $[\ell \Delta, (\ell+1) \Delta]$.
	For $\Delta > 0$ small enough, and take $n \ge 1$ large enough, 
	one can assume without loss of generality that 
	$$
		\{t_k~:k \ge 0\} ~\subset \{ \ell \Delta ~:\ell \ge 0\} ~\subset~ \{ t_{k,i} ~:i = 0, \cdots, n-1, ~k \ge 0 \},
	$$
	where we recall that $(t_k)_{k \ge 0}$ is the discrete time grid on which the relaxed control process $m^*$ is piecewise constant.
	Then, on each time interval $[t_{k,i}, t_{k,i+1}]$, the drift and volatility coefficients of $X^{\Delta}$ and $X^{n, \Delta}$,
	and the control processes are all frozen.
	At the same time, with the definition of $B^{*,n}$, one can easily check that
	$$
		\sum_{i=1}^2 \int_{t_{k,i}}^{t_{k,i+1}} q^i_s ds = t_{k,i+1} - t_{k,i}
		~~~\mbox{and}~~
		\sum_{i=1}^2 \int_{t_{k,i}}^{t_{k,i+1}} \sqrt{q^i_s} dB^i_s
		=
		\int_{t_{k,i}}^{t_{k,i+1}} d B^{n,*}_s,
		~\Pb^*\mbox{-a.s.}
	$$
	This implies that, for $\Delta > 0$ small enough, and then $n\ge 1$ large enough, one has
	$$
		\widehat{X^{\Delta}} = \widehat{X^{n,\Delta}},~\Pb^*\mbox{-a.s.}
	$$
	Together with \eqref{eq:EulerScheme_rate}, and by a simple diagonalization argument, one can then conclude the proof.
	\qed

	\begin{Remark}
		In \cite{ElKaroui_1987}, the authors considered directly the weak limit of $(X^n, B^{*,n})$,
		and proved a weak convergence result of $(m^{*,n}, X^n)$ to $(m^*, X)$.
		The convergence result in Proposition \ref{prop:approximation_general} is in the sense of a.s.
		In particular, Proposition \ref{prop:approximation_general} includes the convergence of the stopping time,
		which would be useful to study the mixed control/stopping problems.
	\end{Remark}

	\begin{Remark} \label{rem:approximation_lsc}
		Let $L: \R_+ \x \Om \x U \longrightarrow \Rb$
		be such that $L(t, \xb, u) = L(t, \xb_{t \wedge \cdot}, u)$, for all $(t,\xb, u) \in \Rb_+ \x \Om \x U$.
		Let $\tau^*$ be a $\F^*$-stopping time taking value in $[0,T]$.
		
		\vspace{0.5em}
		
		\noindent \rmi In the context of Proposition \ref{prop:cvg_un_control_vol}, 
		where $m^{*,n} \longrightarrow m^*$, $\P^*$-a.s. and $\Pb^{*,n} \longrightarrow \Pb^*$,
		one can then apply similar arguments as in Proposition \ref{prop:qmn_cvg} to deduce that,
		when $L$ is uniformly bounded and uniformly continuous in all its arguments,
		\begin{align*}
			\lim_{n \longrightarrow \infty}
			\E^{\Pb^{*,n}} \Big[ \int_0^{\tau^*} \!\!\! L(s, X, \nu^{*,n}_s) ds  \Big]
			=
			\E^{\Pb^{*}} \Big[ \int_0^{\tau^*} \!\!\!\! \int_U L(s, X, u) m^{*}(du, ds)  \Big].
		\end{align*}
		When $L$ is bounded from below and is lower semi-continuous, there is a sequence of Lipschitz functions $(L_k)_{k \ge 1}$ such that $L_k \nearrow L$ point-wise.
		Thus
		\begin{align*}
			&\liminf_{n \longrightarrow \infty}
			\E^{\Pb^{*,n}} \Big[ \int_0^{\tau^*} \!\!\! L(s, X, \nu^{*,n}_s) ds  \Big]
			~\ge~
			\lim_{k \longrightarrow \infty}
			\liminf_{n \longrightarrow \infty}
			\E^{\Pb^{*,n}} \Big[ \int_0^{\tau^*} \!\!\! L_k(s, X, \nu^{*,n}_s) ds  \Big]\\
			=&
			\lim_{k \longrightarrow \infty}
			\E^{\Pb^{*}} \Big[ \int_0^{\tau^*} \!\!\!\! \int_U L_k(s, X, u) m^{*}(du, ds)  \Big]
			~=~
			\E^{\Pb^{*}} \Big[ \int_0^{\tau^*} \!\!\!\! \int_U L(s, X, u) m^{*}(du, ds)  \Big].
		\end{align*}
		
		\noindent \rmii Let us stay in the context of Proposition \ref{prop:approximation_general},
		where $(\tau^{*,n}, m^{*,n}, X^n) \longrightarrow (\tau^*, m^*, X)$, $\Pb^*$-a.s.,
		one can deduce similarly that, 
		when $L$ is uniformly bounded and uniformly continuous in all its arguments,
		$$
			\lim_{n \longrightarrow \infty} \E^{\Pb^{*}} \Big[ \int_0^{\tau^{*,n}} \!\!\!  L(s, X^n, \nu^{*,n}_s) ds \Big] 
			=
			\E^{\Pb^{*}} \Big[ \int_0^{\tau^*} \!\!\!\! \int_U L(s, X, u) m^{*}(du, ds) \Big].
		$$
		When $L$ is  bounded from below and is lower semi-continuous, by the same arguments as above, one has
		$$
			\liminf_{n \longrightarrow \infty} \E^{\Pb^{*}} \Big[ \int_0^{\tau^{*,n}} \!\!\!  L(s, X^n, \nu^{*,n}_s) ds \Big] 
			\ge
			\E^{\Pb^{*}} \Big[ \int_0^{\tau^*} \!\!\!\! \int_U L(s, X, u) m^{*}(du, ds) \Big].
		$$
	\end{Remark}

\subsection{Equivalence of the optimal stopping problems}
\label{subsec:equiv_opt_stop}

\paragraph{On the canonical space}

	Recall that $\Om = \D(\R_+, \E)$ denotes the canonical space, with canonical process $X$ and canonical filtration $\F = (\Fc_t)_{t \ge 0}$.
	Let $\P$ be a fixed probability space, so that $(\Om, \Fc_{\infty}, \P)$ is a fixed probability space,
	we denote by $\F^{\P}$ the completed filtration and by $\F^{\P}_+ = (\Fc^{\P}_{t+})_{t \ge 0}$ the augmented filtration;
	denote also by $\Tc^{\P}$ (resp. $\Tc^{\P}_+$) the class of all  $\F^{\P}$ (resp. $\F^{\P}_+$) -stopping times.
	Let $\tau \in \Tc^{\P}_+$, then the couple $(\tau, X)$ induces a probability measure on $\Omh := \Rb_+ \x \Om $.
	We hence consider the enlarged canonical space $\Omh$, with canonical element $(\Theta_{\infty}, X)$ and canonical filtration $\Fh = (\Fch_t)_{t \ge 0}$ with $\Fch_t := \sigma( X_s, \Theta_s, s \in [0,t] \cap \R_+)$, with $\Theta_s := \Theta_{\infty} \wedge s$, for $t \in \Rb_+$.
	Denote by by $\Fh^X = (\Fch^X_t)_{t \ge 0}$ the filtration generated by $X$ on $\Omh$,
	and
	$$
		\Pch_0 ~:=~ \Big\{ 
			\Ph ~: \Ph \big|_{\Om} = \P,
			~\mbox{and}~ \E^{\Ph} \big[\1_{\Theta_{\infty} \le t} \big| \Fch^X_{\infty} \big]
			= \E^{\Ph} \big[\1_{\Theta_{\infty} \le t} \big| \Fch^X_t \big], ~\Ph \mbox{-a.s.}~\forall t \ge 0
		\Big\},
	$$
	and
	$$
		\Pch^+_0 ~:=~ \Big\{ 
			\Ph ~: \Ph \big|_{\Om} = \P, 
			~\mbox{and}~ \E^{\Ph} \big[\1_{\Theta_{\infty} \le t} \big| \Fch^X_{\infty} \big]
			= \E^{\Ph} \big[\1_{\Theta_{\infty} \le t} \big| \Fch^X_{t+} \big], ~\Ph \mbox{-a.s.}~\forall t \ge 0
		\Big\}.
	$$
	\begin{Proposition} \label{prop:eqiv_opt_stopping}
		Let $\Phi: \Omh \to \R$ satisfy $\Phi(t, \om) = \Phi(t, \om_{t \wedge \cdot})$ for all $(t, \om) \in \Omh$.
		We then have the equivalence of the two different formulations of the optimal stopping problem
	$$
		\sup_{\tau \in \Tc^{\P}} \E^{\P} \big[ \Phi(\tau, X) \big]
		=
		\sup_{\Ph \in \Pch_0} \E^{\Ph} \big[ \Phi(\Theta_{\infty}, X) \big],
		~
		\sup_{\tau \in \Tc_+^{\P}} \E^{\P} \big[ \Phi(\tau, X) \big]
		=
		\sup_{\Ph \in \Pch_0^+} \E^{\Ph} \big[ \Phi(\Theta_{\infty}, X) \big].
	$$	
	\end{Proposition}
	\proof We only prove the first equivalence, the second follows by the same arguments.
	
	\vspace{0.5em}

	\noindent \rmi Let $\tau \in \Tc^{\P}$ be a $\F^{\P}$-stopping time,
	then it is clear that, under $\P$, $(\tau, X)$ induces a probability measure in $\Pch_0$,
	we then have a first inequality 
	$
		\sup_{\tau \in \Tc^{\P}} \E^{\P} \big[ \Phi(\tau, X_{ \cdot}) \big]
		\le
		\sup_{\Ph \in \Pch_0} \E^{\Ph} \big[ \Phi \big(\Theta_{\infty}, X_{ \cdot} \big) \big].
	$
	
	\vspace{0.5em}
	
	\noindent \rmii
	Next, let $\Ph \in \Pch_0$, we denote by $(\Ph_{\om})_{\om \in \Om}$ a family of conditional probability measures of $\Ph$ w.r.t. $\Fch^X_{\infty}$,
	and denote $F_{\om}(t) := \Ph_{\om} \big[ \Theta \le t\big]$, which is right-continuous and $\F^{\P}$-adapted since for any $t \ge 0$,
	$$
		F_{\om}(t) 
		~=~
		\E^{\Ph} \big[ \1_{\Theta_{\infty} \le t} \big| \Fch^X_{\infty} \big] (\om)
		~=~
		\E^{\Ph} \big[ \1_{\Theta_{\infty} \le t} \big| \Fch^X_t \big] (\om),
		~~\mbox{for}~\P\mbox{-a.e.}~ \om \in \Om.
	$$
	Denote by $F^{-1}_{\om}: [0,1] \to \Rb_+$ the right-continuous inverse function of $x \mapsto F_{\om}(x)$,
	it follows that for any $u \in [0,1]$, one has $\{ \om ~: F^{-1}_{\om}(u) \le t \} = \{ \om ~: F_{\om}(t) \le u \} \in \Fc^{\P}_{t} $, and hence $\om \mapsto F^{-1}_{\om} (u)$ is a $\F^{\P}$-stopping time.
	Therefore, one obtains the inequality
	\b*
		&& \E^{\Ph} \big[ \Phi(\Theta_{\infty}, X \big] 
		~=~ 
		\E^{\Ph} \Big[ \E^{\Ph} \big[ \Phi( \Theta_{\infty}, X) \big| \Fcb^X_{\infty} \big] \Big]
		~=~
		\E^{\P} \Big[ \int_{[0,{\infty}]} \Phi(\theta, X) F_X (d \theta) \Big]\\
		&=&
		\E^{\P} \Big[ \int_{[0,1]} \Phi \big( F^{-1}_X(z), X \big) dz \Big]
		~\le~
		\sup_{\tau \in \Tc^{\P}} \E^{\P} \big[ \Phi(\tau, X) \big].
	\e*
	Together with the inequality in Item $\mathrm{(i)}$, it concludes the proof.
	\qed

	\begin{Remark} \label{rem:PropertyK}
		Suppose that, in the filtered probability space $(\Om, \Fc_{\infty}^{\P}, \F^{\P}_+, \P)$, $X$ is a Markov process; and let $\Ph$ be a probability measure on $\Omh$ under which $X$ is still a Markov process w.r.t. $\Fh^{\Ph}_+$ with the same generator.
		Then it is easy to check that $\Ph \in \Pch_0^+$.
	\end{Remark}

\paragraph{A more general equivalence result}

	The above condition 
	$\E^{\Ph} \big[\1_{\Theta_{\infty} \le t} \big| \Fch^X_{\infty} \big]
	= \E^{\Ph} \big[\1_{\Theta_{\infty} \le t} \big| \Fch^X_t \big]$ 
	is also called Property (K) in the context of optimal control/stopping problems, or called Hypothesis (H) in the context of filtration enlargement problems.
	It can be formulated in a more abstract context, where the above equivalence result holds still true.
	Let $(\Om^*, \Fc^*, \P^*, \F^*)$ be a filtered probability space, where the filtration $\F^* = (\Fc^*_t)_{t \ge 0}$ satisfies the usual conditions.
	Denote by $\Tc^*_{\infty}$ the class of all finite $\F^*$-stopping times.
	Further, let $\G^* = (\Gc^*_t)_{t \ge 0}$ be another filtration satisfying the usual conditions,
	such that $\Gc^*_t \subseteq \Fc^*_t$ for all $t \ge 0$,
	we denote by $\Tc^*_{\infty}(\G^*)$ the collection of all finite $\G^*$-stopping times.
	A reward process $Y$ is assumed to be $\G^*$-optional, l\`adl\`ag, and of class (D),
	we then have the following equivalence result by Szpirglas and Mazziotto \cite{SzpirglasMazziotto}.

	\begin{Theorem}
		Suppose that the filtered probability space $(\Om^*, \Fc^*, \P^*, \F^*, \G^*)$ satisfies Property (K), i.e.
		for all $t \ge 0$ and all $\Fc^*_t$-measurable bounded random variable $\xi$,
		$$
			\E^{\P^*} \big[ \xi \big| \Gc^*_{\infty} \big] 
			~=~
			\E^{\P^*} \big[ \xi \big| \Gc^*_t \big],
			~~\P^*-a.s. 			
		$$
		Then, one has the equivalence of the following two optimal stopping problems:
		$$
			\sup_{\tau \in \Tc^*_{\infty}} ~\E^{\P^*} \big[ Y_{\tau} \big]
			~~=~~
			\sup_{\tau \in \Tc^*_{\infty}(\G^*)} \E^{\P^*} \big[ Y_{\tau} \big].
		$$
	\end{Theorem}

\subsection{Equivalence of the controlled/stopped diffusion processes problems}

	Let us stay in the context of the controlled/stopped diffusion processes problem as presented in Section \ref{subsec:diffusion},
	and study the equivalence of different formulations of the problem.
	Recall that, in this context, one has $E \equiv \R^d$,
	and one is given the drift and diffusion coefficient functions $(\mu, \sigma) : \R_+ \x \Om \x U \longrightarrow \R^d \x \S^d$,
	satisfying Assumption \ref{assu:bbd_Lip_coef}.
	We will consider a pure control problem, where the reward functions are given by $L: \R_+ \x \Om \x U \longrightarrow \R$ and $\Phi_1: \Om \longrightarrow \R$,
	and also a mixed control/stopping problem, where  the reward function is given by $\Phi_2: \Rb_+ \x \Om \longrightarrow \R$.
	Moreover, let us assume that $L(t, \xb, u) = L(t, \xb_{t \wedge \cdot},u)$ and $\Phi_2(t, \xb) = \Phi_2(t, \xb_{t \wedge \cdot})$ for all $(t,\xb, u) \in \R_+ \x \Om \x U$.
	
	\vspace{0.5em}
	
	Let us recall quickly from Section \ref{subsec:diffusion} the strong, weak and relaxed formulations of the controlled/stopped diffusion processes problem.
	First, in a probability space $(\Om^*,\Fc^*, \P^*)$ equipped with a Brownian motion $B$ and the Brownian filtration $\F^* = (\Fc^*_t)_{t \ge 0}$,
	we denote by $\Tc$ the collection of all $\F^*$-stopping times.
	Let us denote by $\Uc$ the collection of all $U$-valued $\F^*$-predictable process, and by $\Uc_0$ the subset of all piecewise constant control processes $\nu \in \Uc$.
	Then given a control process $\nu \in \Uc$, $X^{\nu}$ is the corresponding controlled process defined as the unique strong solution to SDE \eqref{eq:controlledSDE} with a fixed initial condition $x_0 \in \R^d$.
	Let us define the value of the strong formulation of the control or control/stopping problem by
	\begin{equation} \label{eq:def_VS0}
		V_1^{S}
		~:=~
		\sup_{\nu \in \Uc}
		\E \Big[ \int_0^{\infty} L(s, X^{\nu}, \nu_s) ds + \Phi_1 \big( X^{\nu} \big) \Big],
	\end{equation}
	and
	\begin{equation}
		V_2^{S}
		~:=~
		\sup_{\nu \in \Uc} ~ \sup_{\tau \in \Tc}
		\E \Big[ \int_0^{\tau} L(s,X^{\nu}, \nu_s)ds + \Phi_2 \big(\tau, X^{\nu} \big) \Big].
	\end{equation}
	Next, without fixing the probability space and the filtration, the set $\Ac_W$ of weak controls $\alpha = (\Om^{\alpha}, \Fc^{\alpha}, \P^{\alpha}, \F^{\alpha}, \tau^{\alpha}, X^{\alpha}, B^{\alpha}, \nu^{\alpha})$ 
	and the set $\Ac_R$ of relaxed controls $\alpha = (\Om^{\alpha}, \Fc^{\alpha}, \P^{\alpha}, \F^{\alpha}, \tau^{\alpha}, X^{\alpha}, M^{\alpha}, \widehat M^{\alpha})$ are given in Definitions \ref{def:weak_control} and \ref{def:relaxed_control}.
	Let us denote by $\Ac_{W_0}$ the subset of weak controls $\alpha \in \Ac_W$ such that $\nu^{\alpha}$ is piecewise constant.
	We then obtain the value of the weak formulation of the control, or control/stopping problem:
	\begin{equation} \label{eq:def_VW}
		V^W_1 := \sup_{\alpha \in \Ac_W} \E^{\Pb^{\alpha}} \Big[ \int_0^{\infty} L(s, X^{\alpha}, \nu^{\alpha}_s) ds + \Phi_1(X^{\alpha}) \Big],
	\end{equation}
	and
	\begin{equation}
		V^W_2 := \sup_{\alpha \in \Ac_W} \E^{\Pb^{\alpha}} \Big[ \int_0^{\tau^{\alpha}} L(s, X^{\alpha}, \nu^{\alpha}_s) ds + \Phi_2(\tau^{\alpha}, X^{\alpha}) \Big].
	\end{equation}
	Similarly, one has the value of the relaxed formulation of the control, or control/stopping problem:
	$$
		V^R_1 := \sup_{\alpha \in \Ac_R} \E^{\Pb^{\alpha}} \Big[ \int_0^{\infty} \int_U L(s, X^{\alpha},u)  M^{\alpha}_s(du) ds + \Phi_1(X^{\alpha}) \Big],
	$$
	and
	$$
		V^R_2 := \sup_{\alpha \in \Ac_R} \E^{\Pb^{\alpha}} \Big[ \int_0^{\tau^{\alpha}} \int_U L(s, X^{\alpha}, u) M^{\alpha}_s(du) ds + \Phi_2(\tau^{\alpha}, X^{\alpha}) \Big].
	$$
	Finally, replacing $\Uc$ by $\Uc_0$ in the definition of $V_1^{S}$ and $V_2^{S}$,
	and replacing $\Ac_W$ by $\Ac_{W_0}$ in the definition of $V_1^W$ and $V_2^W$, one defines similarly
	$$
		V_1^{S_0}, ~~V_2^{S_0},~~ V_1^{W_0} ~\mbox{and}~V_2^{W_0}.
	$$

	Our main result in this part is then the following equivalence of different formulations of the controlled/stopped diffusion processes problem.

	\begin{Theorem} \label{theo:equiv_ctrl}
		\rmi Let Assumption \ref{assu:bbd_Lip_coef} hold true.
		Then 
		$$
			V^{S_0}_1
			~=~
			V^{W_0}_1
			~~\mbox{and}~~
			V^{S_0}_2
			~=~
			V^{W_0}_2.
		$$

		\noindent \rmii Assume in addition that $L$, $\Phi_1$ and $\Phi_2$ are all lower-semicontinuos and bounded from below.
		Then one has the equivalence
		$$
			V^{S_0}_i
			~=~
			V^{W_0}_i
			~=~
			V^S_i
			~=~
			V^W_i
			~=~
			V^R_i,
			~i=1,2.
		$$
	\end{Theorem}
	
	\begin{Remark}

		\noindent \rmi One can relax the boundedness condition on $(\mu, \sigma)$ in Assumption \ref{assu:bbd_Lip_coef} by truncating unbounded coefficient functions.
		In particular, in the context where one replaces the boundedness condition on $(\mu, \sigma)$  in Assumption \ref{assu:bbd_Lip_coef} by \eqref{eq:mu_sigma_integ_p},
		or in the context of Section \ref{subsubsec:unbound_mu_sigma} with integrability conditions on control processes,
		one can consider the optimal control/stopping with truncated coefficient functions $(\mu_n, \sigma_n) := ( -n \vee \mu_i \wedge n, -n \vee \sigma_{i,j} \wedge n)_{1 \le i,j \le d}$.
		Next, by considering the corresponding values of the control/stopping problem with index $n \ge 1$, 
		and under mild conditions on $L$ and $\Phi_i$,
		one can show the convergence
		$$
			\big(V^{S_0}_{i,n}, V^{W_0}_{i,n}, V^S_{i,n}, V^W_{i,n}, V^R_{i,n}\big)
			\longrightarrow
			\big(V^{S_0}_{i}, V^{W_0}_{i}, V^S_i, V^W_i, V^R_i\big),
			~i=1,2,
		$$
		and then obtains the same equivalence results.
		
		\vspace{0.5em}
		
		\noindent \rmii The boundedness from below condition in $\mathrm{(i)}$ is only used to apply the approximation results in Propositions \ref{prop:approximate_q} and \ref{prop:approximation_general},
		in order to show that $V^{W_0}_i = V^R_i$, $i=1,2$. 
		This boundedness condition on $L$, $\Phi_1$ and $\Phi_2$ can be replaced by some uniform integrability conditions so that the approximation argument still works,
		and the equivalence result hold still true.
	\end{Remark}

	Let us first provide a technical lemma.
	Let $\alpha^* = (\Om^*, \Fc^*, \P^*, \F^*, X^*,  \tau^*, \nu^*, W^*) $ be a weak control with piecewise constant control process $\nu^*$, i.e. $\alpha^* \in \Ac_{W_0}$,
	For simplicity and without loss of generality, 
	we assume that $\Om^*$ is a metric space and $\Fc^*$ is its Borel $\sigma$-field,
	and $\nu^*$ is piecewise constant over a deterministic time grid $0 = t_0 < t_1 < \cdots < t_n < \infty$, 
	so that
	$\nu^*_s = u^*_i$ for $s \in (t_i, t_{i+1}]$, where $u^*_i$ is a $\Fc^*_{t_i}$-measurable random variable.
	Further, let us enlarge the space $\Om^*$ to $\Omt^* := \Om^* \x [0,1]^{n+1}$, on which we obtain an independent sequence of i.i.d. random variables $(Z_k)_{ 0 \le k \le n}$ of uniform distribution on $[0,1]$.
	Let us denote the enlarged probability space by $(\Omt^*, \Fct^*, \Pt^*)$.

	\begin{Lemma} \label{lemma:cst_weak_ctrl_repres}
		There are measurable functions $(\Psi_k)_{0 \le k \le n-1}$ 
		($\Psi_k : C([0,t_k], \R^d) \x [0,t_k] \x [0,1]^{k+1} \to U$) 
		such that for every $0 \le k \le n$,
		\be \label{eq:cst_control_rcpd}
			&& \Pt^* \circ \Big( \tau^*, ~W^*_{[0,t_k]}, ~\big(\Psi_i(W_{[0,t_i]}, \tau^* \wedge t_i, Z_0, \ldots, Z_i) \big)_{0 \le i \le k} \Big)^{-1} \nonumber \\
			&=&
			\P^* \circ \Big( \tau^*, ~W^*_{[0,t_k]}, ~\big(u^*_i \big)_{0 \le i \le k} \Big)^{-1}.
		\ee
	\end{Lemma}
	\proof First, we suppose that $U \subseteq [0,1]$ without loss of generality, since  any Polish space is isomorphic to a Borel subset of $[0,1]$.
	Let $x \mapsto F_0(x)$ be the cumulative distribution function of $u_0^{*}$ 
	and $F^{-1}_0$ be its inverse function. 
	It follows that \eqref{eq:cst_control_rcpd} holds true in the case $k =0$ with 
	$\Psi_0(\0, 0, z) := F^{-1}_0( z)$.

	\vspace{0.5em}

	Next, let us prove the lemma by induction. 
	Suppose that \eqref{eq:cst_control_rcpd} holds true for some $k < n$ with measurable functions $(\Psi_i)_{0 \le i \le k}$, we shall show that it is also true for the case $k+1$. 
	Let $\big( \P^{*}_{(\xb,s, u)} ~: (\xb,s, u) \in C([0,t_{k+1}], \R^d) \x \big([0, t_{k+1}] \big) \x U^{k+1} \big)$ be a family of regular conditional distribution probability of $\P^{*}$ 
	w.r.t. the $\sigma-$field generated by $W^{*}_{[0,t_{k+1}]}$, $\tau^* \wedge t_{k+1}$ and $(u^{*}_i)_{0 \le i \le k}$, 
	and denote by $F_{k+1}(\xb, s, u, x)$ the cumulative distribution function of $u^{*}_{k+1}$ under $\P^{*}_{(\xb,s, u)}$. 
	Let $F_{k+1}^{-1}(\xb, s, u, x)$ be the inverse function of $x \mapsto F_{k+1}(\xb,s, u, x)$ and
	\b*
		\Psi_{k+1}(\xb, s, x_0, \cdots, x_k, z) 
		~:=~ 
		F_{k+1}^{-1}\big(\xb, s, \Psi_0(x_0), \cdots, \Psi_k(\xb, s \wedge {t_k}, x_0, \cdots, x_{k-1}), z \big).
	\e*
	One can check that \eqref{eq:cst_control_rcpd} holds still true for the case $k+1$ with the given $(\Psi_i)_{0 \le i \le k}$ and $\Psi_{k+1}$ defined above, and we hence conclude the proof. \qed

	\vspace{0.5em}

	\noindent {\sc Proof of Theorem \ref{theo:equiv_ctrl}}
	We will only prove the equality between $V^{S_0}_2$, $V^{W_0}_2$, $V^S_2$, $V^W_2$ and $V^R_2$, while the other equivalence follows in the same (but easier) way.
	
	\vspace{0.5em}
	
	\noindent \rmi 
	Let us fix an arbitrary weak control
	$\alpha^* = (\Om^*, \Fc^*, \P^*, \F^*, X^*,  \tau^*, \nu^*, W^*) $ with piecewise constant control process $\nu^*$, i.e. $\alpha^* \in \Ac_{W_0}$,
	so that one can construct the functionals $(\Psi_k)_{0 \le k \le n-1}$ as in Lemma \ref{lemma:cst_weak_ctrl_repres}.
	Following the notations therein,
	in the probability space $(\Omt^*, \Fct^*, \Pt^*)$, 
	let us define $\tilde \nu^*_s := \tilde u^*_k$ for all $s \in (t_k, t_{k+1}]$ with $\tilde u^*_k = \Psi_k(W^*_{[0,t_k]}, t_k, Z_0, \cdots, Z_k)$,
	and a process $\widetilde X^*$ by
	\be \label{eq:def_X0}
		\widetilde X^*_t 
		~=~ 
		\int_0^t \mu(s,\widetilde X^*_{s \wedge \cdot}, \tilde \nu^*_s) ds  
		~+~ 
		\int_0^t \sigma(s, \widetilde X^*_{s \wedge \cdot}, \tilde \nu^*_s) dW^*_s,
		~~\Pt^* \mbox{-a.s.}
	\ee
	Notice that the law $\P^* \circ \big (\tau^*, W^*, \1_{s \le \tau^*} \delta_{\nu^*_s}(du) ds \big) ^{-1} 
	= \Pt^* \circ \big(\tau^*, W^*, \1_{s \le \tau^*} \ \delta_{\tilde \nu^*_s}(du) ds \big)^{-1}$,
	then $\P^* \circ \big(\tau^*, X^*_{\tau^* \wedge \cdot} \big)^{-1} = \Pt^* \circ \big(\tau^*, \widetilde X^*_{\tau^* \wedge \cdot} \big)^{-1}$. 
	
	\vspace{0.5em}

	Let $(\Pt^*_z)_{z \in [0,1]^{n+1}}$ be a family of r.c.p.d. of $\Pt^*$ w.r.t. the $\sigma$-field generated by $(Z_k)_{0 \le k \le n}$. 
	Then there is a $\Pt^*$-null set $N \subset [0,1]^n$ such that for each $z \in [0,1]^n \setminus N$, 
	under $\Pt^*_z$, $W^*$ is still a Brownian motion and \eqref{eq:def_X0} holds true
	(see Section 4 of Claisse, Talay and Tan \cite{JDX} for some technical subtitles).
	Notice that $\tilde \nu^*$ is adapted to the (augmented) Brownian filtration under $\Pt^*_z$,
	using Proposition \ref{prop:eqiv_opt_stopping}, it follows that for each $z \in [0,1]^n \setminus N$, one has
	\b*
		\E^{\Pt^*_z} \Big[ \Phi_2(\tau^*, \widetilde X^*_{\tau^*\wedge \cdot}) \Big] &\le& V^{S_0}_2.
	\e*
	And hence
	\b*
		\E^{\Pt^*} \Big[\Phi_2 \big(\tau^*, \widetilde X^*_{\tau^* \wedge \cdot} \big) \Big]
		&=&
		\int_{[0,1]^{n+1}} \E^{\Pt^*_z} ~ \big[\Phi_2(\tau^*, \widetilde X^*_{\tau^* \wedge \cdot}) \big] ~dz
		~~\le~~
		V^{S_0}_2.
	\e*
	This is enough to prove that $V^{W_0}_2 = W^{S_0}_2$. 
	
	\vspace{0.5em}
	
	\noindent \rmii Given the trivial inequalities
	\b*
		V^{S_0}_2 ~~\le~~ V^S_2 ~~\le~~ V^W_2 ~~\le~~ V^R_2,
		~~~\mbox{and}~~
		V^{W_0}_2 ~~\le~~ V^W_2,
	\e*
	and that $V^{S_0}_2 = V^{W_0}_2$,
	it is enough to prove that $V^{W_0}_2 = V^R_2$ to obtain their equivalence.
	In fact, the equivalence $V^{W_0}_2 = V^R_2$
	is a direct consequence of Propositions \ref{prop:approximate_q} and \ref{prop:approximation_general}, 
	together with the semicontinuity and boundedness from below of $L$, $\Phi_2$ (as in Remark \ref{rem:approximation_lsc}).
	\qed

\section{Conclusions}

	We studied a general controlled/stopped martingale problem and showed its dynamic programming principle under the abstract framework given in our previous work \cite{ElKarouiTan1}.
	In particular, to derive the DPP, we don't need uniqueness of control/stopping rules, neither the existence of the optimal control/stopping rules.
	Restricted to the controlled/stopped diffusion processes problem, 
	we obtained the dynamic programming principle for different formulations of the control/stopping problem,
	including the relaxed formulation, weak formulation, and the strong formulation, where in the last one the probability space together with the Brownian motion is fixed. 
	Moreover, under further regularity conditions, we obtained a stability result as well as the equivalence of the value function of different formulations of the control/stopping problem.

\end{document}